\documentclass[runningheads]{llncs}
\usepackage[T1]{fontenc}
\usepackage{graphicx}
\usepackage{amsmath}
\usepackage{blindtext}
\usepackage{caption}
\usepackage[pass]{geometry}
\usepackage{latexsym}
\usepackage{placeins}

\numberwithin{equation}{section}
\numberwithin{figure}{section}
\numberwithin{table}{section}
\numberwithin{theorem}{section}
\numberwithin{lemma}{section}
\numberwithin{remark}{section}

\DeclareSymbolFont{AMSb}{U}{msb}{m}{n}
\DeclareMathSymbol{\B}{\mathalpha}{AMSb}{"42}
\DeclareMathSymbol{\F}{\mathalpha}{AMSb}{"46}
\DeclareMathSymbol{\Q}{\mathalpha}{AMSb}{"51}
\DeclareMathSymbol{\R}{\mathalpha}{AMSb}{"52}
\DeclareMathSymbol{\Z}{\mathalpha}{AMSb}{"5A}

\newcommand\AHull[1]{{\mathrm{aff}(#1)}}
\newcommand\CHull[1]{{\mathrm{conv}(#1)}}
\newcommand\TSP[1]{{\mathrm{TSP}(#1)}}
\newcommand\STHGP[1]{{\mathrm{STHGP}(#1)}}
\newcommand\STGP[1]{{\mathrm{STGP}(#1)}}
\newcommand\Bell{\mathrm{Bell}}

\newcommand{\MyFrac}[2]{{\genfrac{}{}{0.5pt}{0}{#1}{#2}}}

\newcommand\PLUS{{\, + \,}}
\newcommand\MINUS{{\, - \,}}
\newcommand\TIMES{{\,}}

\spnewtheorem{MyLemma}[theorem]{Lemma}{\bfseries}{\itshape}
\spnewtheorem*{MyProof}{Proof}{\bfseries}{\upshape}

\newcommand{\MyBeginFig}{\begin{figure}[htbp]\centering}
\newcommand{\MyEndFig}{\end{figure}}
\newcommand{\MyBeginTable}{\begin{table}[htbp]\centering}
\newcommand{\MyEndTable}{\end{table}}

\newcommand{\MyAtop}[2]{{\genfrac{}{}{0pt}{0}{#1}{#2}}}

\begin{document}
\title{Quantitative Indicators for Strength of Inequalities with
  Respect to a Polyhedron}%
\subtitle{Part I: Theory}
\titlerunning{Strength of Inequalities with Respect to a Polyhedron, Part I: Theory}
\author{David M. Warme%
%\inst{1}
\orcidID{0009-0004-2307-9812}}
\authorrunning{D. M. Warme}
\institute{Group W, Vienna VA 22180, USA
\email{dwarme@groupw.com}%
%\url{http://www.springer.com/gp/computer-science/lncs}
}%
\maketitle              % typeset the header of the contribution

\begin{quotation}
\noindent
Measure what is measurable, and make measurable what is not so. \\
\hspace*{\fill}Gallileo Galilei (1564--1642)
\end{quotation}

\begin{abstract}

We study the notion of ``strength'' of inequalities used in integer
and mixed-integer programming, and the branch-and-cut algorithms used
to solve such problems.
Strength is an ethereal property lacking any good formal definition,
but crucially affects speed of computations.
We review several quantitative indicators proposed in the literature
that we claim provide a measure of the relative strength of
inequalities with respect to a given polyhedron.
We evaluate two of these indicators
(extreme point ratio (EPR) and centroid distance (CD))
on various facet classes for both
the traveling salesman polytope $\TSP{n}$, and the spanning tree in
hypergraph polytope $\STHGP{n}$, obtaining closed-form expressions for
each indicator on each facet class.
Within each facet class, the two indicators yield strikingly similar
strength rankings, with excellent agreement on which facets are
strongest and which are weakest.
Both indicators corroborate all known computational experience with
both polytopes.
The indicators also reveal properties of $\STHGP{n}$ subtours that
were previously neither known nor suspected.

We also evaluate these indicators for the subtour inequalities of the
spanning tree in graphs polytope $\STGP{n}$, obtaining surprising and
unexpected results that (at least for STGP and STHGP subtours) lead us
to believe EPR to be a more accurate estimate of strength than CD.

Applications include:
 comparing the relative strength of different classes of inequalities;
 design of rapidly-converging separation algorithms;
 and design or justification for constraint strengthening procedures.

The companion paper exploits one of the newly revealed properties of
$\STHGP{n}$ subtours in GeoSteiner, presenting detailed computational
results.
Across all distance metrics and instances studied, these results are
remarkable --- culminating with an optimal solution of a 1,000,000
terminal random Euclidean instance.
This confirms these indicators to be highly predictive and strongly
correlated with actual computational strength.

\keywords{Integer programming	\and%
  Combinatorial optimization	\and%
  Steiner tree			\and%
  Traveling salesman problem	\and%
  Enumerative Combinatorics	\and%
  Extreme point ratio		\and%
  Centroid distance}%
{%
\renewcommand{\keywordname}{{\bf Mathematics Subject Classification (2020):}}
\keywords{90C10 \and 90C11 \and 90C27 \and 05A15}
}%
\end{abstract}

\section{Introduction}

The effect studied in this paper is strength of inequalites used in
integer and mixed-integer programming.  Strong inequalities cause
computations to finish quickly, whereas weak inequalities cause
computations to progress more slowly, even to the point where
computation is abandoned.

This characterization, while descriptive, lacks both rigor and
practical utility.
We would prefer notions of strength for which the use of inequalities
deemed to be ``strong'' generally result in rapid computational
solutions.
This is quite challenging in practice,
especially in a discipline where computational performance is achieved
only through a careful blend of many factors that interact in
enormously complex ways.
For example, a set of inequalities that collectively close a large
fraction of the gap might nonetheless introduce a high condition
number into the LP basis matrix, or have a very high density of
non-zeros.
Either of these situations can obstruct the LP solver, significantly
slowing all subsequent computation.

The strength of an inequality has previously been characterized in
several ways:
\begin{enumerate}
  \item
	\label{enum:goemans}
	Best-case improvement in the objective value when the
	inequality is added to a given relaxed polyhedron.
	This was proposed by Goemans~\cite{Goemans1995}.
  \item
	\label{enum:CG-rank}
	Chv\'{a}tal-Gomory rank~\cite{Gomory1958,Chvatal1973}.
  \item
	\label{enum:shooting}
	Shooting experiment size, proposed by
        Kuhn~\cite{Kuhn1964,Kuhn1991}, with further development by
        Gomory, Johnson, Evans~\cite{GomoryJohnsonEvans2003} and
        Hunsaker~\cite{Hunsaker2003}.
  \item
	\label{enum:prob-int-optimum}
	Probability of an integer optimum when the inequality is added
        to a given relaxed polyhedron~\cite{Hunsaker2003}.
\end{enumerate}
Such measures pose difficulties:
\ref{enum:goemans} and~\ref{enum:prob-int-optimum} depend upon the
objective function.
In case~\ref{enum:goemans}, the objective function offering the
greatest improvement is assumed and says nothing about the behavior of
the inequality away from this objective.
Case~\ref{enum:prob-int-optimum} requires choosing a particular
distribution of objective functions.
Computing the Chv\'{a}tal-Gomory rank is very difficult in general,
although Hunsaker~\cite{Hunsaker2003} gives a method for classifying
an inequality as rank 0, rank 1, or rank $>1$.
(This does not provide a very fine gradation of strength.)
Shooting (case~\ref{enum:shooting}) depends upon choice of shooting
point.
All of these except~\ref{enum:shooting} depend upon a particular
relaxation.

More recent work on quantifying strength of inequalities:
Walter~\cite{Walter2021} studies general-purpose cutting planes used
by SCIP while solving practical instances, computationally assessing
the dimensionality of the induced face as a measure of strength;
Basu, et al.~\cite{BasuConfortiSummaZambelli2018} quantify strength of
valid inequalities from Gomory and Johnson's finite and infinite group
relaxations, maximizing the volume of the nonnegative orthant cut off
by a valid inequality.
Dey and Molinaro~\cite{DeyMolinaro2018} discuss issues faced by
general purpose MIP solvers, including various classes of inequalities
available, and difficulties involved in {\em selecting} small,
computationally effective subsets of inequalities from a larger pool
of violated inequalities (requiring strength assessment of individual
inequalities);
Fitisone and Zhou~\cite{FitisoneZhou2023} present a method to compute
solid angle measure of polyhedral cones using multivariate
hypergeometric series, potentially yielding closed-forms.

It is tempting to consider the objective function when
considering strength, for by doing so we can obtain stronger,
more context-sensitive notions of strength.
Taking this idea to its logical conclusion, we note that a specific
objective implies a specific optimal solution (or face) $x$,
and can then assert that all facets incident to $x$
% (or a best sufficient subset thereof --- where best means fewest
% total non-zeros or some other attractive measure)
are strong, and all others are weak.
Such a definition for strength may be breathtakingly strong and
precise, but has no practical algorithmic utility since few (if any)
of these strong inequalities might be known nor does it provide
any assistance in finding them during actual computations.

We take the exact opposite approach, completely ignoring the
objective function.
The indicators presented here depend only upon the inequality itself
and the polyhedron being optimized over, and are based solely on
combinatorial and/or geometric properties of the polyhedron.
We present several such indicators, and study two of them
(extreme point ratio (EPR) and centroid distance (CD)) in more detail,
calculating both indicators for several well-known facet classes of
both the Traveling Salesman Polytope (TSP) and the Spanning Tree in
Hypergraph Polytope (STHGP).
Despite their very different definitions, these two indicators
strongly agree with each other and corroborate all that we have
learned (through extensive computational experience with both
polytopes) about the relative strengths of these inequalities in
practice.
Although general in nature, the companion
paper~\cite{WarmeIndicators2} demonstrates that these indicators
successfully specialize to given objective functions, yielding
computationally meaningful quantifications of strength.
% in the context of that given objective.

These indicators assume that the objective is linear, so that the set
$X$ of extreme points of the polyhedron is sufficient to contain at
least one optimal solution for each possible objective function
(except those for which the optimal objective is unbounded).
Non-linear objectives can force even unique optimal solutions to
reside anywhere within the polyhedron, including the interior.
% It is an open question whether these indicators usefully generalize to
% non-linear objectives, but caution is advised.

We believe strength indicators will become a natural and standard tool
used by researchers to provide insight, answer important questions,
and lead them in successful directions --- but researchers will first
need to gain experience with and trust in the predictive power of
these indicators via a large body of supportive computational
evidence.

When designing separation algorithms or constraint-strengthening
procedures, no particular objective function should be assumed nor
should one assume particular cuts already generated (although we
recognize that in some cases correctness may depend upon certain
properties of the objective function, e.g., satisfying the triangle
inequality or residing in the non-negative orthant).

The hallmark of an {\em effective} separation algorithm is that it
{\em converge rapidly}, i.e., using few optimize/separate iterations,
and that these iterations run quickly.
There can be many orders of magnitude difference in the convergence
rates of simple/naive versus highly tuned separation algorithms.
Such tuning is often experimental --- more art than science.
If a body of evidence confirming the predictive power of these
indicators appears, we assert that strength indicators can replace
this art and experimentation with quantitative data and mathematical
insight, leading the reasearcher toward highly tuned separation
algorithms and constraint strengthening procedures in a principled
manner.
In effect, we propose an {\em empirical} approach to designing
effective separation and constraint strengthening algorithms.
% Empiricism is the hallmark of science, and the aim should be to
% establish a rigorous, principled and more methodical path toward
% rapidly-converging separation algorithms.

All constraints used in integer and mixed-integer programming must
of course be {\em valid}.
Beyond that, inequalities are either (a) {\em facet-defining} or
(b) {\em not facet-defining}.
These two classes of inequalities can be viewed as a
``practically complete'' taxonomy for two reasons:
(1) Facet-defining inequalities are the strongest possible in the sense
that none of them can be tightened without cutting off at least one
integer feasible solution;
(2) When using separation procedures, there are sound reasons for
avoiding {\em any} constraints that are not facet-defining.
Introducing such ``merely valid'' constraints into a computation
causes adjacent facet-defining inequalities to become
``nearly satisfied,'' greatly delaying their discovery by separation
procedures that typically look for maximally-violated inqualities.
Because of this, using non-facet-defining inequalities can drastically
slow the convergence of separation algorithms for inequalities that
are facet-defining.
This philosophy can be relaxed in cases where separation algorithms
for facet-defining classes of inequalities are {\em not} being used,
or when all pertinent facet classes have already been satisfied.
In such cases, any tighter approximation to the integer hull can have
computational merit.

Once a class of inequalities has been proved to be facet-defining,
many researchers cease further polyhedral investigations.
We aim to demonstrate the value of the additional analysis needed to
calculate the indicators.
% A separation algorithm may have many violated inequalities to choose
% from, all of them facet-defining --- but some may be vastly more
% effective than others at improving the dual bound.
% Even within a single facet class, strength can vary widely.
% The indicators presented here show precisely which inequalities are
% strong and which are weak.
% The indicators can provide much insight into how to find violations that
% are strong, or how to split, bend or contort weak inequalities into
% one or more vastly stronger inequalities.

This paper studies EPR and CD indicators in depth.
These are not original.
EPR is a minor adjustment to an indicator suggested by
Naddef and Rinaldi~\cite{NaddefRinaldiCrown1992}.
CD has previously been studied for the TSP by
Queyranne (personal communication).
We believe this to be the first large-scale and systematic study of
these indicators, however.

This paper makes the following contributions:
\begin{itemize}
 \item Formal definitions of EPR and CD indicators.
 \item Exact closed-form expressions for:
   \begin{itemize}
	\item EPR and CD of TSP Non-negativity and subtour
          inequalities (not original, but included for completeness).
	\item CD of TSP 3-toothed comb inequalities.
	\item EPR and CD of STHGP non-negativity and subtour inequalities.
	\item EPR and CD of STGP subtour inequalities.
   \end{itemize}
 \item Previously unknown properties of $\STHGP{n}$ subtours that can
   be exploited algorithmically.
 \item A surprising result for STGP subtours casting suspicion on CD.
 \item The greatest disagreement between EPR and CD indicators for
   $\STHGP{n}$ subtours is regarding the weakest subtours.
 \item A formula giving the angle between two arbitrary $\STHGP{n}$
   subtours.
 \item Example polytopes showing these indicators are not infallible.
 \item A call to replace trial and error approaches to tuning
   separation and constraint strengthening algorithms with a more
   disciplined, empirical approach informed by strength indicators.
\end{itemize}

In the companion paper~\cite{WarmeIndicators2}, we exploit these new
theoretical insights algorithmically in GeoSteiner to great
computational effect.
Results for the large suite of test instances (expanded over that
of~\cite{JuhlWarmeWinterZachariasen2018}) demonstrate a 10-fold
increase in the size of practically solvable instances, culminating in
the optimal solution of a random Euclidean instance having 1,000,000
terminals.
(Condensed versions of this paper and its companion have been
submitted for journal publication.)

The rest of this paper is structured as follows:
Section~\ref{sec:notation} defines notation and definitions used
throughout the paper.
Section~\ref{sec:indicators} defines all of the indicators proposed
in this paper, two of which are further explored in subsequent
sections.
Section~\ref{sec:result-summary}
presents a brief summary of our main results.
Section~\ref{sec:tsp-polytope} presents our results for TSP.
Section~\ref{sec:sthgp-polytope} presents our results for STHGP.
Section~\ref{sec:stgp-polytope} presents our results for STGP.
Section~\ref{sec:validation} presents the methods we have used to
independently validate our indicator formulae.
Section~\ref{sec:counter-examples}
shows that neither EPR nor CD are infallible, presenting example
polytopes upon which each indicator performs poorly.
Section~\ref{sec:open-problems} presents some future research and open
problems.
Section~\ref{sec:conclusion} concludes.

\section{Notation and Definitions}
\label{sec:notation}

\noindent
{\bf Basic Notation:}
We denote the integers as $\Z$ and the real numbers as $\R$.
The corresponding non-negative sets are denoted
$\Z^+$ and $\R^+$, respectively.
We define $\B = \lbrace 0, 1 \rbrace$.
For $x,\,y \in \R^n$, we define
$x \cdot y \,\equiv\, \sum_{i=1}^n x_i \, y_i$.
For $x \in \R^n$, we define
$||x|| \,\equiv\, ||x||_2 \,\equiv\, \sqrt{x.x}$.

\strut

\noindent
{\bf Hulls:}
Let $X \subset \R^n$.
We define the affine hull of $X$:
\begin{displaymath}
\AHull{X} \equiv \lbrace \sum_i a_i \, x_i :
 \forall i ((a_i \in \R) \wedge (x_i \in X))
 \wedge \sum_i a_i = 1 \rbrace.
\end{displaymath}
We define the convex hull of $X$:
\begin{displaymath}
\CHull{X} \equiv \lbrace \sum_i a_i \, x_i :
 \forall i ((a_i \in \R^+) \wedge (x_i \in X))
 \wedge \sum_i a_i = 1 \rbrace.
\end{displaymath}

\noindent
{\bf Incidence Vectors:}
Let $P$ be a finite set, and let $n=|P|$.
Let $x \in \R^n$.
For any $Q \subseteq P$,
we say that $x$ is the {\em incidence vector of $Q$}
if for all $p \in P$,
$((p \in Q \Longrightarrow x_p=1$) and
($p \notin Q \Longrightarrow x_p=0))$.
Given an incidence vector $x \in \B^n$, the
corresponding set $Q$ satisfies
$Q = \lbrace p \in P : x_p = 1 \rbrace$.
There is therefore a one-to-one correspondence between incidence
vectors and the subsets they represent.

\strut

\noindent
{\bf Graphs and Hypergraphs:}
We carefully define both graphs and hypergraphs because hypergraphs
may be less familiar and also to highlight the nuances between them,
such as paths in graphs versus chains in hypergraphs.

\strut

\noindent
{\bf Graphs:}
A graph $G=(V,E)$ consists of a finite set $V$ of {\em vertices},
and a set $E \subseteq V \times V$ of {\em edges}.
The vertices $u$ and $v$ of edge $(u,v) \in E$ are not considered
to be listed in any particular order --- edge $(u,v)$ is considered to
be the same edge as $(v,u)$.
Edges of the form $(u,u)$ are called {\em self loops}, and are
generally disallowed.
An edge $(s,t) \in E$ is said to be {\em incident} to vertex $v \in V$
if $v = s$ or $v = t$.
The {\em degree} of a vertex $v \in V$ is defined to be
$\mathrm{deg}(v) =
 |\lbrace e \in E : \hbox{$e$ is incident to $v$} \rbrace|$.
A graph is said to be {\em complete} if for every $u,v \in V$ such
that $u \ne v$, edge $(u,v) \in E$.
The complete graph $G=(V,E)$ such that $|V|=n$ is denoted $K_n$.
A graph $G'=(V',E')$ is said to be a subgraph of graph $G=(V,E)$ if
$V' \subseteq V$ and $E' \subseteq E$.
A {\em path} consists of an alternating sequence
$v_0,\, e_1,\, v_1,\, e_2,\, v_2,\, \ldots, e_{k-1},\, v_{k-1},\, e_k,\, v_k$
such that $(v_{i-1},v_i) = e_i \in E$ for all $1 \le i \le k$,
and $e_{i-1} \ne e_i$ for all $2 \le i \le k$.
This can also be called a {\em path from $v_0$ to $v_k$}.
The {\em length} of such a path is $k$.
A path of length 0 therefore consists of a single vertex $v_0$.
A path of length $k \ge 3$ such that all $v_i$ are distinct except
that $v_0 = v_k$ is called a {\em cycle}, or {\em cycle of length $k$}.
A graph is {\em connected} if for every pair $u,v \in V$ there is
a path from $u$ to $v$.
A graph that contains no cycles is said to be {\em acyclic}.
A graph is a {\em tree} if it is both connected and acyclic.
A {\em Hamiltonian cycle} of a graph $G=(V,E)$ is a cycle of length
$|V|$ containing every $v \in V$.
A Hamiltonian cycle is also called a {\em tour}.

\strut

\noindent
{\bf Hypergraphs:}
A hypergraph $H=(V,E)$ consists of a finite set $V$ of {\em vertices},
and a set $E \subset 2^V$ such that $|e| \ge 2$ for all $e \in E$.
(Edges $e$ such that $|e|=1$ are
called {\em self loops}, and are generally disallowed.  The empty edge
$e = \emptyset$ is also generally disallowed.)  If $e \in E$, then
$e \subseteq V$ is a set and we assume no particular ordering of its
members.
An edge $e \in E$ is said to be {\em incident} to vertex $v \in V$
if $v \in e$.
The {\em degree} of a vertex $v \in V$ is defined to be
$\mathrm{deg}(v) = |\lbrace e \in E : v \in e \rbrace|$.
A hypergraph $H=(V,E)$ is said to be {\em complete} if
$E = \lbrace e \subseteq V : |e| \ge 2 \rbrace$.
A hypergraph $H'=(V',E')$ is said to be a {\em subhypergraph} of
hypergraph $H=(V,E)$ if $V' \subseteq V$ and
for all $e' \in E'$ there exists $e \in E$ such that
$e' \,=\, e \cap V'$.
A {\em chain} consists of an alternating sequence
$v_0,\, e_1,\, v_1,\, e_2,\, v_2,\, \ldots, v_{k-1},\, e_k,\, v_k$
such that
$v_{i-1} \ne v_i$,
$v_{i-1} \in e_i$ and
$v_i \in e_i$
for all $1 \le i \le k$, and
$e_{i-1} \ne e_i$ for all $2 \le i \le k$.
This can also be called a {\em chain from $v_0$ to $v_k$}.
The {\em length} of such a chain is $k$.
A chain of length 0 therefore consists of a single vertex $v_0$.
A chain of length $k \ge 2$ such that all $e_i$ are distinct
and all $v_i$ are distinct except that $v_0 = v_k$ is called a
{\em cycle}, or {\em cycle of length $k$}.
A hypergraph is {\em connected} if for every pair $u,v \in V$ there is
a chain from $u$ to $v$.
A hypergraph that contains no cycles is said to be {\em acyclic}.
A hypergraph is a {\em tree} if it is both connected and acyclic.
Note that the literature contains multiple (incompatible) definitions
for {\em tree} in a hypergraph.

\strut

\noindent
{\bf Common Graph / Hypergraph Notation:}
The following notations are used for both graphs and hypergraphs.
For any $S \subseteq V$, let:
\begin{eqnarray*}
  E(S) &=& \lbrace e \in E : |e \cap S| \ge 2 \rbrace, \\
  \delta(S) &=&
	\lbrace e \in E : \emptyset \subset (e \cap S) \subset e \rbrace.
\end{eqnarray*}
We extend the latter definition to also apply to single vertices so that
for any $v \in V$,
$\delta(v) = \delta(\lbrace v \rbrace)$.
For any $x \in \R^{|E|}$,
and $F \subseteq E$,
$x(F)$ denotes $\sum_{e \in F} x_e$.
Let $S,T \subseteq V$ such that $S \cap T = \emptyset$.
We define
\begin{eqnarray*}
	(S:T) &=& \lbrace e \in E :
		(e \cap S \ne \emptyset) \wedge
		(e \cap T \ne \emptyset)
		\rbrace.
\end{eqnarray*}
Let $x \in \R^{|E|}$.
We define
\begin{eqnarray*}
	x(S)	&=& \sum_{e \in E} \max(|e \cap S| \,-\, 1,\, 0)\, x_e, \\
	x(S:T)	&=& \sum_{e \in (S:T)} x_e.
\end{eqnarray*}

\noindent
{\bf Generating Functions:}
Let $A \,=\, \lbrack a_0,\, a_1,\, \ldots,\, a_n,\, \ldots \rbrack$ be
an infinite sequence of real numbers.
The {\em generating function} for sequence $A$ is the function $f(z)$
such that
\begin{displaymath}
	f(z) ~=~ \sum_{i\ge0} a_i \, z^i.
\end{displaymath}
The {\em exponential generating function} for sequence $A$ is the
function $f(z)$ such that
\begin{displaymath}
	f(z) ~=~ \sum_{i \ge 0} a_i \, {{z^i} \over {i!}}.
\end{displaymath}
Let $f(z)$ be a series in powers of $z$.  Then
$\lbrack z^n \rbrack \,\, f(z)$
denotes the coefficient of $z^n$ in the series $f(z)$.  If $\mu$
is any non-zero real number, then
$\lbrack {{z^n} \over {\mu}} \rbrack \,\, f(z)$
denotes $\mu \,\, \lbrack z^n \rbrack \,\, f(z)$.

\strut

\noindent
{\bf Probability and Random Variables:}
Let $S$ be some finite or countably infinite set, and $X \in S$ be a
random variable over $S$.
Then for all $s \in S$ we define $P(X = s)$ to be the probability that
$X$ has value $s$.
We define the {\em expected value} of $X$ (also called the {\em mean}
of $X$) as
$E[X] \,\equiv\,\sum_{s \in S} s \, P(X=s)$.
We define the {\em $k$-th moment} of $X$ to be
$E[X^k] \,\equiv\, \sum_{s \in S} s^k \, P(X=s)$.
A random variable $X_{\lambda} \in \Z^+$ is a
{\em Poisson random variable with parameter $\lambda$} if for all
$k \in \Z^+$,
$P(X_{\lambda} \,=\, k) \,=\, {{\lambda^k \, e^{-k}} \over {k!}}$.
It is easily shown that $E[ X_{\lambda}] \,=\, \lambda$.

\strut

\noindent
{\bf Combinatorial Functions:}
We define ${n \brace k}$ to be the Stirling numbers of the second
kind, i.e., the number of distinct ways that $n$ items can be
partitioned into $k$ non-empty subsets.
They can be defined by the following
recurrence~\cite{GrahamKnuthPatashnik}:
\begin{eqnarray*}
{0 \brace 0}\,\,	&=& 1, \\
{n \brace 0}\,\,	&=& 0 \hbox{\ \ for $n \ge 1$}, \\
{n \brace k}\,\,	&=& 0 \hbox{\ \ for $n < k$}, \\
{n \brace k}\,\,	&=&
	k \,\, {n-1 \brace k} + {n-1 \brace k-1}
	\hbox{\ \ for $1 \le k \le n$}.
\end{eqnarray*}
For $k > 0$, let $\Bell(k)$ be the $k$-th Bell number ($\Bell(k)$ is
the number of ways of partitioning $k$ items into non-empty subsets).
The Bell numbers can be expressed in terms of the Stirling numbers of
the second kind:
\begin{displaymath}
\label{eq:bell-as-stirling-sum}
\Bell(k) = \sum_{i=1}^k {k \brace i}.
\end{displaymath}

\section{Strength Indicators}
\label{sec:indicators}

We seek quantitative indicators $I$ having at least two properties:
(1) Some natural, intuitive explanation why $I$ should be a good
estimate of computational strength; and
(2) Extensive computational evidence that $I$ is highly-corellated
with actual computational strength.
Only this second property is essential, although the first property
can help identify potentially useful indicators and provide insight
into why they work.

The strength of a given inequality only has meaning with respect to a
particular polyhedron over which the objective is being optimized.
Therefore, each indicator we propose is a function
$m : (H,P) \mapsto \R^+$,
where $H$ is the hyperplane of an inequality, and $P$ is a polyhedron.

\subsection{Extreme Point Ratio}
\label{sec:indicator-extreme-points}

The first indicator that we study is to simply count the number of
extreme points of polyhedron $P$ that are incident to hyperplane $H$.
This indicator was proposed by
Naddef and Rinaldi~\cite{NaddefRinaldiCrown1992},
``without any claim of generality.''
One aim of this paper is to begin exploring the extent to which their
indicator might generalize to other inequalities and polyhedra.

Let $X$ be the set of all extreme points of polyhedron $P$.
It is convenient to divide Naddef and Rinaldi's indicator by the total
number of extreme points of $P$, yielding the Extreme Point Ratio
(EPR) indicator:
\begin{displaymath}
	epr(P,H) ~=~ {{|H \cap X|} \over {|X|}}.
\end{displaymath}
We assert that larger EPR indicates stronger inequalities.
This ratio can be interpreted as the probability that inequality
defined by $H$
will be incident to the optimal solution (under the unlikely assumption
that all extreme points are equiprobable).
Furthermore, optimal solution $x$ is generally identified by finding
sufficiently many hyperplanes incident to $x$ that their intersection
is just $x$ --- hyperplanes incident to more extreme points are likely
to be more effective in this regard.

One disadvantage of this indicator is that it should probably only be
used for inequalities that are facet-defining.
For example, one could easily relax a facet by some small
$\epsilon > 0$
on the right hand side, obtaining a constraint that is of similar
actual strength, but which this indicator assigns a strength of zero.
Still, one might consider using it on inequalities that form high
dimensional faces of the polyhedron, e.g., {\em ridges} (which have
dimensionality one less than that of a facet), but it is probably
meaningless to compare ratios of hyperplanes defining faces of
different dimensionality.
(Let $0 < d_1 < d_2$ be integers.
Unless you can show, e.g., that the average number of extreme points
per face of dimension $d_1$ is roughly equal to that for faces of
dimension $d_2$ (unlikely for most polyhedra), counting extreme points
for faces of dimensions $d_1$ and $d_2$ will yield an
``apples vs oranges'' comparison.)

Note that some applications (e.g., heuristics) could benefit from
counting {\em all} feasible integer lattice points (incident to $H$
and in $P$) not just the extreme points.

% We note EPR are especially easy to calculate for various facets of
% $\TSP{n}$, and many have likely done so, including Queyranne
% (personal communication).
%%% FIXME -- find references!

\subsection{Centroid Distance}
\label{sec:indicator-centroid-distance}

Our second proposed indicator requires first computing what we shall
call the ``centroid'' $C$ of the given polyhedron $P$.
Let $X$ be the set of all extreme points of $P$, and
$P' = \CHull{X}$.
Then
\begin{equation}
\label{eq:centroid-def}
	C ~=~ {1 \over {|X|}} \sum_{x \in X} x.
\end{equation}
(We consider only extreme points $X$: there is no convenient way to
incorporate $P$'s extreme rays into $C$; and $X$ suffices to describe
all but unbounded optimal solutions.)
Unless polytope $P'$ is degenerate (i.e., a single point), then $C$ is
located in the interior of all full-dimensional $P'$, and in the
interior (relative to the affine hull) of all $P'$ that are not
full-dimensional.
CD indicators measure the distance between $C$ and $H$, although
care must be exercised for this to be well-defined.
Figure~\ref{fig:non-affine-dist} shows a polyhedron that is not
of full dimension, residing in the affine hull $A$, with a hyperplane
$H$
that can be tilted arbitrarily with respect to the affine hull,
thereby allowing the distance from $C$ to $H$ to be any arbitrarily
small $\epsilon > 0$.
Figure~\ref{fig:affine-dist} illustrates a consistent distance to
measure --- along a line residing entirely within the affine hull (if
not entirely within the polyhedron itself).
\MyBeginFig
 \begin{center}
  \begin{minipage}[t]{2.25in}
   \includegraphics[width=2.25in,clip=]{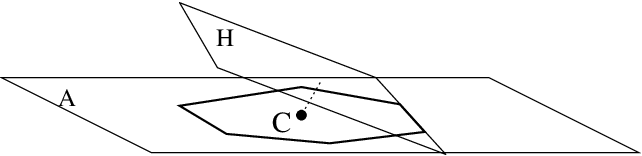}
   \caption{Distance outside the affine hull}%
   \label{fig:non-affine-dist}%
  \end{minipage}
  \quad
  \begin{minipage}[t]{2.25in}
   \includegraphics[width=2.25in,clip=]{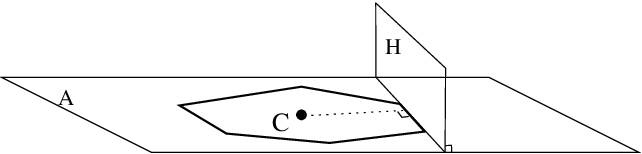}
   \caption{Distance inside the affine hull}%
   \label{fig:affine-dist}%
  \end{minipage}
 \end{center}
\MyEndFig

\noindent
We assert that smaller CD indicates stronger inequalities.
Constraints close to $C$ tend to exclude more volume than those
further away, nipping at the outer edges of the polytope.
We distinguish two versions of this indicator, the {\em normal}
indicator
$d(H,P)$,
and a {\em weak} version
$d_w(H,P)$, defined as follows:
\begin{eqnarray*}
d(H,P) &=& \min_{x \in (H \cap P)} ||C \, - \, x||_2 \\
d_w(H,P) &=& \min_{x \in (H \, \cap \, \AHull{P})} ||C - x||_2
\end{eqnarray*}
When $H \cap P = \emptyset$, we define $d(H,P) = \infty$.
When $H \cap \AHull{P} = \emptyset$, we define $d_w(H,P) = \infty$.
%
% Since $P$ is a convex body, $d(H,P)$ gives the shortest distance
% (along a line lying entirely within the polyhedron) from the centroid
% $C$ to hyperplane $H$.  Note that $d(H,P)$ is not a useful indicator
% for inequalities that do not intersect polyhedron $P$.
%
The line along which indicator $d_w(H,P)$ measures distance need not
reside entirely within polyhedron $P$.
% --- it need only reside within the affine hull of $P$.
This indicator may therefore be useful for inequalities that have no
intersection with $P$.

The weak version $d_w(H,P)$ may be easier to calculate than $d(H,P)$
in certain cases, although there are some cases where they can be
proven to be identical.
Figure~\ref{fig:weak-dist-different} illustrates a situation where
$d(H,P)$ and $d_w(H,P)$ are different.
\MyBeginFig
 \begin{minipage}[t]{4.0in}
  \begin{center}
   \includegraphics[width=1.5in,clip=]{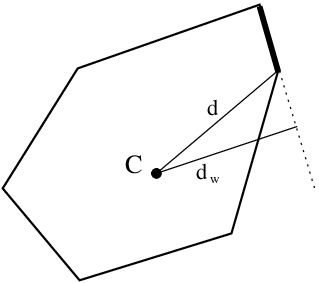}
   \caption{$d(H,P) \, \ne \, d_w(H,P)$}%
   \label{fig:weak-dist-different}%
  \end{center}
 \end{minipage}
\MyEndFig

\noindent
We note that CD indicators have previously been calculated for various
facets of $\TSP{n}$ (e.g., Queyranne, personal communication).

\subsection{Calculating Indicators}
\label{sec:calculating-indicators}

There are at least three approaches to calculating these indicators:
(1) numerically, (2) exact closed-forms, and (3) asymptotic forms.
For sufficiently small instances, it may be feasible to recursively
generate the entire set $X$ of extreme points for a polyhedron.
It is then a simple matter to:
(a) calculate the centroid $C$ using~(\ref{eq:centroid-def});
(b) generate a set $F \subseteq X$ such that $x \in F$
if-and-only-if $x$ is incident to hyperplane $H$ by testing the
members of $X$ exhaustively.
EPR is simply $|F|/|X|$.
The square of CD (both $d$ and $d_w$) can be
computed using the following quadratic program:
\begin{eqnarray}
\label{eq:qp-objective}
\hbox{Minimize:} && z ~=~ y \cdot y \\
\nonumber
\hbox{Subject to:} \\
\label{eq:y=x-centroid}
&& y = x \,-\, C, \\
\label{eq:y-in-hull}
&& x = \sum_{f \in F} a_f f, \\
\label{eq:affine-constraint}
&& \sum_{f \in F} a_f ~=~ 1, \\
\label{eq:convex-bounds}
&& a_f \ge 0~~~~\hbox{for all $f \in F$}
\end{eqnarray}%
where $C$ is the centroid,
and $f \in F$ are the extreme points on $H$.
Bounds~(\ref{eq:convex-bounds}) are optional, yielding normal CD.
Excluding them yields weak CD.
% if included, then $x \in \CHull{F}$ (strong CD);
% if omitted, then $x \in \AHull{F}$ (weak CD).
% Thus, by changing only these bounds, the same formulation can compute
% both the normal CD (convex hull) and the weak CD (affine hull).

Closed-forms have numerous advantages:
they can often be evaluated numerically for large $n$;
they may provide additional insight into structural aspects of the
problem, or establish connections to other mathematical areas or
techniques.
In contrast, Naddef and Rinaldi~\cite{NaddefRinaldiCrown1992}
recursively enumerated all extreme points, counting those incident to
the hyperplane --- a technique limited to small $n$.

Calculating exact closed-forms for these indicators in general
requires use of various mathematical techniques, (e.g.,
enumerative combinatorics) to calculate $|X|$, $|F|$ and $C$.
This paper demonstrates this is feasible for each polytope and facet
class studied herein, and likely feasible for many others as well.
(Complexity class $\#P$, however, argues that closed-forms are likely
unachievable for many polyhedra and/or inequalities.)

An asymptotic form $f(n)$ approximates some function $g(n)$
``in the limit'' as $n \to \infty$.
($f(n) \sim g(n)$ if-and-only-if
$\lim_{n \to \infty} f(n)/g(n) \,=\, 1$.)
Such forms may be useful in cases that resist exact closed-forms.
We do not further consider asymptotic forms.

\subsection{Other Indicators}
\label{sec:indicator-others}

Other indicators are certainly possible.
Suppose $F = H \cap P$ defines a facet $F$.
One could in principle measure the solid angle subtended by $F$ from
the centroid $C$.
(It would in most cases be useful to divide this by the solid angle of
the entire $n$-sphere, yielding a constant reasonably interpreted as a
probability.)
Similarly, one could in principle measure the surface area of $F$ and
divide by the entire surface area of $P$.
Such computations are in general intractible, but could perhaps be
feasible for certain polyhedra of low dimension or very special
structure.

Lee and Morris~\cite{LeeMorris1992} present an approach for comparing
the relative volumes of two polytopes as the difference in radius of
two spheres having equivalent volumes.

Hunsaker~\cite{Hunsaker2003} surveys these measures (and others),
showing that shooting experiments are equivalent to optimizing over
the polar polytope, and presenting new results for master cyclic group
polyhedra, and knapsack, matching and node-packing polytopes.

\section{Summary of Results}
\label{sec:result-summary}

\begin{figure}[!ht]
\begin{center}
\begin{minipage}[t]{2.25in}
\includegraphics[width=2.25in,clip=]{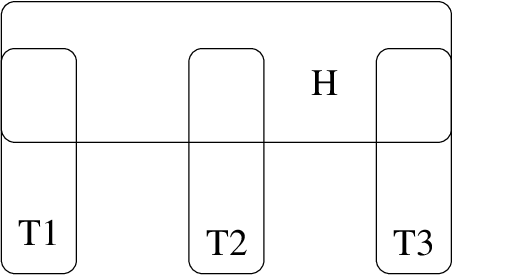}
 \captionof{figure}{Standard TSP 3-toothed comb}%
 \label{fig:standard-TSP-3comb}
\end{minipage}
\quad
\begin{minipage}[t]{2.25in}
\includegraphics[width=2.25in,clip=]{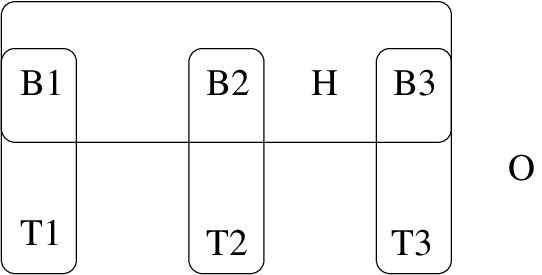}
 \captionof{figure}{Our partition of TSP 3-toothed comb into 8 disjoint sets}%
 \label{fig:our-TSP-3comb}
\end{minipage}
\end{center}
\end{figure}

In this section we briefly present our main results.
For the Traveling Salesman polytope, we study the non-negativity,
subtour and 3-toothed comb inequalities.
Figure~\ref{fig:standard-TSP-3comb} shows a TSP 3-toothed comb
inequality using the standard notation for the handle $H$, and teeth
$T_i$ wherein each tooth $T_i$ must have at least one vertex inside
$H$ and at least one vertex outside of $H$.
Figure~\ref{fig:our-TSP-3comb} illustrates a different notation that
is more convenient for our purposes.
We partition all of the vertices (both inside and outside the comb)
into 8 mutually disjoint sets
$B_1,\, T_1,\, B_2,\, T_2,\, B_3,\, T_3,\, H,$ and $O$.
This notation splits each tooth $i$ into a
``base'' $B_i$ and a
``tip'' $T_i$, where the base portion resides inside the handle
and the tip portion resides outside the handle of the comb.
We also define $H$ to represent only those vertices of the handle that
do {\em not} reside in any of the $B_i$.
Finally, we define set $O$ to be all vertices residing outside of the
comb.
We define $b_i = |B_i|$, $t_i = |T_i|$ for $1 \le i \le 3$,
$h=|H|$ and $o=|O|$,
and note that
$n \,=\,
 (b_1 \,+\, t_1) \,+\, (b_2 \,+\, t_2) \,+\, (b_3 \,+\, t_3)
 \,+\, h \,+\, o$.
Valid 3-toothed comb inequalities require that
$b_i,\, t_i \,>\, 0$ for $1 \le i \le 3$,
and $h,\, o \,\ge\, 0$.

\noindent
Table~\ref{tab:tsp-extreme-point-ratios} presents EPR and
Tables~\ref{tab:tsp-centroid-distances}
and~\ref{tab:3comb-dist-subexprs}
presents CD (squared) for these three families of facet-defining TSP
inequalities.
No correct closed-form is presently known for EPR of TSP
3-toothed comb inequalities.

\begin{table}[!ht]
\begin{center}
\renewcommand{\arraystretch}{2.0} % More space between lines of this table!
\captionof{table}{EPR for TSP Facets}%
\begin{tabular}{|l|c|} \hline
TSP Facet Class & EPR \\ \hline
Non-negativity & ${{n-3} \over {n-1}}$ \\
Subtour $k=|S|$ & ${n \over {n \choose k}}$ \\
3-toothed Combs &
~~~~~~~~ Open ~~~~~~~~ \\ \hline
\end{tabular}
\label{tab:tsp-extreme-point-ratios}
\end{center}
\begin{center}
\renewcommand{\arraystretch}{2.0} % More space between lines of this table!
\captionof{table}{CD for TSP Facets}%
\begin{tabular}{|l|c|l|} \hline
TSP Facet Class & CD Squared $d_w^2$ & $d = d_w$? \\ \hline
Non-negativity & ${4 \over {(n-1) \, (n-3)}}$ & Yes \\
Subtour $k=|S|$ & ${{2 \, (k-1) \, (n-2) \, (n-k-1)}
		\over {k \, (n-1) \, (n-k)}}$ & Yes \\
3-toothed Combs &
	${{2 \, (n-2) \, A^2} \over {(n-1) \, B}}$
 & No \\
& (See Table~\ref{tab:3comb-dist-subexprs} for $A$ and $B$.) & \\ \hline
\end{tabular}
\label{tab:tsp-centroid-distances}
\end{center}
\end{table}
\begin{table}[!ht]
\begin{center}
\captionof{table}{$A$ and $B$ for TSP 3-toothed comb CD
    (see text for special meaning of subscripts appearing in
    these expressions)}%
\begin{tabular}{|p{0.975\linewidth}|}\hline
\begin{eqnarray*}
A &=&
 (b_i \TIMES{} t_i)
		\PLUS{} 2 \TIMES{} b_i \TIMES{} b_j
		\PLUS{} 3 \TIMES{} b_i \TIMES{} t_j
		\PLUS{} 2 \TIMES{} t_i \TIMES{} t_j
		\PLUS{} b_i \TIMES{} h
		\PLUS{} 2 \TIMES{} t_i \TIMES{} h
		\PLUS{} 2 \TIMES{} b_i \TIMES{} o
		\PLUS{} t_i \TIMES{} o
		\PLUS{} h \TIMES{} o
		\MINUS{} 5 \TIMES{} (n-1), \\
&& \\
B &=&
	(b_i^2 \TIMES{} t_i^2)
	\PLUS{} 2 \TIMES{} (b_i^2 \TIMES{} t_i) \TIMES{} t_j
	\PLUS{} 4 \TIMES{} b_i^2 \TIMES{} b_j^2
	\PLUS{} 2 \TIMES{} b_i^2 \TIMES{} b_j \TIMES{} b_k
	\PLUS{} 2 \TIMES{} b_i^2 \TIMES{} b_j \TIMES{} t_k
	\PLUS{} 9 \TIMES{} b_i^2 \TIMES{} t_j^2
	\PLUS{} 8 \TIMES{} b_i^2 \TIMES{} t_j \TIMES{} t_k \\
&&
	\PLUS{} 2 \TIMES{} (b_i \TIMES{} t_i^2) \TIMES{} b_j
	\PLUS{} 4 \TIMES{} t_i^2 \TIMES{} t_j^2
	\PLUS{} 8 \TIMES{} t_i^2 \TIMES{} b_j \TIMES{} b_k
	\PLUS{} 2 \TIMES{} t_i^2 \TIMES{} b_j \TIMES{} t_k
	\PLUS{} 2 \TIMES{} t_i^2 \TIMES{} t_j \TIMES{} t_k
	\PLUS{} 12 \TIMES{} (b_i \TIMES{} t_i) \TIMES{} (b_j \TIMES{} t_j) \\
&&
	\PLUS{} 8 \TIMES{} (b_i \TIMES{} t_i) \TIMES{} b_j^2
	\PLUS{} 4 \TIMES{} (b_i \TIMES{} t_i) \TIMES{} b_j \TIMES{} b_k
	\PLUS{} 4 \TIMES{} (b_i \TIMES{} t_i) \TIMES{} b_j \TIMES{} t_k
	\PLUS{} 4 \TIMES{} (b_i \TIMES{} t_i) \TIMES{} t_j \TIMES{} t_k
	\PLUS{} 8 \TIMES{} (b_i \TIMES{} t_i) \TIMES{} t_j^2
	\PLUS{} b_i^2 \TIMES{} h^2 \\
&&
	\PLUS{} 4 \TIMES{} t_i^2 \TIMES{} h^2
	\PLUS{} 2 \TIMES{} (b_i \TIMES{} t_i) \TIMES{} h^2
	\PLUS{} 2 \TIMES{} t_i \TIMES{} t_j \TIMES{} h^2
	\PLUS{} 2 \TIMES{} b_i^2 \TIMES{} b_j \TIMES{} h
	\PLUS{} 2 \TIMES{} b_i^2 \TIMES{} t_j \TIMES{} h
	\PLUS{} 4 \TIMES{} b_i \TIMES{} t_j \TIMES{} t_k \TIMES{} h
	\PLUS{} 2 \TIMES{} (b_i \TIMES{} t_i^2) \TIMES{} h \\
&&
	\PLUS{} 4 \TIMES{} (b_i \TIMES{} t_i) \TIMES{} b_j \TIMES{} h
	\PLUS{} 4 \TIMES{} (b_i \TIMES{} t_i) \TIMES{} t_j \TIMES{} h
	\PLUS{} 8 \TIMES{} t_i^2 \TIMES{} b_j \TIMES{} h
	\PLUS{} 2 \TIMES{} t_i^2 \TIMES{} t_j \TIMES{} h
	\PLUS{} 2 \TIMES{} (b_i \TIMES{} t_i) \TIMES{} o^2 \\
&&
	\PLUS{} 4 \TIMES{} b_i^2 \TIMES{} o^2
	\PLUS{} 2 \TIMES{} b_i \TIMES{} b_j \TIMES{} o^2
	\PLUS{} t_i^2 \TIMES{} o^2
        \PLUS{} 2 \TIMES{} b_i \TIMES{} h \TIMES{} o^2
	\PLUS{} h^2 \TIMES{} o^2
	\PLUS{} 2 \TIMES{} (b_i^2 \TIMES{} t_i) \TIMES{} o
	\PLUS{} 2 \TIMES{} b_i^2 \TIMES{} b_j \TIMES{} o \\
&&
	\PLUS{} 8 \TIMES{} b_i^2 \TIMES{} t_j \TIMES{} o
	\PLUS{} 2 \TIMES{} t_i^2 \TIMES{} b_j \TIMES{} o
	\PLUS{} 2 \TIMES{} t_i^2 \TIMES{} t_j \TIMES{} o
	\PLUS{} 4 \TIMES{} b_i \TIMES{} b_j \TIMES{} t_k \TIMES{} o
	\PLUS{} 4 \TIMES{} (b_i \TIMES{} t_i) \TIMES{} b_j \TIMES{} o
	\PLUS{} 4 \TIMES{} (b_i \TIMES{} t_i) \TIMES{} t_j \TIMES{} o \\
&&
	\PLUS{} 4 \TIMES{} (b_i \TIMES{} t_i) \TIMES{} h \TIMES{} o
	\PLUS{} 2 \TIMES{} b_i^2 \TIMES{} h \TIMES{} o
	\PLUS{} 4 \TIMES{} b_i \TIMES{} t_j \TIMES{} h \TIMES{} o
	\PLUS{} 2 \TIMES{} t_i^2 \TIMES{} h \TIMES{} o
	\PLUS{} 2 \TIMES{} t_i \TIMES{} h^2 \TIMES{} o
	\MINUS{} (b_i^2 \TIMES{} t_i)
	\MINUS{} (b_i \TIMES{} t_i^2) \\
&&
	\MINUS{} 4 \TIMES{} b_i^2 \TIMES{} b_j
	\MINUS{} 9 \TIMES{} b_i^2 \TIMES{} t_j
	\MINUS{} 9 \TIMES{} t_i^2 \TIMES{} b_j
	\MINUS{} 4 \TIMES{} t_i^2 \TIMES{} t_j
	\MINUS{} 6 \TIMES{} b_i \TIMES{} b_j \TIMES{} b_k
	\MINUS{} 10 \TIMES{} (b_i \TIMES{} t_i) \TIMES{} b_j
	\MINUS{} 10 \TIMES{} (b_i \TIMES{} t_i) \TIMES{} t_j \\
&&
	\MINUS{} 12 \TIMES{} b_i \TIMES{} b_j \TIMES{} t_k
	\MINUS{} 12 \TIMES{} b_i \TIMES{} t_j \TIMES{} t_k
	\MINUS{} 6 \TIMES{} t_i \TIMES{} t_j \TIMES{} t_k
	\MINUS{} b_i \TIMES{} h^2
	\MINUS{} 4 \TIMES{} t_i \TIMES{} h^2
	\MINUS{} h^2 \TIMES{} o
	\MINUS{} b_i^2 \TIMES{} h
	\MINUS{} 4 \TIMES{} b_i \TIMES{} b_j \TIMES{} h \\
&&
	\MINUS{} 4 \TIMES{} t_i^2 \TIMES{} h
	\MINUS{} 10 \TIMES{} b_i \TIMES{} t_j \TIMES{} h
	\MINUS{} 4 \TIMES{} (b_i \TIMES{} t_i) \TIMES{} h
	\MINUS{} 6 \TIMES{} t_i \TIMES{} t_j \TIMES{} h
	\MINUS{} 4 \TIMES{} b_i \TIMES{} o^2
	\MINUS{} t_i \TIMES{} o^2
	\MINUS{} h \TIMES{} o^2
	\MINUS{} 4 \TIMES{} b_i^2 \TIMES{} o
	\MINUS{} t_i^2 \TIMES{} o \\
&&
	\MINUS{} 4 \TIMES{} (b_i \TIMES{} t_i) \TIMES{} o
	\MINUS{} 6 \TIMES{} b_i \TIMES{} b_j \TIMES{} o
	\MINUS{} 10 \TIMES{} b_i \TIMES{} t_j \TIMES{} o
	\MINUS{} 4 \TIMES{} t_i \TIMES{} t_j \TIMES{} o
	\MINUS{} 4 \TIMES{} b_i \TIMES{} h \TIMES{} o
	\MINUS{} 4 \TIMES{} t_i \TIMES{} h \TIMES{} o
	\PLUS{} 4 \TIMES{} b_i \TIMES{} b_j
	\PLUS{} (b_i \TIMES{} t_i) \\
&&
	\PLUS{} 9 \TIMES{} b_i \TIMES{} t_j
	\PLUS{}4 \TIMES{} t_i \TIMES{} t_j
	\PLUS{} b_i \TIMES{} h
	\PLUS{} 4 \TIMES{} t_i \TIMES{} h
	\PLUS{} 4 \TIMES{} b_i \TIMES{} o
	\PLUS{} t_i \TIMES{} o
        \PLUS{} h \TIMES{} o
\end{eqnarray*}
\label{tab:3comb-dist-subexprs} \\
\hline
\end{tabular}
\end{center}
\end{table}

Note that in Table~\ref{tab:3comb-dist-subexprs} we use some special
conventions within expressions
$A$ and $B$ appearing in CD for TSP 3-toothed combs.
First, we use a special subscript convention to reduce the number of
terms.  The subscripts $i$, $j$, and $k$ run over all teeth 1, 2 and
3, with the restriction that $i \ne j$, $i \ne k$ and $j \ne k$.
We furthermore restrict the subscript assignment so that no identical
term is generated more than once.
For example,
$b_i \TIMES{} t_j$
expands to 6 terms
$	b_1 \TIMES{} t_2
	\PLUS{} b_1 \TIMES{} t_3
	\PLUS{} b_2 \TIMES{} t_1
	\PLUS{} b_2 \TIMES{} t_3
	\PLUS{} b_3 \TIMES{} t_1
	\PLUS{} b_3 \TIMES{} t_2$,
whereas
$b_i \TIMES{} b_j$
expands to only 3 terms
$b_1 \TIMES{} b_2 \PLUS{} b_1 \TIMES{} b_3 \PLUS{} b_2 \TIMES{} b_3$.
Secondly, we emphasize sub-expressions such as
$(b_i \TIMES{} t_i^2)$
and
$(b_j \TIMES {} t_j)$
that refer to whole tooth structures by placing parentheses around
them.
If you ignore the $-5 \TIMES{} (n-1)$ term, then $A$ has 28
terms when fully expanded.
If you substitute
$n \,=\, (b_1+t_1)+(b_2+t_2)+(b_3+t_3)+h+o$,
then $A$ has 37 terms when fully expanded.
$B$ has 348 terms when fully expanded.

For the Spanning Tree in Hypergraph Polytope, we study the
non-negativity and subtour inequalities.
Table~\ref{tab:sthgp-extreme-point-ratios} presents EPR and
Table~\ref{tab:sthgp-centroid-distances} presents CD
for both families of facet-defining inequalities.
Table~\ref{tab:sthgp-subtour-angle} presents the closed-form giving
the interior angle between two arbitrary STHGP subtours
$S_1$ and $S_2$.
Table~\ref{tab:stgp-subtour-indicators} presents EPR
and CD (squared) for the subtours of the
Spanning Tree in Graphs Polytope.
\begin{table}[!ht]
\begin{center}
\renewcommand{\arraystretch}{2.0} % More space between lines of this table!
\captionof{table}{EPR for STHGP Facets}%
\begin{tabular}{|l|c|} \hline
$\MyAtop{\STHGP{n}}{\mathrm{Facet\phantom{j}}}$ & EPR \\ \hline
$\MyAtop{x_e \ge 0,}{~~k=|e|}$ &
 \begin{minipage}[c]{3.2in}
  \begin{displaymath}
	1 - k \, {{E[X_n^{n-k}]} \over {E[X_n^{n-1}]}}
  \end{displaymath}
 \end{minipage} \\
$\MyAtop{\mathrm{Subtour,}}{~~k=|S|}$ &
\begin{minipage}[c]{3.2in}
 \begin{displaymath}
     \MyAtop{{n \over {k \, (n-k)\,E[X_n^{n-1}]}}
	     \sum_{i=0}^{k-1} {{k-1} \brace i} \, k^i
		~~~~~~~~~~~~~~~~~~~~~~~~~~~}%
  	    {~~~
	     \left \lbrack
		\sum_{j=1}^{n-k} {{n-k} \choose j}
			 \, j
			 \, E[X_{n-k}^{n-k-j}]
			 \sum_{p=0}^j {j \choose p} E[X_k^p] \, i^{j-p}
	     \right \rbrack}
 \end{displaymath}
\end{minipage} \\
&	\hbox{\small where $X_n =$ Poisson random variable with mean $n$.} \\
\hline
\end{tabular}
\label{tab:sthgp-extreme-point-ratios}
\end{center}
\begin{center}
\renewcommand{\arraystretch}{2.0} % More space between lines of this table!
\captionof{table}{CD for STHGP Facets}%
{
\begin{tabular}{|l|c|} \hline
$\MyAtop{\STHGP{n}}{\mathrm{Facet\phantom{j}}}$ & Weak CD $(d_w)$ \\ \hline
$\MyAtop{x_e \ge 0,}{~~k=|e|}$
&
\begin{minipage}[c]{3.3in}
 \begin{displaymath}
   {{k \, E[X_n^{n-k}]} \over {E[X_n^{n-1}]}}
   \, \sqrt{{\alpha(n)}
	    \over
	    {\alpha(n) \,-\, (k-1)^2}}
 \end{displaymath}
\end{minipage} \\
$\MyAtop{\mathrm{Subtour,}}{~~k=|S|}$
&
\begin{minipage}[t]{3.3in}
 \begin{displaymath}
   {{(n-k) \, t(n,k)} \over {E[X_n^{n-1}]}} \,
	\sqrt{{\alpha(n) \over {\alpha(n) \, \gamma(n,k) \,-\, \beta^2(n,k)}}}
 \end{displaymath}
\hbox to 3.0in{where\hfill}
\vspace{-0.125in}
 \begin{eqnarray*}
  %\small
X_n &=& \hbox{Poisson random variable with mean $n$}, \\
t(n,k) &=& E[X_n(X_n+1)^{n-2}] - E[X_n^k(X_n+1)^{n-k-1}], \\
\alpha(n) &=& (n^2 \,-\, 3n \,+\, 4) \, 2^{n-2} \,-\, 1, \\
\beta(n,k) &=& \gamma(n,k) \,+\, b(n \,-\, k) \, d(k), \\
\gamma(n,k) &=& 2^{n-k} \, \alpha(k), \\
b(n) &=& n \, 2^{n-1}, \\
d(n) &=& 1 \,+\, n\,2^{n-1} \,-\, 2^n.
 \end{eqnarray*}%
\end{minipage} \\
\hline
\end{tabular}
\label{tab:sthgp-centroid-distances}
}
\end{center}
\end{table}
\begin{table}[!ht]
\begin{center}
 \captionof{table}{Interior angle between $\STHGP{n}$ subtours $S_1$ and $S_2$}%
 \begin{tabular}{|c|} \hline
  \begin{minipage}[c]{4.5in}%
   \begin{displaymath}
	\theta(n,p,q,r) ~=~ \pi \,-\, \mathrm{acos}(f(n,p,q,r))
   \end{displaymath}
  \hbox to 4.0in{where\hfill} \\
   \begin{eqnarray*}
	n &=& |V|, \\
	p &=& |S_1 \setminus S_2|, \\
	q &=& |S_2 \setminus S_1|, \\
	r &=& |S_1 \cap S_2|, \\
	f(n,p,q,r) &=&
	{{\alpha(n) \,
	  w(n,p,q,r)
	 - \beta(n,p+r) \beta(n,q+r)
	 }
	 \over
	 {\sqrt{\mu(n,p+r) \, \mu(n,q+r)}}}, \\
	w(n,p,q,r) &=&
		2^{n-p-q-r} \,
		\left \lbrack
			2^{p+q}\alpha(r)
			+ d(p) d(q)
			+ b(p) b(q) c(r)
			+ b(p+q) d(r)
		\right \rbrack, \\
	\alpha(n) &=& (n^2 \,-\, 3n \,+\, 4) \, 2^{n-2} \,-\, 1, \\
	\beta(n,k) &=& \gamma(n,k) \,+\, b(n \,-\, k) \, d(k), \\
	\gamma(n,k) &=& 2^{n-k} \, \alpha(k), \\
	\mu(n,k) &=& \alpha(n) \, \gamma(n,k) \,-\, \beta^2(n,k), \\
	b(n) &=& n \, 2^{n-1}, \\
	c(n) &=& 2^n \,-\ 1, \\
	d(n) &=& 1 \,+\, n\,2^{n-1} \,-\, 2^n.
   \end{eqnarray*}
  \end{minipage} \\
 \hline
 \end{tabular}
  \label{tab:sthgp-subtour-angle}%
\end{center}
\end{table}
\begin{table}[!ht]
\begin{center}
 \renewcommand{\arraystretch}{2.0} % More space between lines of this table!
 \captionof{table}{Indicators for $k$-vertex subtours of $\STGP{n}$}%
 \begin{tabular}{|ccccc|}\hline
  \hbox to 0.25in{} &
  EPR &
   \hbox to0.5in{} &
  $(k/n)^{k-1}$ &
   \hbox to0.25in{} \\
   \hbox to0.25in{} &
  CD Squared &
   \hbox to0.5in{} &
  ${{2 \, (k \,-\, 1)\, (n \,-\, 1)\, (n \,-\, k)} \over
    {k \, n \, (n \,+\, k \,-\, 1)}}$ &
   \hbox to0.25in{} \\
 \hline
 \end{tabular}
 \label{tab:stgp-subtour-indicators}%
\end{center}
\end{table}

To gain more insight into CDs for 3-toothed combs, we examine the
smallest such combs, where
$b_i \,=\, t_i \,=\, 1$, we also substitute $o = n - h - 6$ to obtain
a distance squared of:
\begin{displaymath}
%\small
{{2 (n-2) \left \lbrack (h + 4) n
			 - (h^2 + 6h + 16)
	  \right \rbrack^2}
 \over
 {(n-1)
  \left \lbrack f_2(h) n^2 - f_1(h) n + f_0(h)
  \right \rbrack}}
\end{displaymath}
where
\begin{eqnarray*}
f_2(h) &=& h^2 + 5h + 12 \\
f_1(h) &=& 2h^3 + 17h^2 + 59h + 96 \\
f_0(h) &=& h^4 + 12h^3 + 65h^2 + 174h + 228
\end{eqnarray*}
The limit of this as $n \rightarrow \infty$ is
\begin{displaymath}
	{{2 \, h^2 \,+\, 16\,h \,+\, 32}
	 \over
	 {h^2 \,+\, 5\,h \,+\, 12}}
\end{displaymath}
which  gives distance squared values of
8/3, 25/9, 36/13, 49/18, and 8/3, for $h=0,1,2,3,4$.  It thereafter
decreases to 2 as $h \rightarrow \infty$.  So apparently the comb
weakens abruptly between $h=0$ and $h=1$, has the same strength for
$h=0$ and $h=4$, and gradually gets stronger as the handle gets
bigger (at least according to CD).

If we examine the ratio of distances squared of the
$b_i=t_i=1$
comb to that where we instead let
$b_i=t_i=2$ and $o=n-h-12$,
and take the limit of this ratio as
$n \rightarrow \infty$,
we obtain
\begin{displaymath}
	{{(h \,+\, 4)^2 \, (h^2 \,+\, 11\,h \,+\, 78)}
	 \over
	 {(h \,+\, 13)^2 \, (h^2 \,+\, 5\,h \,+\, 12)}}
\end{displaymath}
which is always in the interval
$\lbrack 8/13, 1)$,
but has some interesting behavior for
$0 \le h \le 6$
before increasing asmptotically to 1 as
$h \rightarrow \infty$.
Because this ratio is less than 1, however, we conclude that
$b_i=t_i=2$
teeth are weaker than the smaller
$b_i=t_i=1$
teeth.
One could of course perform more nuanced analyses comparing the
relative strength of less symmetric configurations of teeth.

\begin{figure}[!ht]
\begin{center}
\begin{minipage}[t]{2.25in}
 \includegraphics[width=2.25in,clip=]{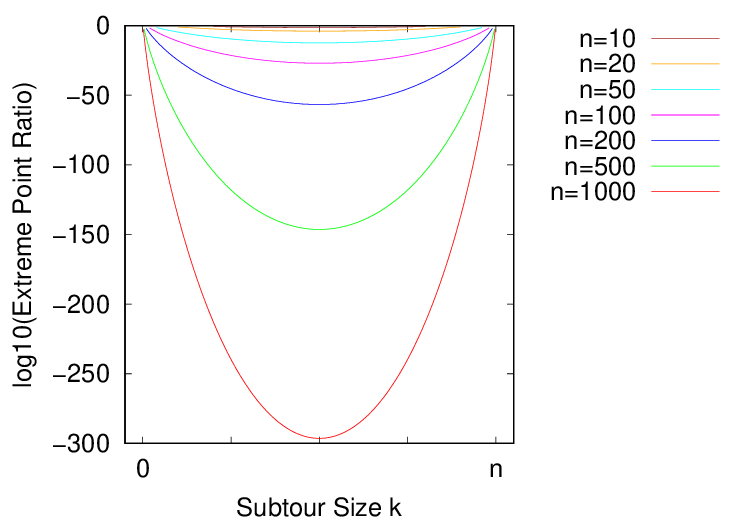}
 \captionof{figure}{Base 10 log of TSP subtour EPR (larger
     ratio = stronger)}%
 \label{fig:tsp-sec-extreme-point-ratio-plots-all}
\end{minipage}
   ~
\begin{minipage}[t]{2.25in}
 \includegraphics[width=2.25in,clip=]{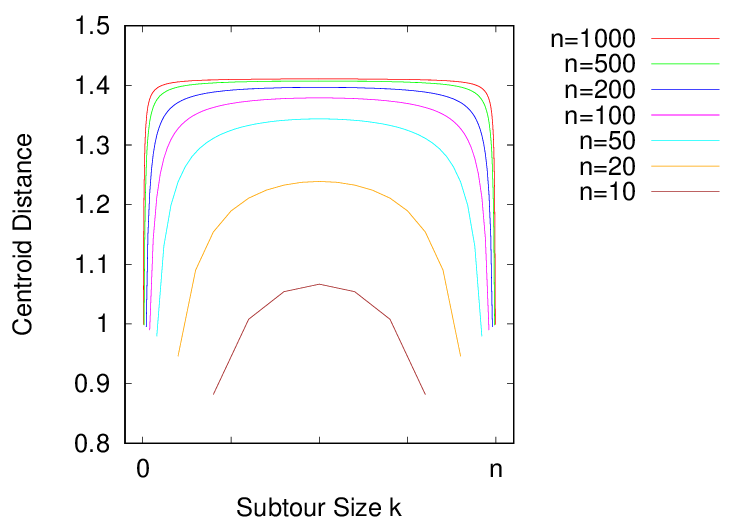}
 \captionof{figure}{TSP subtour CD $d(n,k)$ (smaller distance
     = stronger)}%
 \label{fig:tsp-sec-dist-plots-all}
\end{minipage}
\end{center}
\end{figure}

For TSP subtours,
Figure~\ref{fig:tsp-sec-extreme-point-ratio-plots-all}
plots $\log(\mathrm{EPR})$,
and Figure~\ref{fig:tsp-sec-dist-plots-all}
plots CD
as functions of subtour size $k$, for several $n$, $10 \le n \le
1000$.
Predicted strength is highly correlated between the two indicators,
which show the smallest and largest subtours are (symmetrically)
strongest.
For large $n$, only a 
vanishiningly small fraction of the subtours have any significant
strength,
with a vast region of very weak and ineffective inequalities
lying in between.
The largest subtour CD at $d(n,n/2)$
approaches $\sqrt{2}$ as $n \rightarrow \infty$.

\begin{figure}[!ht]
\begin{center}
\hbox{
\begin{minipage}[t]{2.25in}
 \includegraphics[width=2.25in,clip=]{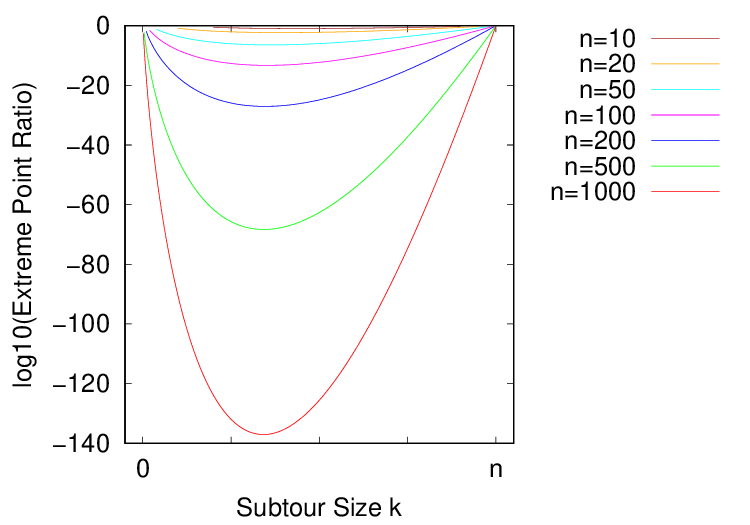}
 \captionof{figure}{Base 10 log of STHGP subtour EPR (larger
     ratio = stronger)}%
 \label{fig:sthgp-sec-extreme-point-ratio-plots-all}
\end{minipage}
   ~
\begin{minipage}[t]{2.25in}
 \includegraphics[width=2.25in,clip=]{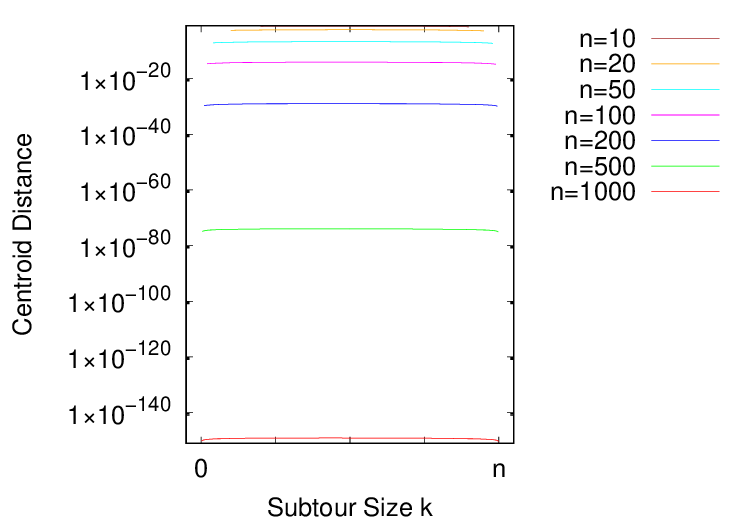}
 \captionof{figure}{STHGP subtour CD $d_w(n,k)$ (smaller distance
     = stronger)}
 \label{fig:sthgp-sec-dist-plots-all}
\end{minipage}
}
\end{center}
\end{figure}

Figure~\ref{fig:sthgp-sec-extreme-point-ratio-plots-all}
presents a plot of EPRs of STHGP subtour
inequalities, plotting the ratio as a function of $k$ (the subtour
size) for several values of $n$ from $n=10$ to $n=1000$.
(Note that we plot the base ten logarithm of EPRs because they get
so small.)

The corresponding plot for STHGP subtour CDs appears in
Figures~\ref{fig:sthgp-sec-dist-plots-all},
which is plotted with a log scale for the $y$ coordinate.
These clearly show the distances rapidly decreasing with $n$, but
cause the details of each curve to be blurred to nearly flat lines.

\begin{figure}[!ht]
\begin{center}
 \begin{minipage}[t]{2.25in}
  \includegraphics[width=2.25in,clip=]{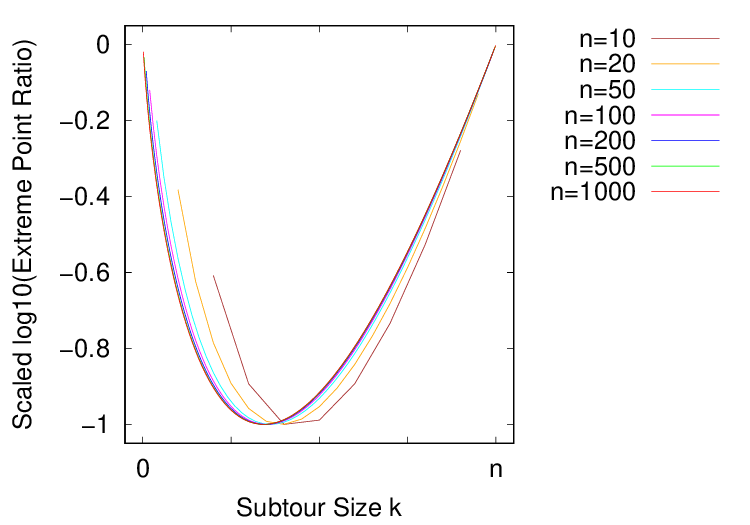}
  \captionof{figure}{Scaled STHGP subtour log10(EPR) (larger ratio = stronger)}%
  \label{fig:sthgp-sec-ratio-scaled-plot-all}
 \end{minipage}
 ~
 \begin{minipage}[t]{2.25in}
  \includegraphics[width=2.25in,clip=]{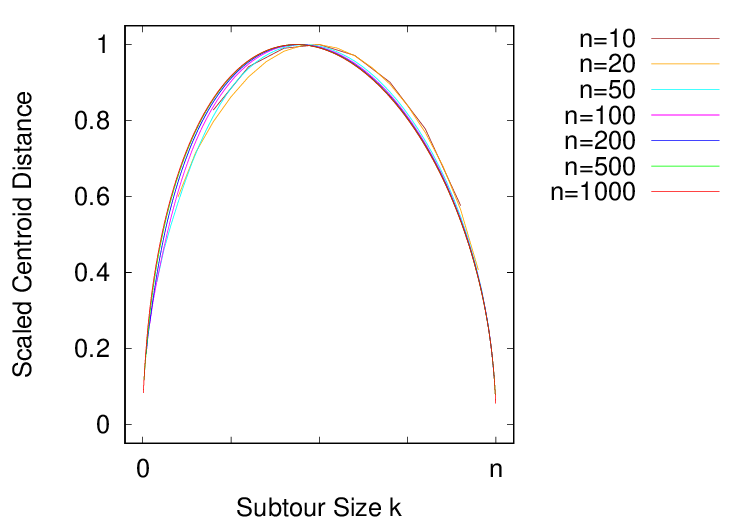}
  \captionof{figure}{Scaled STHGP subtour CD (smaller distance = stronger)}%
  \label{fig:sthgp-sec-dist-scaled-plot-all}
 \end{minipage}
\end{center}
\end{figure}

To remedy this, we plot EPR of STHGP subtours in
Figure~\ref{fig:sthgp-sec-ratio-scaled-plot-all}
with the $y$ coordinates of each curve scaled to have minimum
log(ratio) of -1.
The corresponding plot of STHGP subtour CDs with the $y$ coordinates
of each curve scaled to have maximum distance of 1 appears in
Figure~\ref{fig:sthgp-sec-dist-scaled-plot-all}.
Note how the scaled curves retain their shape as $n$ increases.

Predicted strength is highly correlated between the two indicators,
which clearly show that the very smallest and largest
subtours are the strongest.
Note that the very strongest subtour is actually for $k \,=\, n-1$.
By subtracting this inequality
from~(\ref{eq:sthgp-total-degree-equation}), one finds that this
strongest subtour is equivalent to
$x(\delta(v)) \ge 1$
for the vertex $v \notin S$, which intuitively must be a very
strong inequality.
Quite surprising, however, the next-strongest subtour is
$k \,=\, n-2$ (even stronger than the smallest subtour $k=2$).

Both indicators for TSP subtours are symmetric with respect to
$k$ and $n-k$,
as expected.
The indicators for STHGP subtours, however, agree that the weakest
subtours are shifted slightly left of center, and that the position of
the weakest subtour shifts slowly to the left as $n \rightarrow \infty$.
Furthermore, large subtours are actually {\em stronger} than the
corresponding small subtours,
as seen in Figures~\ref{fig:reflected-epr-1000}
and~\ref{fig:reflected-cd-1000}
wherein the right half of each curve has been reflected upon the
left.
The region between the two curves indicates how much
stronger large-cardinality subtours are than the corresponding
small-cardinality subtours.

\begin{center}
\begin{minipage}[1]{4.75in}
 \begin{center}
 \begin{minipage}[t]{2.25in}
   \includegraphics[width=2.25in,clip=]{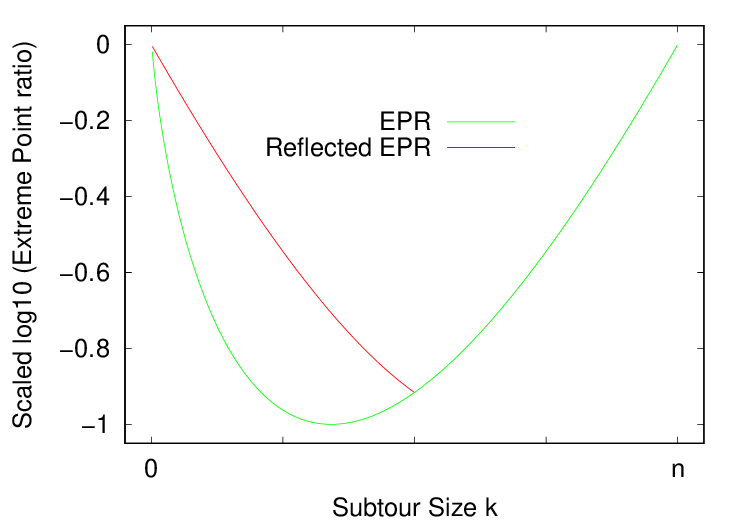}
   \captionof{figure}{Large STHGP subtours stronger than small
     subtours: EPR, $n=1000$ ~~(larger = stronger)}
   \label{fig:reflected-epr-1000}
 \end{minipage}
 \quad
 \begin{minipage}[t]{2.25in}
  \includegraphics[width=2.25in,clip=]{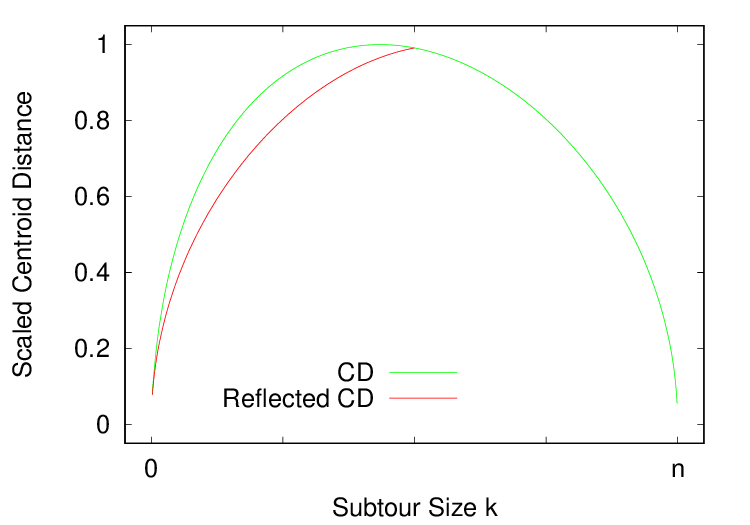}
  \captionof{figure}{Large STHGP subtours stronger than small
    subtours: CD, $n=1000$ ~~(smaller = stronger)}
  \label{fig:reflected-cd-1000}
 \end{minipage}
 \end{center}
\end{minipage}
\end{center}

These properties of STHGP subtours were unknown prior to the
computation of these indicators, and explain why computations that
require many subtour optimize/separate iterations often produce many
$k$-vertex subtours where
$0.25 \, n \,\le\, k \,\le\, 0.45 \, n$.

All previous constraint strengthening procedures in
{\ttfamily GeoSteiner} operate by removing vertices from the violated
subtour $S$~\cite{JuhlWarmeWinterZachariasen2018}.
If $|S| > 0.35 \, n$, this procedure actually {\em weakens} the
subtour (unless and until one removes sufficiently many vertices).
Violations in this vicinity are correctly strengthened by {\em adding}
vertices to them.
The companion paper exploits this property to great computational
effect~\cite{WarmeIndicators2}.

\begin{figure}[!ht]
\begin{center}
 \begin{minipage}[t]{2.25in}
  \includegraphics[width=2.25in,clip=]{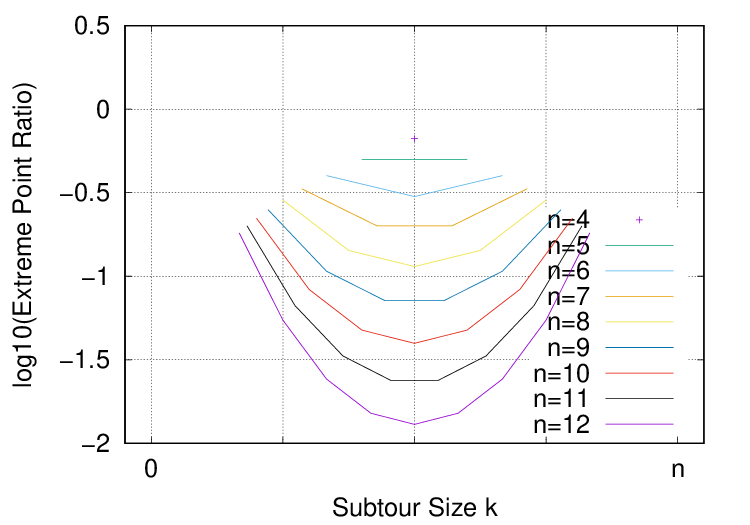}
  \captionof{figure}{EPR of TSP subtours for small $n$}%
  \label{fig:small-n-tsp-epr}
 \end{minipage}
 \quad
 \begin{minipage}[t]{2.25in}
  \includegraphics[width=2.25in,clip=]{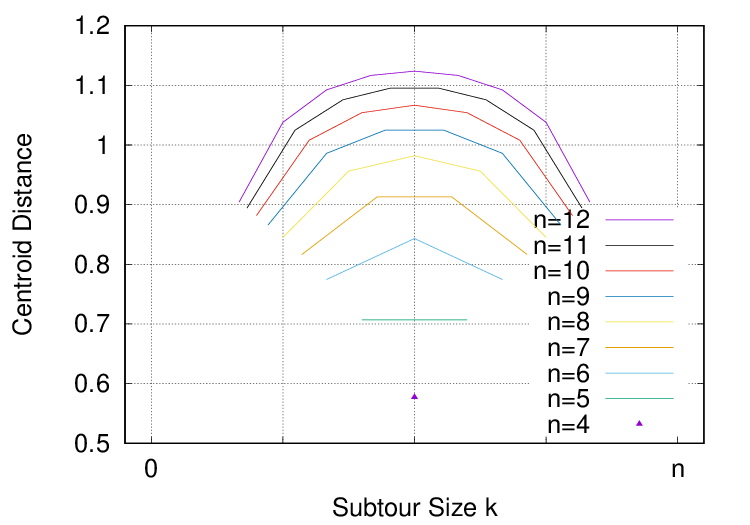}
  \captionof{figure}{CD of TSP subtours for small $n$}%
  \label{fig:small-n-tsp-cd}
 \end{minipage}
 \vspace*{\floatsep}
 \begin{minipage}[t]{2.25in}
  \includegraphics[width=2.25in,clip=]{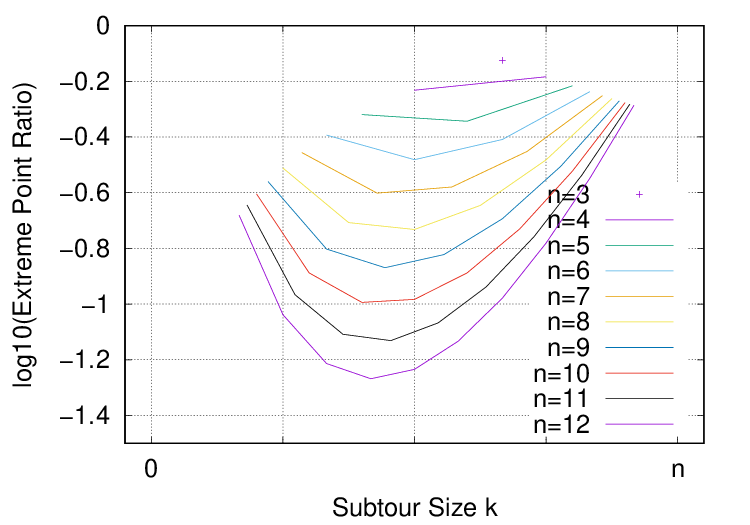}
  \captionof{figure}{EPR of STHGP subtours for small $n$}%
  \label{fig:small-n-sthgp-epr}
 \end{minipage}
 \quad
 \begin{minipage}[t]{2.25in}
  \includegraphics[width=2.25in,clip=]{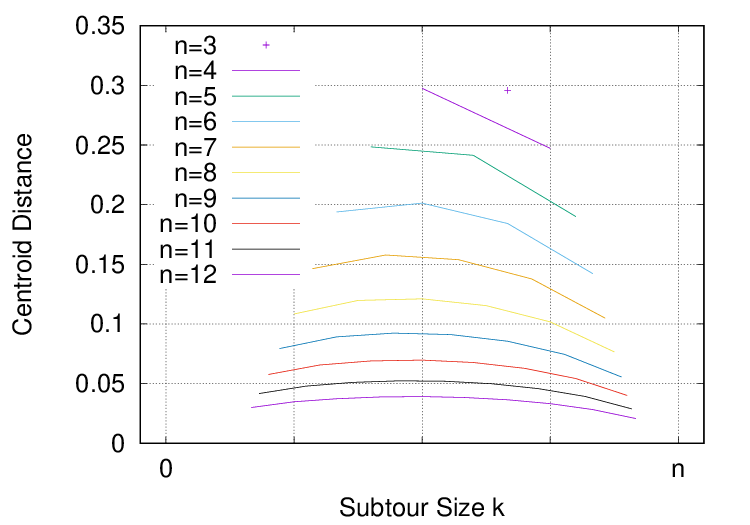}
  \captionof{figure}{CD of STHGP subtours for small $n$}%
  \label{fig:small-n-sthgp-cd}
 \end{minipage}
\end{center}
\end{figure}

Figure~\ref{fig:small-n-tsp-epr} through~\ref{fig:small-n-sthgp-cd}
present plots of EPR and CD of TSP and STHGP subtours for small
values of $n$.
These plots demonstrate that obtaining closed-forms for these
indicators are nice, but not essential.
Computing the indicators for small $n$ (e.g., by brute-force,
recursive enumeration of extreme points) can yield plots that reveal
the general trends for larger $n$.

\begin{figure}[!ht]
\begin{center}
 \begin{minipage}[t]{2.25in}
  \includegraphics[width=2.25in,clip=]{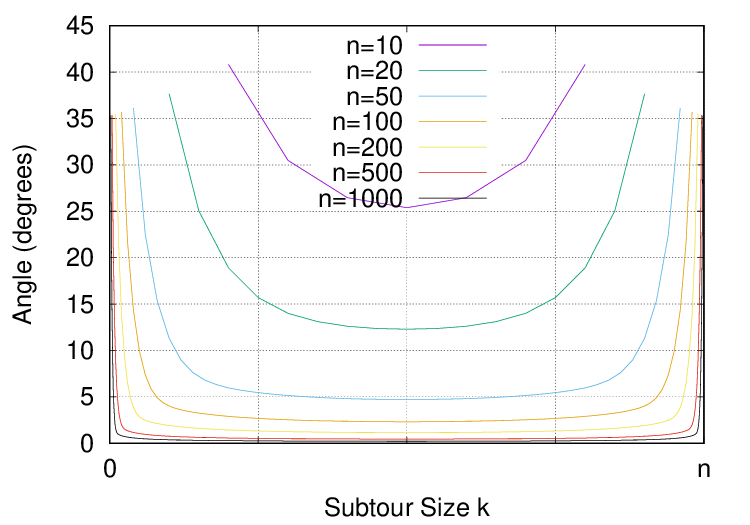}
  \captionof{figure}{Complementary subtour angles ($\theta(n,k,n-k,0)$ for $2 \le k \le n-2$)}%
  \label{fig:complementary-subtour-angles}
 \end{minipage}
 \quad
 \begin{minipage}[t]{2.25in}
  \includegraphics[width=2.25in,clip=]{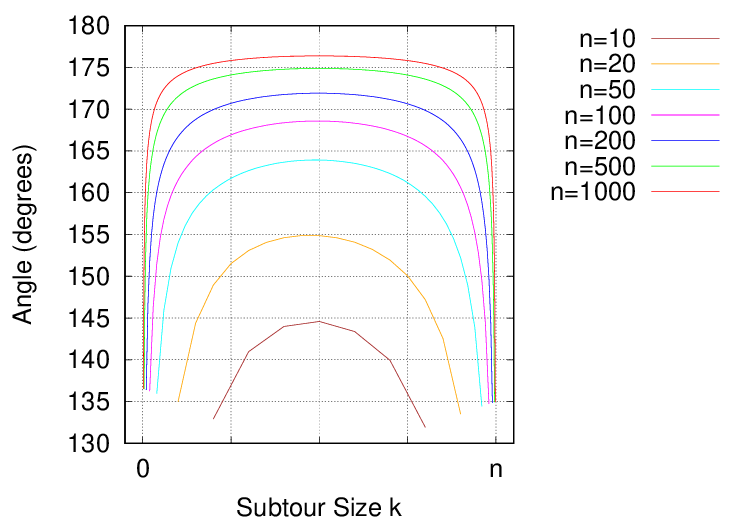}
  \captionof{figure}{Subtour angles: $k$ vs $k+1$ ($\theta(n,0,1,k)$ for $2 \le k \le n-2$)}%
  \label{fig:subtour-angles-1}
 \end{minipage}
\end{center}
\end{figure}

\begin{figure}[!ht]
\begin{center}
 \begin{minipage}[t]{2.25in}
  \includegraphics[width=2.25in,clip=]{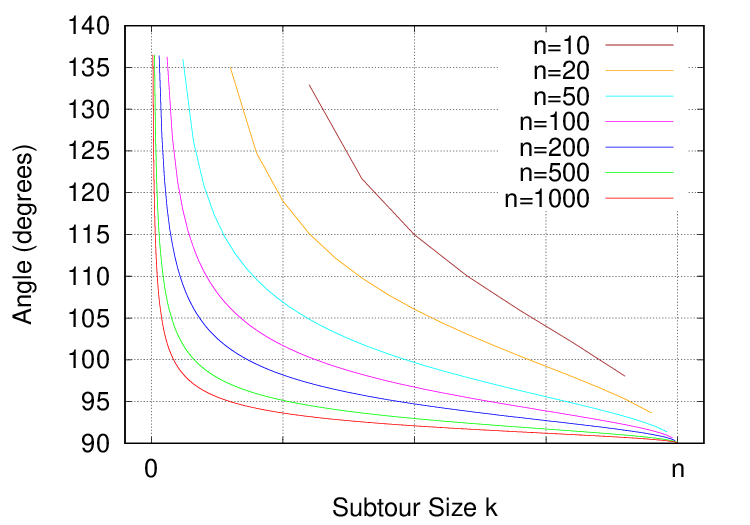}
  \captionof{figure}{Subtour angles: $2$ vs $k$ ($\theta(n,0,k-2,2)$ for $3 \le k \le n-1$)}%
  \label{fig:subtour-angles-2}
 \end{minipage}
 \quad
 \begin{minipage}[t]{2.25in}
  \includegraphics[width=2.25in,clip=]{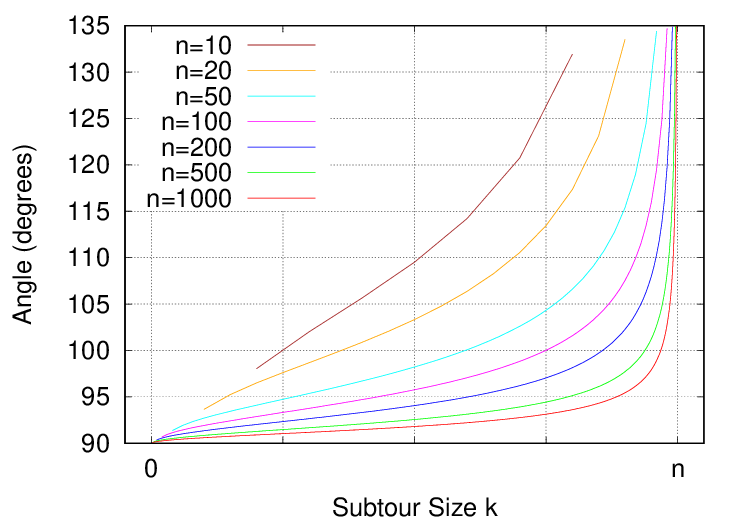}
  \captionof{figure}{Subtour angles: $k$ vs $n-1$ ($\theta(n,0,n-k-1,k)$ for $2 \le k \le n-2$)}%
  \label{fig:subtour-angles-3}
 \end{minipage}
\end{center}
\end{figure}

\begin{figure}[!ht]
\begin{center}
 \begin{minipage}[t]{2.25in}
  \includegraphics[width=2.25in,clip=]{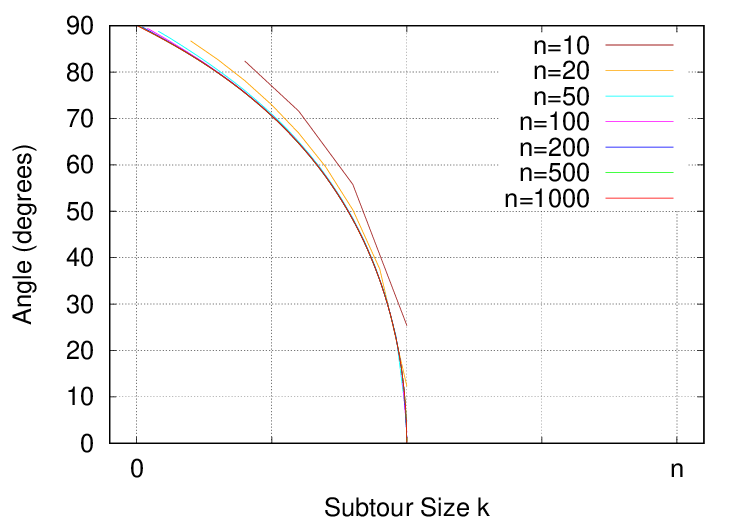}
  \captionof{figure}{Subtour angles: $k$ vs $k$, independent
    ($\theta(n,k,k,0)$ for $2 \le k \le \lfloor n/2 \rfloor$)}%
  \label{fig:subtour-angles-4}
 \end{minipage}
 \quad
 \begin{minipage}[t]{2.25in}
  \includegraphics[width=2.25in,clip=]{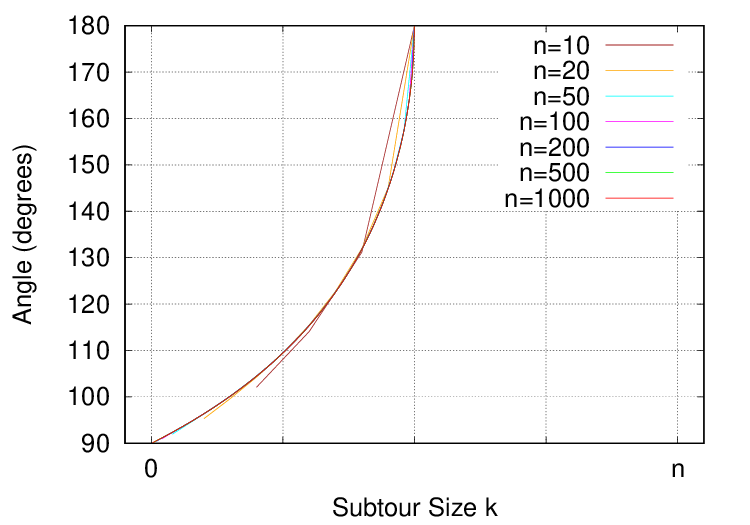}
  \captionof{figure}{Subtour angles: $k$ vs $n-k$, nested
    ($\theta(n,0,n-2k,k)$ for $2 \le k \le \lfloor n/2 \rfloor$)}%
  \label{fig:subtour-angles-5}
 \end{minipage}
\end{center}
\end{figure}

\begin{figure}[!ht]
\begin{center}
 \begin{minipage}[t]{2.25in}
  \includegraphics[width=2.25in,clip=]{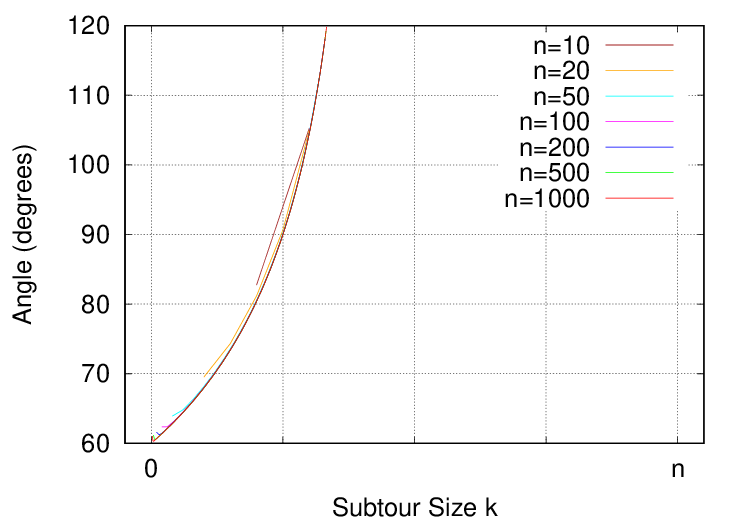}
  \captionof{figure}{Subtour angles: $2k$ vs $n-2k$, overlap by $k$
    ($\theta(n,k,n-3k,k)$ for $2 \le k \le \lfloor n/3 \rfloor$)}%
  \label{fig:subtour-angles-6}
 \end{minipage}
 \quad
 \begin{minipage}[t]{2.25in}
  \includegraphics[width=2.25in,clip=]{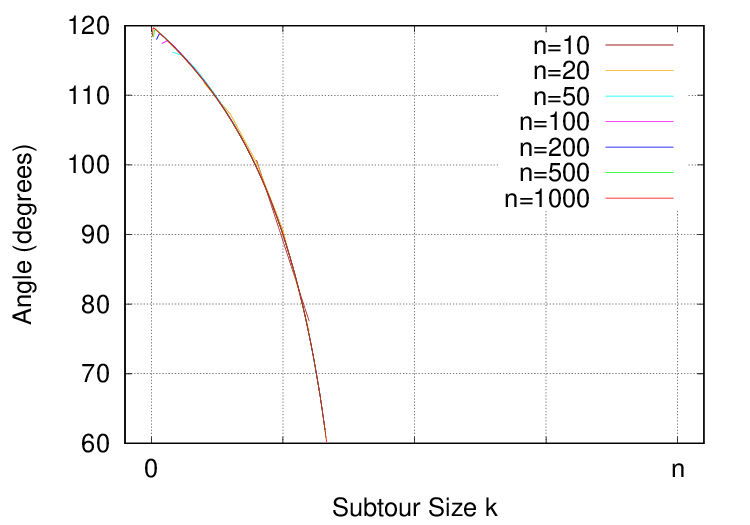}
  \captionof{figure}{Subtour angles: $2k$ vs $2k$, overlap by $k$
    ($\theta(n,k,k,k)$ for $2 \le k \le \lfloor n/3 \rfloor$)}%
  \label{fig:subtour-angles-7}
 \end{minipage}
\end{center}
\end{figure}

Figure~\ref{fig:complementary-subtour-angles} plots the interior angle
between complementary subtours $S$ and $V-S$ as a function of $k=|S|$
for several values of $n$.
(These angles are given by $\theta(n,k,n-k,0)$ for $2 \le k \le n-2$.)
Figure~\ref{fig:subtour-angles-1} plots the interior angle between
subtours $k$ vs $k+1$
($\theta(n,0,1,k)$, for $2 \le k \le n-2$).
Figure~\ref{fig:subtour-angles-2} plots the interior angle between
subtours $2$ vs $k$
($\theta(n,0,k-2,2)$, for $3 \le k \le n-1$).
Figure~\ref{fig:subtour-angles-3} plots the interior angle between
subtours $k$ vs $n-1$
($\theta(n,0,n-k-1,k)$, for $2 \le k \le n-2$).
Figure~\ref{fig:subtour-angles-4} plots the interior angle between
two independent subtours of cardinality $k$
($\theta(n,k,k,0)$, for $2 \le k \le \lfloor n / 2 \rfloor$).
Figure~\ref{fig:subtour-angles-5} plots the interior angle between
a subtour of cardinality $k$ nested within a subtour of cardinality
$n-k$
($\theta(n,0,n-2k,k)$ for $2 \le k \le \lfloor n / 2 \rfloor$).
Figure~\ref{fig:subtour-angles-6} plots the interior angle between
subtours of cardinality $2k$ and $n-2k$ that overlap by $k$ vertices
($\theta(n,k,n-3k,k)$ for $2 \le k \le \lfloor n / 3 \rfloor$).
Figure~\ref{fig:subtour-angles-7} plots the interior angle between
subtours of cardinality $2k$ and $2k$ that overlap by $k$ vertices
($\theta(n,k,k,k)$ for $2 \le k \le \lfloor n / 3 \rfloor$).

\begin{figure}[!ht]
\begin{center}
 \begin{minipage}[t]{2.25in}
  \includegraphics[width=2.25in,clip=]{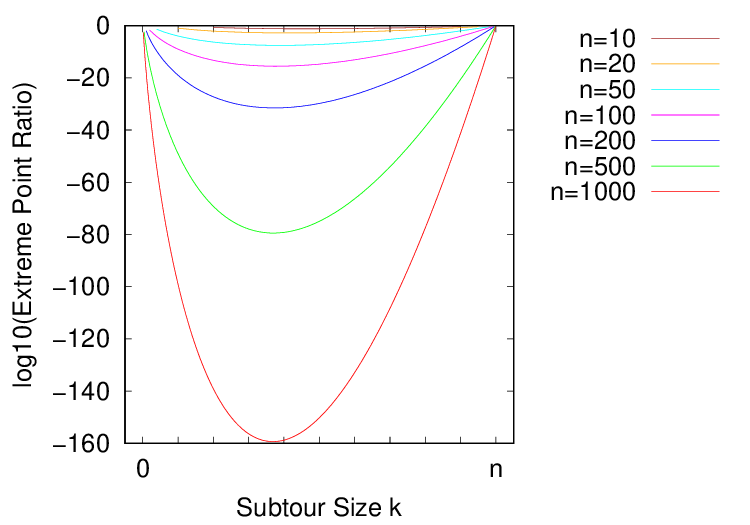}
  \captionof{figure}{Base 10 log of EPR for $k$-vertex
      subtours of $n$-vertex spanning tree in graph polytope}%
  \label{fig:stgp-subtour-epr-plot-all}
 \end{minipage}
 \quad
 \begin{minipage}[t]{2.25in}
  \includegraphics[width=2.25in,clip=]{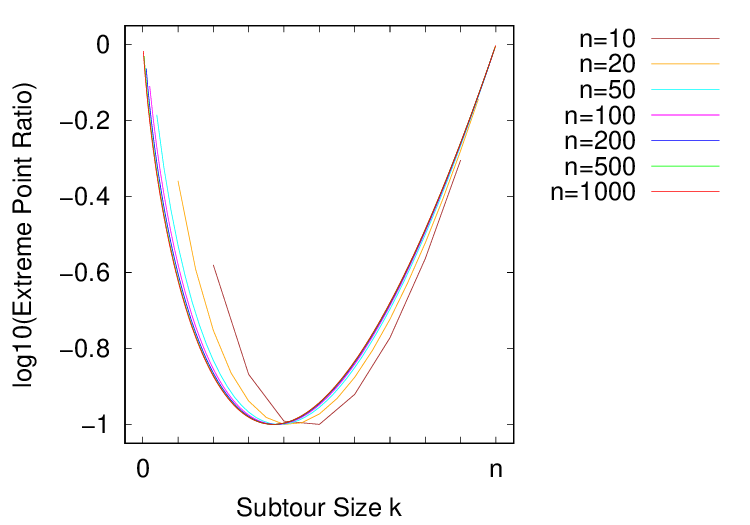}
  \captionof{figure}{EPR (scaled to minimum $y=-1$) for
      $k$-vertex subtours of $n$-vertex spanning tree in graph
      polytope}%
  \label{fig:stgp-subtour-scaled-epr-plot-all}
 \end{minipage}
\end{center}
\end{figure}

The spanning tree in graph polytope is well known.
Its complete linear description is:
\begin{eqnarray}
\label{eq:stgp-affine-hull}
	x(V) &=& |V| \,-\, 1, \\
\label{eq:stgp-subtours}
	x(X) &\le& |S| \,-\, 1 ~~~~
		\hbox{for all $S \subset V$, $|S|\ge 2$}, \\
\label{eq:stgp-nonneg}
	x_e &\ge& 0   ~~~~~~ \hbox{for all $e \in E$}.
\end{eqnarray}
where~(\ref{eq:stgp-affine-hull}) is the affine hull,
(\ref{eq:stgp-subtours}) are the subtour inequalities, and
(\ref{eq:stgp-nonneg}) are the non-negativity inequalities.
Figure~\ref{fig:stgp-subtour-epr-plot-all}
presents a plot of EPRs of subtour inequalities of the spanning
tree in graph polytope, which is plotted with a log scale for the $y$
coordinate. 
This clearly shows the ratios rapidly decreasing with $n$, but
cause the details of each curve to be blurred to nearly flat lines.
To remedy this, we present additional
Figure~\ref{fig:stgp-subtour-scaled-epr-plot-all}
in which the $y$ coordinates of each curve are scaled so that its
minimum log(ratio) is -1.
These scaled EPR curves in
Figure~\ref{fig:stgp-subtour-scaled-epr-plot-all}
for spanning tree in graph subtours are seen to be almost
indistinguishable from the scaled EPR curves in
Figure~\ref{fig:sthgp-sec-ratio-scaled-plot-all}
for spanning tree in hypergraph subtours.

\begin{figure}[!ht]
\begin{center}
 \includegraphics[width=3.0in,clip=]{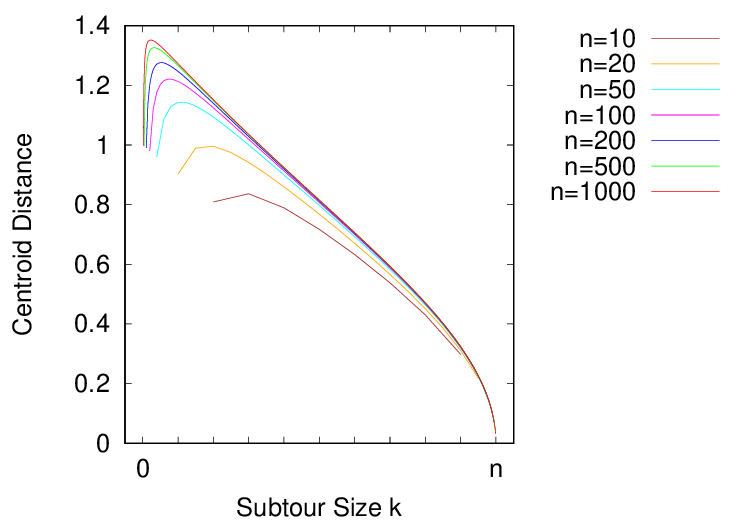}
 \captionof{figure}{CD for $k$-vertex subtours of $n$-vertex
     spanning tree in graph polytope}%
 \label{fig:stgp-subtour-cd-plot-all}
\end{center}
\end{figure}

Figure~\ref{fig:stgp-subtour-cd-plot-all}
presents a plot of CDs of subtour inequalities of the spanning
tree in graph polytope.
Note that these curves are shockingly different in shape from the
corresponding CD curves for subtours of the spanning tree in
hypergraph polytope
(Figure~\ref{fig:sthgp-sec-dist-scaled-plot-all}).
Figures~\ref{fig:stgp-subtour-cd-plot-all}
and~\ref{fig:stgp-subtour-scaled-epr-plot-all}
display contradictory estimates of STGP subtour strength (as a
function of subtour cardinality).
EPR and CD cannot therefore both be good estimates of actual STGP
subtour strength.

We find it unlikely that small cardinality STGP subtours are actually
as weak as CD predicts:
(1) Although one does not typically use an LP solver and separation
algorithm to compute MSTs in conventional graphs, this is precisely
how computational strength of $\STGP{n}$ subtours would necessarily be
evaluated;
and
(2) Computational evidence with
GeoSteiner~\cite{JuhlWarmeWinterZachariasen2018}
demonstrates small cardinality STHGP subtours are quite strong on
hypergraphs that are ``almost'' conventional graphs
($|V| < |E| < 1.6|V|$,
there are $n-1$ edges of cardinality 2 forming an MST and few edges
have cardinality > 5).
In Section~\ref{sec:unexpected-stgp-cd-behavior} we present some
evidence that dimensionality of the space plays a role in the
unexpected shape of the curves shown in
Figure~\ref{fig:stgp-subtour-cd-plot-all}
versus those for STHGP shown in
Figure~\ref{fig:sthgp-sec-dist-scaled-plot-all}.
Note further that EPR is independent of the dimensionality of the
polyhedral space.
The surprising and unexpected CD results for STGP subtours lead us to
conclude EPR provides a better strength estimate than CD.
We therefore believe that
Figure~\ref{fig:reflected-epr-1000}
is a more accurate estimate of the stregth advantage of large
cardinality STHGP subtours (over those of correspondingly small
cardinality) than
Figure~\ref{fig:reflected-cd-1000}.

\subsection{Comparison of Indicators for STHGP Subtours}
\label{sec:compare-indicators}

It is clear from Figures~\ref{fig:sthgp-sec-ratio-scaled-plot-all}
and~\ref{fig:sthgp-sec-dist-scaled-plot-all} that EPR and CD
indicators are in remarkable agreement regarding the relative
(predicted) strength of STHGP subtours.
We now examine the extent of this agreement more closely.
Table~\ref{tab:weakest-sthgp-subtours-by-indicator} presents the
cardinality of the subtour that each indicator predicts to be weakest,
for several $n$.
As $n$ increases, EPR and CD indicators appear to increasingly
disagree as which subtours are weakest.
\begin{table}[!ht]
\setlength{\tabcolsep}{3.0pt}
\begin{center}
 \captionof{table}{Weakest STHGP subtours, as predicted by EPR and CD
   indicators}%
\begin{tabular}{|r|r|r|} \hline
\multicolumn{1}{|c}{n} &
\multicolumn{1}{|c}{Weakest EPR} &
\multicolumn{1}{|c|}{Weakest CD} \\
\hline
10 & 4 & 5 \\
100 & 35 & 45 \\
1000 & 342 & 434 \\
\hline
\end{tabular}
\label{tab:weakest-sthgp-subtours-by-indicator}
\end{center}
\end{table}

\begin{figure}[!ht]
\begin{center}
 \includegraphics[width=2.25in,clip=]{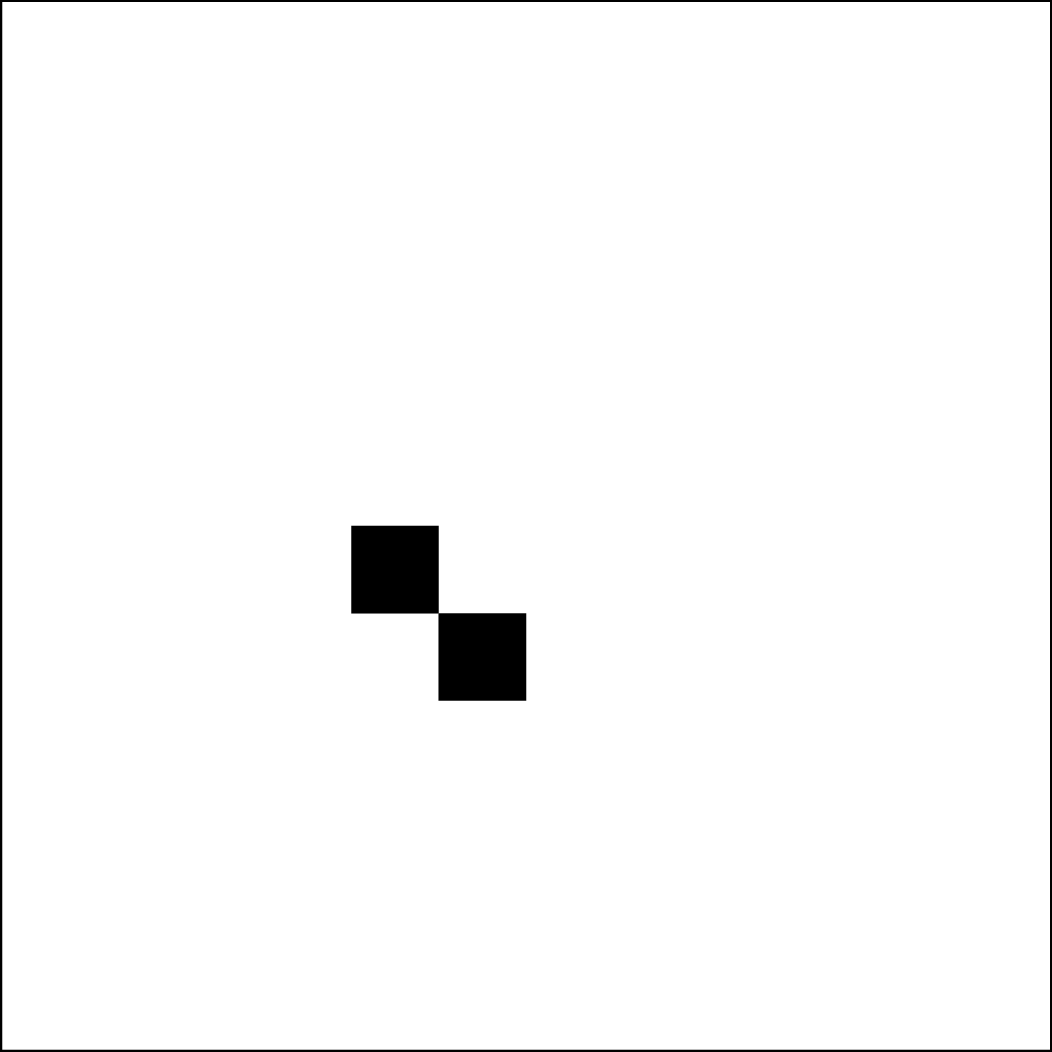}
 \captionof{figure}{Detailed comparison of EPR vs CD indicators on
   STHGP subtours $(n=10)$, both axes are subtour size,
   regions of disagreement are shown in black}%
 \label{fig:sthgp-compare-10}
\end{center}
\vspace*{\floatsep}
\begin{center}
 \begin{minipage}[t]{2.25in}
  \includegraphics[width=2.25in,clip=]{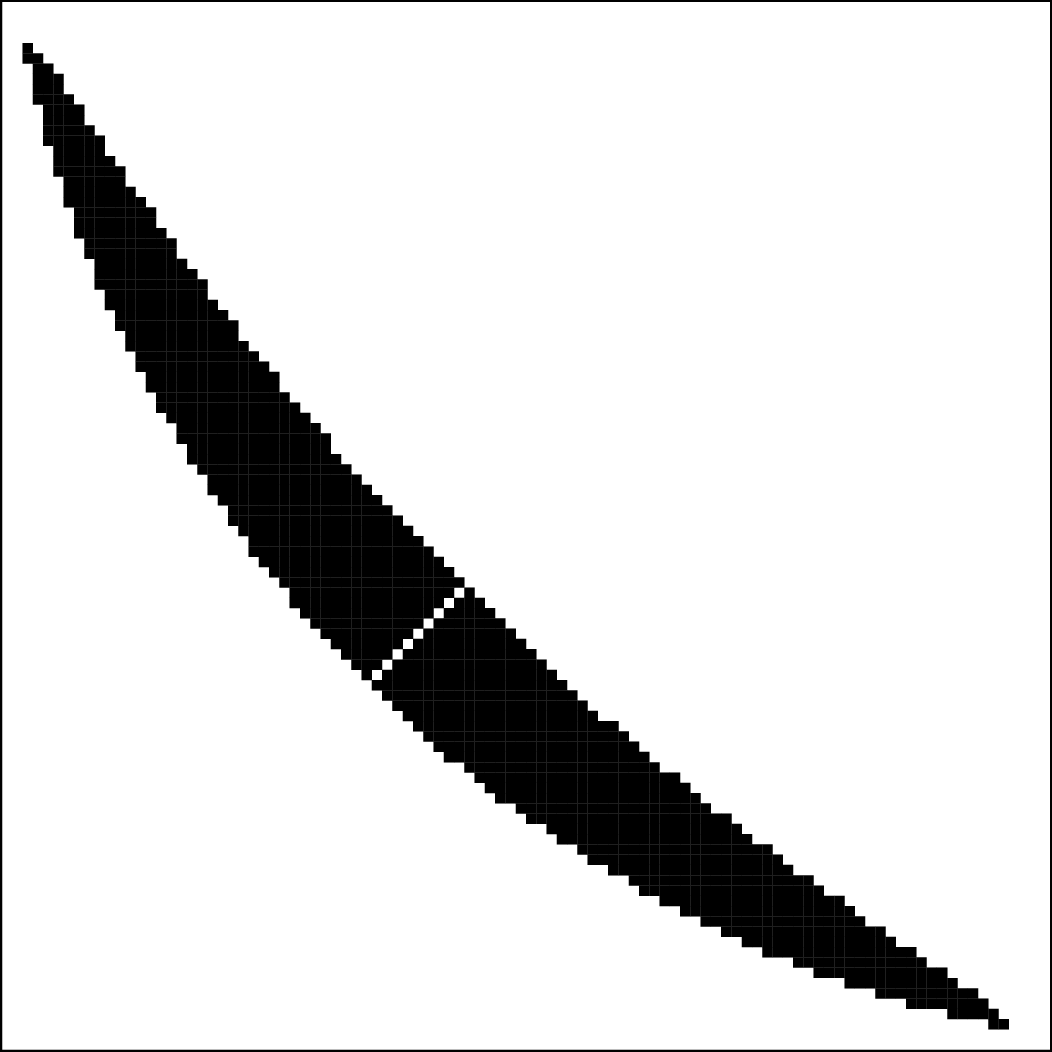}
  \captionof{figure}{Detailed comparison of EPR vs CD indicators on
    STHGP subtours $(n=100)$, both axes are subtour size,
    regions of disagreement are shown in black}%
  \label{fig:sthgp-compare-100}
 \end{minipage}
 \quad
 \begin{minipage}[t]{2.25in}
  \includegraphics[width=2.25in,clip=]{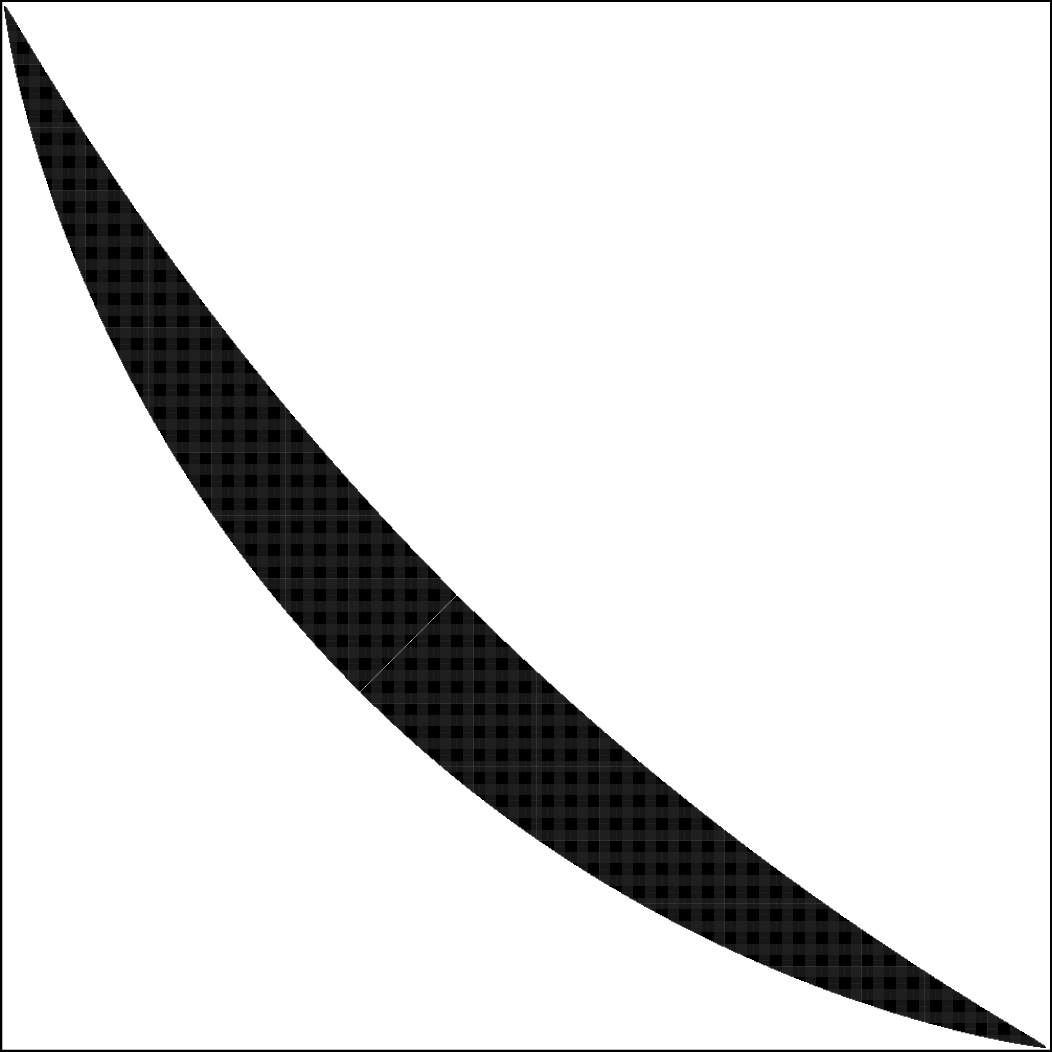}
  \captionof{figure}{Detailed comparison of EPR vs CD indicators on
    STHGP subtours $(n=1000)$, both axes are subtour size,
    regions of disagreement are shown in black}%
  \label{fig:sthgp-compare-1000}
 \end{minipage}
\end{center}
\end{figure}

\noindent
We now consider a more detailed comparison of the two indicators.
Let $n > 3$, and $2 \,\le\, k_1,\,k_2 \,<\, n$.
Let $m_1(n,k)$ denote EPR, and
$m_2(n,k)$ denote CD.
For EPR we have the following cases:
\begin{center}
\begin{tabular}{cc}
  $m_1(n,k_1) \,<\, m_1(n,k_2)$ & $k_1$ weaker than $k_2$ \\
  $m_1(n,k_1) \,=\, m_1(n,k_2)$ & $k_1$ same strength as $k_2$ \\
  $m_1(n,k_1) \,>\, m_1(n,k_2)$ & $k_1$ stronger than $k_2$
\end{tabular}
\end{center}
For CD we have the following cases:
\begin{center}
\begin{tabular}{cc}
  $m_2(n,k_1) \,<\, m_2(n,k_2)$ & $k_1$ stronger than $k_2$ \\
  $m_2(n,k_1) \,=\, m_2(n,k_2)$ & $k_1$ same strength as $k_2$ \\
  $m_2(n,k_1) \,>\, m_2(n,k_2)$ & $k_1$ weaker than $k_2$
\end{tabular}
\end{center}
Figures~\ref{fig:sthgp-compare-10}--\ref{fig:sthgp-compare-1000} plot
a black square for each pair $(k_1,k_2)$ such that $m_1$ and $m_2$
{\em disagree} regarding the relative strength of subtours of size
$k_1$ versus subtours of size $k_2$.
The two indicators agree much more often than they disagree, and their
largest region of disagreement concerns the very weakest subtours.
This is fortunate ---
we want to use the strongest possible constraints when optimizing, and
can forgive quibbles over which subtours are the weakest, so long as
we can confidently identify those subtours that are strongest.

\section{Traveling Salesman Polytope: $TSP(n)$}
\label{sec:tsp-polytope}

Let $n \in \Z$, $n \ge 3$.  Let $G = (V,E) = K_n$ be the complete
graph on $n$ vertices.
Let $m = |E| = {n \choose 2}$.
We work in the space $\R^m$ with incidence vectors that represent
subsets of $E$.
Let $T(n)$ be the set of all incidence vectors $x \in \B^m$
that correspond to tours of $G$
(i.e., each incidence vector $x \in T(n)$ corresponds to an
$E' \subseteq E$ such that $|E'| = n$,
and the edges $e \in E'$ form a Hamiltonian cycle in $G$).
We define the {\em Traveling Salesman Polytope}
\begin{displaymath}
\TSP{n} \,=\, \CHull{T(n)}.
\end{displaymath}

\noindent
The subtour relaxation of $\TSP{n}$ is defined as
\begin{eqnarray}
\label{eq:tsp-degree-two-equations}
	x(\delta(v)) &=& 2
		~~~~~~~~~ \forall v \in V, \\
\label{eq:tsp-subtour-facet-def}
	x(\delta(S)) &\ge& 2
  ~~~~~~~\hbox{for all $S \subset V$ such that $2 \le |S| \le |V|-2$}, \\
\label{eq:tsp-nonneg-facet-def}
	x_e &\ge& 0 ~~~~~~\hbox{for all $e \in E$}.
\end{eqnarray}
Equations~(\ref{eq:tsp-degree-two-equations}) are known as the
{\em degree two} equations, since they force each vertex $v \in V$ to
have degree 2.
These $n$ equations are the only equations satisfied by all tours,
so they define the affine hull of $\TSP{n}$.
(\ref{eq:tsp-subtour-facet-def}) are the subtour elimination
constraints (more simply called subtours or SECs).
(\ref{eq:tsp-nonneg-facet-def}) are the non-negativity constraints.

Many classes of facet-defining inequalities are known for $\TSP{n}$.
This paper considers three of these classes:
(1) the non-negativity constraints~(\ref{eq:tsp-nonneg-facet-def}),
(2) the subtours~(\ref{eq:tsp-subtour-facet-def}),
and
(3) the comb inequalities
\begin{displaymath}
x(E(H)) \,+\, \sum_{i=1}^t x(E(T_i)) ~\le~
	|H| \,+\, \sum_{i=1}^t |T_i| \,-\, {{3 t \,+\, 1} \over 2}
\end{displaymath}
where $t \ge 3$ is an odd integer,
all $T_i \subset V$ are disjoint,
$H \cap T_i$ and $T_i \setminus H$ are non-empty for all $1 \le i \le t$.

We now prove EPR
(as presented in Table~\ref{tab:tsp-extreme-point-ratios}) for the
various TSP facets.
We start with a few simple preliminary results.

\begin{theorem}
\label{th:tsp-n-tours}
Let $n \ge 3$.
Let $G = (V,E) = K_n$ be the complete graph on $n$ vertices.
There are $(n-1)!/2$ tours in $G$ and
extreme points in $\TSP{n}$.
\end{theorem}
\begin{MyProof}
The extreme points of $\TSP{n}$ correspond to the tours of $K_n$.
There are $n!$ distinct permutations of $n$ vertices.
We divide by $n$ because the vertex on which a tour begins and ends
does not matter.
We also divide by 2 because the clockwise and counter-clockwise
traversals of a tour yield different vertex sequences for the same tour.
There are no other symmetries to divide out.
\qed
\end{MyProof}

\begin{figure}[!ht]
\begin{center}
  \begin{minipage}[t]{2.25in}
   \includegraphics[width=2.25in,clip=]{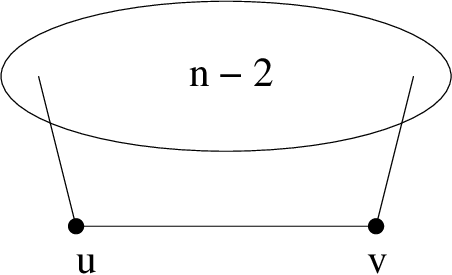}
   \captionof{figure}{Tours containing edge $e=(u,v)$}%
   \label{fig:tours-with-edge-e}
  \end{minipage}
  \quad
  \begin{minipage}[t]{2.25in}
   \includegraphics[width=2.25in,clip=]{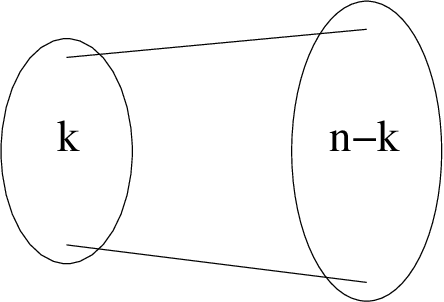}
   \captionof{figure}{Tours incident to $k$-vertex subtour}%
   \label{fig:tsp-k-subtour}
  \end{minipage}
\end{center}
\end{figure}

\begin{theorem}
\label{th:tsp-tours-for-edge-e}
Let $n \ge 3$,
$G = (V,E) = K_n$, and
$e = (u,v) \in E$.
There are $(n-2)!$ distinct tours that contain edge $e$.
\end{theorem}
\begin{MyProof}
As shown in Figure~\ref{fig:tours-with-edge-e}, there are $n-2$
vertices besides $u$ and $v$.
There are $(n-2)!$ permutations of these vertices, and each such
permutation gives rise to a distinct tour.
\qed
\end{MyProof}

\subsection{EPR for TSP Non-Negativity Inequalities}
\label{sec:TSP-nonneg-epr}

\begin{theorem}
\label{th:tsp-nonneg-extreme-point-ratio}
Let $n \ge 3$, $G = (V,E) = K_n$, and $e \in E$.
Then EPR for the TSP facet $x_e \,\ge\, 0$ is
\begin{displaymath}
r(n)	~=~
	{{n-3} \over {n-1}}.
\end{displaymath}
\end{theorem}
\begin{MyProof}
By Theorem~\ref{th:tsp-n-tours}, the total number of tours is
$(n-1)!/2$.
The number of extreme points for which $x_e = 0$ is just
the total number extreme points minus the number of extreme points
containing edge $e$, which by Theorem~\ref{th:tsp-tours-for-edge-e}
is $(n-2)!$.
EPR is therefore
\begin{displaymath}
r(n)	~=~
    ~=~ {{ {{(n-1)!} \over 2} \,-\, (n-2)! } \over { {{(n-1)!} \over 2} }}
    ~=~ 1 \,-\, {2 \over {n-1}}
    ~=~ {{n-3} \over {n-1}}.
\tag*{\qed}
\end{displaymath}
\end{MyProof}

\subsection{EPR for TSP Subtour Inequalities}
\label{sec:TSP-subtour-epr}

\begin{theorem}
\label{th:tsp-subtour-extreme-point-ratio}
Let $n \ge 3$,
$G = (V,E) = K_n$,
$S \subset V$ such that $|S| \ge 2$,
and $k = |S|$.
Then EPR for the TSP subtour facet
$x(\delta(S)) \ge 2$ is
\begin{displaymath}
r(n,k) ~=~
	{n \over {n \choose k}}.
\end{displaymath}
\end{theorem}
\begin{MyProof}
As shown in Figure~\ref{fig:tsp-k-subtour}, there are
${{k! \, (n-k)!} \over 2}$
tours incident to a $k$-vertex subtour.
(We must divide by 2 because the permutations produce both a clockwise
and a counter-clockwise traversal of each tour.)
This yields an EPR of
\begin{displaymath}
r(n,k)
    ~=~ {{ {{k! \, (n-k)!} \over 2} } \over { {{(n-1)!} \over 2} }}
    ~=~ {{ n \, k! \, (n-k)!} \over {n!}}
    ~=~ {n \over {n \choose k}}.
\tag*{\qed}
\end{displaymath}
\end{MyProof}

No closed-form is currently known for EPR of TSP 3-toothed comb
inequalities.

We now prove CD (as presented in
Table~\ref{tab:tsp-centroid-distances}) for the various TSP facets.

\subsection{Centroid of $\TSP{n}$}
\label{sec:TSP-centroid}

\begin{theorem}
\label{th:tsp-centroid}
Let $C$ be the centroid of $\TSP{n}$.
Then for each $e \in E$ we have $C_e \,=\, {2 \over {n-1}}$
and $C$ can itself be expressed as
\begin{displaymath}
	C ~=~ {2 \over {n-1}} \, {e_m}
\end{displaymath}
where $e_m$ is the $m$-element column vector of all ones.
\end{theorem}
\begin{MyProof}
For each $e \in E$, there are
$(n-2)!$
tours containing edge $e$ and
${{(n-1)!} \over 2}$ tours total.
The corresponding element $C_e$ of the centroid is therefore
$(n-2)! / ((n-1)!/2) \,=\, {2 \over {n-1}}$.
\qed
\end{MyProof}

\subsection{CD of TSP Non-Negativity Inequalities}
\label{sec:TSP-nonneg-cd}

The affine hull of $\TSP{n}$ is the intersection of $n$
equations (a variable number).
When working with symbolic closed-forms, eliminating variables using a
{\em variable} number of equations requires different techniques than
for eliminating variables using a fixed number of equations.
We address this by appeal to the symmetry of $\TSP{n}$.

\begin{theorem}
Let $n \ge 3$, $G = (V,E) = K_n$, and $e \in E$.
Then CD squared for the TSP facet $x_e \,\ge\, 0$
is
\begin{displaymath}
	d^2(n) ~=~ {4 \over {(n-1) \, (n-3)}}.
\end{displaymath}
\end{theorem}
\begin{MyProof}
Refer to
Figure~\ref{fig:tsp-nonneg-fdist}, where edge $e=(u,v)$ is clearly
missing.
\MyBeginFig
 \begin{center}
  \includegraphics[width=1.5in,clip=]{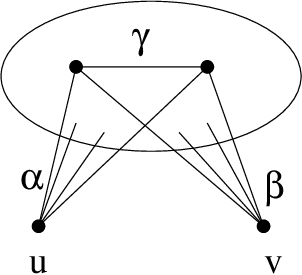}
  \captionof{figure}{CD for $x_e \ge 0$ for edge $e=(u,v)$}
  \label{fig:tsp-nonneg-fdist}
 \end{center}
\MyEndFig
Consider the point $y \in H \cap \AHull{P}$ that is closest to the
centroid $C$.
Point $y$ is unique because the point closest to $C$ in any hyperplane
is unique.
$G$ has 4 kinds of edges in this figure:
edge $e=(u,v)$ (having weight $x_e = 0$);
those incident to vertex $u$;
those incident to vertex $v$;
and all others.

Let $T = V - \lbrace u,\, v \rbrace$.
By symmetry of $\TSP{n}$, $y$ must remain fixed over all
permutations of vertices $T$.
Let $p,q \in T$.
Swapping $p$ and $q$:
(1) swaps edges $e_1 = (u,p)$ and $e_2 = (u,q)$ so we must have
$y_{e_1} = y_{e_2}$;
(2) swaps edges $e_3 = (v,p)$ and $e_4 = (v,q)$ so we must have
$y_{e_3} = y_{e_4}$;
(3) for edge $e_5 = (p,q)$, we trivially must have
$y_{e_5} = y_{e_5}$.
Let $\alpha$, $\beta$ and $\gamma$ be the common $y$-component value
for edges
$\lbrace (u,p) : p \in T \rbrace$,
$\lbrace (v,p) : p \in T \rbrace$, and
$\lbrace (p,q) : p,q \in T \rbrace$, respectively.
We now write down three of the {\em degree 2}
equations~(\ref{eq:tsp-degree-two-equations}): one for vertex $u$, one
for vertex $v$, and a third equation for some vertex in $T$
\begin{displaymath}
 \begin{array}{rrrcc}
   (n-2) \, \alpha &		     &				&=& 2 \\
		   & (n-2) \, \beta  &				&=& 2 \\
	    \alpha &	 + \, \beta  & + \, (n-3) \, \gamma	&=& 2 \\
 \end{array}
\end{displaymath}
whose solution is
\begin{displaymath}
    \begin{array}{cc}
	\alpha ~=~ \beta ~=~ {2 \over {n-2}}, &
	~~~~~
	\gamma  ~=~ {{2 \, (n-4)} \over {(n-2) \, (n-3)}}. \\
    \end{array}
\end{displaymath}
CD squared is then
\begin{align*}
d^2(n) &= (n-2) \, (\alpha \,-\, C_e)^2
	\,+\, (n-2) \, (\beta \,-\, C_e)^2
	\,+\, {{n-2} \choose 2} \, (\gamma \,-\, C_e)^2 \\
&~~~	\,+\, 1 \, (0 \,-\, C_e)^2 \\
&=	{4 \over {(n-1) \, (n-3)}}.
\tag*{\qed}
\end{align*}
\end{MyProof}

\subsection{CD of TSP Subtour Inequalities}
\label{sec:TSP-subtour-cd}

\begin{theorem}
Let $n \ge 3$,
$G = (V,E) = K_n$,
$S \subset V$ such that $|S| \ge 2$,
and $k = |S|$.
Then CD squared for the TSP subtour facet
$x(\delta(S)) \ge 2$ is
\begin{displaymath}
d^2(n,k) ~=~
	{{2 \, (k-1) \, (n-2) \, (n-k-1)}
	  \over
	 {k \, (n-1) \, (n-k)}}.
\end{displaymath}
\end{theorem}
\begin{MyProof}
Refer to Figure~\ref{fig:tsp-sec-dist}.
\MyBeginFig
 \begin{center}
  \includegraphics[width=1.75in,clip=]{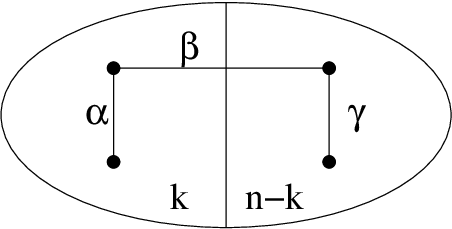}
 \end{center}
 \begin{center}
  \caption{CD for $k$ vertex TSP subtour}%
  \label{fig:tsp-sec-dist}
 \end{center}
\MyEndFig
Consider the point $y \in H \cap \AHull{P}$ that is closest to $C$.
Since $y$ is unique, is must be fixed under all permutations of
vertices in $S$, and all permutations of vertices $V-S$.
For distince edges $e_1, e_2 \in E(S)$, there are 2 cases:
(1) $e_1 = (t,u)$ and $e_2 = (u,v)$ for which swapping $t$ and $v$
swaps edges $e_1$ and $e_2$; and
(2) $e_1 = (t,u)$ and $e_2 = (v,w)$, for which swapping $t$ with $v$
and $u$ with $w$ swaps edges $e_1$ and $e_2$.
This implies that $y_{e_1} = y_{e_2}$ for all $e_1,e_2 \in E(S)$.
Similar arguments apply for edges in $E(S:V-S)$ and $E(V-S)$.
Let $\alpha$, $\beta$ and $\gamma$ be the common $y$ component values
for edges in $E(S)$, $E(S:V-S)$ and $E(V-S)$, respectively.
We write down 3 simultaneous equations:
the {\em degree 2} equation for a vertex in the subtour of
size $k$,
the {\em degree 2} equation for a vertex outside the subtour (in the
$n-k$ subset), and the subtour inequality
\begin{displaymath}
 \begin{array}{rrrcc}
	(k-1) \, \alpha & + \, (n-k) \, \beta & &=& 2 \\
	& k \, \beta & + \, (n-k-1) \, \gamma &=& 2 \\
	& k \, (n-k) \, \beta & &=& 2 \\
 \end{array}
\end{displaymath}
whose solution is
\begin{displaymath}
  \begin{array}{ccc}
	\alpha	~=~ {2 \over k}, &
	~~~~
	\beta	~=~ {2 \over {k \, (n-k)}}, &
	~~~~
	\gamma	~=~ {2 \over {(n-k)}}. \\
  \end{array}
\end{displaymath}
CD squared is then
\begin{align*}
d^2(n,k)
 &= {k \choose 2} (\alpha \,-\, C_e)^2
     \,+\, k \, (n-k) \, (\beta \,-\, C_e)^2
     \,+\, {{n-k} \choose 2} \, (\gamma \,-\, C_e)^2 \\
 &= {{2 \, (k-1) \, (n-2) \, (n-k-1)}
      \over
      {k \, (n-1) \, (n-k)}}.
\tag*{\qed}
\end{align*}
\end{MyProof}
It is not difficult to show that
$y \,=\, {1 \over \lambda} \sum_{x \in (H \cap X) } x$
where
$\lambda \,=\, |H \cap X| \,=\, {{k! \, (n-k)!} \over 2} \,>\, 0$.
(For example, for $e \in E(S)$
there are $(k-1)!$ Hamiltonian paths in $E(S)$ containing $e$, each
having $(n-k)!$ tour completions in $E(V-S)$, giving
$y_e = (k-1)!(n-k)!/\lambda = 2/k = \alpha$.)
This establishes that for facet $F = \CHull{H \cap X}$:
(1) $y \in F$;
(2) CD and weak CD are equal; and
(3) $y$ is the centroid of $F$.

\subsection{CD of TSP 3-Toothed Comb Inequalities}
\label{sec:TSP-3-toothed-comb-CD}

Recall Figure~\ref{fig:standard-TSP-3comb} showing a 3-toothed comb
using the standard notation for the handle $H$, and teeth $T_i$.
Recall also Figure~\ref{fig:our-TSP-3comb}, and the 8 parameters
$b_1,\, t_1,\, b_2,\, t_2,\, b_3,\, t_3,\, h$ and $o$ that we use to
describe a particular 3-toothed comb inequality.

\begin{theorem}
\label{th:tsp-3-toothed-comb-centroid-distance}
Let $n \ge 6$,
$G = (V,E) = K_n$.
Let $B_1,\, T_1,\, B_2,\, T_2,\, B_3,\, T_3,\, H,\, O \subset V$
be a partition of $V$ such that only $H$ and $O$ are allowed to be
empty.
For $1 \le i \le 3$ let $b_i = |B_i|$ and $t_i = |T_i|$.
Let $h = |H|$ and $o = |O|$.
Then this partition forms a valid TSP 3-toothed comb inequality
\begin{displaymath}
x(E(B_1 \cup B_2 \cup B_3 \cup H))
\,+\, \sum_{i=1}^3 x(E(B_i \cup T_i))
~\le~
|B_1 \cup B_2 \cup B_3 \cup H|
\,+\, \sum_{i=1}^3 |B_i \cup T_i|
\,-\, 5
\end{displaymath}
with teeth
$B_i \cup T_i$ for $1 \le i \le 3$
and handle $H \cup \bigcup_{i=1}^3 B_i$.
CD squared of this 3-toothed comb inequality is
\begin{displaymath}
d_w^2(b_1,\,t_1,\,b_2,\,t_2,\,b_3,\,t_3,\,h,\,o) ~=~
	{{2 \, (n-2) \, A^2} \over {(n-1) \, B}}
\end{displaymath}
where $A$ and $B$ are given in Table~\ref{tab:3comb-dist-subexprs}
(subject to the associated discussion on interpretation of subscripts
in the expressions for $A$ and $B$).
\end{theorem}

We present an overview of this computation, not a detailed proof.
Computation of CD for 3-toothed combs is similar to
that for the non-negativity constraints and the subtours, although the
math itself is considerably more voluminous and requires assistance
from computer algebra tools.
Let $P=\TSP{n}$.
Let $F$ be the hyperplane of the given 3-toothed comb inequality.
There are 8 classes of vertices, as shown in
Figure~\ref{fig:our-TSP-3comb}.
This gives rise to
${{8+1} \choose 2} \,=\, 36$ classes of edges.
Let $y$ be the point in
$F \cap \mathrm{aff}(P)$
that is closest to $C$.
We argue by symmetry of polytope $\TSP{n}$ that
$y_{e_1} = y_{e_2}$ for all edges $e_1$ and $e_2$ that belong to the
same class.
We let $x_{u,v}$ represent this comon value of the various $y_e$ (for
all edges $e$ from vertex class $u$ to vertex class $v$).
We construct a system of 9 equations in 44 unknowns (the 36 $x$
variables, plus the 8 parameters
$b_1$, $t_1$, $b_2$, $t_2$, $b_3$, $t_3$, $h$ and $o$).
These equations are as follows: 8 degree 2 equations (i.e., pick
one vertex from each of the 8 vertex classes, and write down the
degree 2 equation for each of these 8 vertices).
The ninth equation is the hyperplane for the comb inequality.
This results in a system $A \cdot x = b$ where
$A$ is 9 by 36,
$b$ is 9 by 1,
and the coefficients are polynomials in the 8 parameters.
(We multiply the last equation by 2 to eliminate the denominator
in those terms whose coefficient has the form ${k \choose 2}$.)
This system is as follows:
{\allowdisplaybreaks
\begin{align*}
	(b_1 - 1) \, x_{B_1,B_1}
	\,+\, t_1 \, x_{B_1,T_1}
	\,+\, b_2 \, x_{B_1,B_2}
	\,+\, t_2 \, x_{B_1,T_2}
	\,+\, b_3 \, x_{B_1,B_3} \\
	\,+\, t_3 \, x_{B_1,T_3}
	\,+\, h   \, x_{B_1,H}
	\,+\, o   \, x_{B_1,O}
	&= 2, \\
	b_1	  \, x_{B_1,T_1}
	\,+\, (t_1 - 1) \, x_{T_1,T_1}
	\,+\, b_2 \, x_{T_1,B_2}
	\,+\, t_2 \, x_{T_1,T_2}
	\,+\, b_3 \, x_{T_1,B_3} \\
	\,+\, t_3 \, x_{T_1,T_3}
	\,+\, h	  \, x_{T_1,H}
	\,+\, o	  \, x_{T_1,O}
	&= 2, \\
	b_1	  \, x_{B_1, B_2}
	\,+\, t_1 \, x_{T_1, B_2}
	\,+\, (b_2 - 1) \, x_{B_2,B_2}
	\,+\, t_2 \, x_{B_2,T_2}
	\,+\, b_3 \, x_{B_2,B_3} \\
	\,+\, t_3 \, x_{B_2,T_3}
	\,+\, h	  \, x_{B_2,H}
	\,+\, o	  \, x_{B_2,O}
	&= 2, \\
	b_1	  \, x_{B_1,T_2}
	\,+\, t_1 \, x_{T_1,T_2}
	\,+\, b_2 \, x_{B_2,T_2}
	\,+\, (t_2 - 1) \, x_{T_2,T_2}
	\,+\, b_3 \, x_{T_2,B_3} \\
	\,+\, t_3 \, x_{T_2,T_3}
	\,+\, h	  \, x_{T_2,H}
	\,+\, o	  \, x_{T_2,O}
	&= 2, \\
	b_1	  \, x_{B_1, B_3}
	\,+\, t_1 \, x_{T_1, B_3}
	\,+\, b_2 \, x_{B_2,B_3}
	\,+\, t_2 \, x_{T_2,B_3}
	\,+\, (b_3 - 1) \, x_{B_3,B_3} \\
	\,+\, t_3 \, x_{B_3,T_3}
	\,+\, h	  \, x_{B_3,H}
	\,+\, o	  \, x_{B_3,O}
	&= 2, \\
	b_1	  \, x_{B_1,T_3}
	\,+\, t_1 \, x_{T_1,T_3}
	\,+\, b_2 \, x_{B_2,T_3}
	\,+\, t_2 \, x_{T_2,T_3}
	\,+\, b_3 \, x_{B_3,T_3} \\
	\,+\, (t_3 - 1) \, x_{T_3,T_3}
	\,+\, h	  \, x_{T_3,H}
	\,+\, o	  \, x_{T_3,O}
	&= 2, \\
	b_1	  \, x_{B_1,H}
	\,+\, t_1 \, x_{T_1,H}
	\,+\, b_2 \, x_{B_2,H}
	\,+\, t_2 \, x_{T_2,H}
	\,+\, b_3 \, x_{B_3,H}
	\,+\, t_3 \, x_{T_3,H} \\
	\,+\, (h - 1) \, x_{H,H}
	\,+\, o	  \, x_{H,O}
	&= 2, \\
	b_1	  \, x_{B_1,O}
	\,+\, t_1 \, x_{T_1,O}
	\,+\, b_2 \, x_{B_2,O}
	\,+\, t_2 \, x_{T_2,O}
	\,+\, b_3 \, x_{B_3,O}
	\,+\, t_3 \, x_{T_3,O} \\
	\,+\, h	  \, x_{H,O}
	\,+\, (o - 1) \, x_{O,O}
	&= 2, \\
	2 \, b_1 \, (b_1 - 1) \, x_{B_1,B_1}
	+ 2 \, b_1 \, t_1 \, x_{B_1,T_1}
	+ 2 \, b_1 \, b_2 \, x_{B_1,B_2}
	+ 2 \, b_1 \, b_3 \, x_{B_1,B_3} \\
	+ 2 \, b_1 \, h   \, x_{B_1,H}
	+ t_1 \, (t_1 - 1) \, x_{T_1,T_1}
	\,+\, 2 \, b_2 \, (b_2 - 1) \, x_{B_2,B_2} \\
	\,+\, 2 \, b_2 \, t_2 \, x_{B_2,T_2}
	\,+\, 2 \, b_2 \, b_3 \, x_{B_2,B_3}
	\,+\, 2 \, b_2 \, h   \, x_{B_2,H}
	\,+\, t_2 \, (t_2 - 1) \, x_{T_2,T_2} \\
	\,+\, 2 \, b_3 \, (b_3 - 1) \, x_{B_3,B_3}
	\,+\, 2 \, b_3 \, t_3 \, x_{B_3,T_3}
	\,+\, 2 \, b_3 \, h   \, x_{B_3,H} \\
	\,+\, t_3 \, (t_3 - 1) \, x_{T_3,T_3}
	\,+\, h \, (h - 1)     \, x_{H,H}
	&= \\
	4 \, b_1
	\,+\, 2 \, t_1
	\,+\, 4 \, b_2
	\,+\, 2 \, t_2
	\,+\, 4 \, b_3
	\,+\, 2 \, t_3
	\,+\, 2 \, h
	\,- 10.
\end{align*}
}
Now use Gaussian elimination to solve for 9 of the $x$ variables in terms
of the remaining 27.
(Note that one should be careful to pick pivots that cannot be zero
for any valid assignment of the parameter values.)

Next, one forms the expression for distance squared
from the centroid
\begin{displaymath}
\sum_{U,V} \mathrm{NumEdges}(U,V) \, (x_{U,V} \,-\, {2 \over {n-1}})^2
\end{displaymath}
where $\mathrm{NumEdges}(U,V)$ is the number of edges in the class of
all edges from vertex class $U$ to vertex class $V$.
For example, $\mathrm{NumEdges}(T_1,T_1) = {{t_1} \choose 2}$,
and $\mathrm{NumEdges}(B_1,T_2) = b_1 \, t_2$.
Use the solution from the first system of equations to eliminate the 9
$x$ variables from this distance squared expression.
The resulting expression represents distance squared for a point $x$
that is constrained to lie in $F \cap \mathrm{aff}(P)$.
This distance squared is minimized when its gradient (with
respect to the 27 remaining $x$ variables) is identically zero.

This gives rise to a new 27 by 27 system of equations, the $i$-th row
obtained by taking the partial derivative of this distance squared
expression with
respect to the $i$-th $x$ variable (of the remaining 27 $x$ variables).
Clearing denominators results in coefficients that are polynomials in
the 8 parameters.
Solving this system is more challenging.
Once again, one should be careful to avoid choosing pivots that can be
zero for any valid combination of the 8 parameter values.
This is relatively easy to do until the last 4 pivots --- for which
all remaining pivot candidates are large polynomials, some of which
are fairly easy to zero, and others for which neither a simple zero
assignment nor proof of non-zeroness are evident.
Note that recognizing zero pivots in this context is indeed
Hilbert's tenth problem, which Matiyasevich proved has no algorithmic
solution~\cite{Matiyasevich1970},~\cite{Matiyasevich1993}.
We chose pivots for which we were unable to easily obtain zero
assignments.
The solution vector $x$ that we obtain does indeed satisfy all 9 of
the original equations, and identically zeros each of the 27 partial
derivatives.
Sustituting these $x$ values (each a rational function of the 8
parameter values) back into the distance squared formula yields the
final CD squared formula given in
Tables~\ref{tab:tsp-centroid-distances}
and~\ref{tab:3comb-dist-subexprs}.
Both the closest point $x$ (as given by the $x$ variables) and the
final CD formula yield answers that agree exactly with
numeric values obtained independently via quadratic programming for
various small $n$, as described in
Section~\ref{sec:validation}
below.

The careful reader will notice that when $h=0$ or $o=0$, there are no
vertices in these sets, nor are there any edges incident to them.
The corresponding degree-2 equations therefore reduce to
$0 ~=~ 2$,
which certainly calls into question the validity of the general
solution in three special cases
\begin{eqnarray*}
 \text{$h = 0$ and $o > 0$}, \\
 \text{$h > 0$ and $o = 0$}, \\
 \text{$h = 0$ and $o = 0$}.
\end{eqnarray*}
We have verified in each of these three cases that re-solving the
restricted problem (e.g., by omitting the offending vertices, together
with their degree-2 equations and incident edges) yields a solution
that is identical (both for distance squared, and the closest point
$x$) to simply letting $h$ and/or $o$ be zero in the general solution.
The general solution can therefore be safely used without regard to
these concerns.

\section{Spanning Tree in Hypergraph Polytope: $STHGP(n)$}
\label{sec:sthgp-polytope}

Let $n \in \Z$, $n \ge 2$.  Let $H = (V,E)$ be the complete
hypergraph on $n$ vertices.
Let $m = |E| = 2^n-n-1$.
($E$ contains all subsets of $V$ except those of cardinality 0 and 1.)
We work in the space of $\R^m$ with incidence vectors that represent
subsets of $E$.
Let $ST(n)$ be the set of all incidence vectors $x \in \B^m$ that
correspond to spanning trees of $H$
(i.e., each incidence vector $x \in ST(n)$ corresponds to an $E'$ such
that $E' \subseteq E$
and the subhypergraph $H'=(V,E')$ is a tree).
We define the {\em Spanning Tree in Hypergraph Polytope}
\begin{displaymath}
\STHGP{n} \,=\, \CHull{ST(n)}
\end{displaymath}

\noindent
The subtour relaxation of $\STHGP{n}$ is defined as
\begin{eqnarray}
\label{eq:sthgp-total-degree-equation}
\sum_{e \in E} (|e| \,-\, 1) \, x_e &=& |V| \,-\, 1, \\
\label{eq:sthgp-subtour-facet-def}
	\sum_{e \in E} \max(|e \cap S| - 1, 0) \, x_e &\le& |S| - 1
	~~~~~~~
	\hbox{for all $S \subset V$, $|S| \ge 2$}, \\
\label{eq:sthgp-nonneg-facet-def}
	x_e &\ge& 0
	~~~~~~~~~~~~~~
	\hbox{for all $e \in E$}.
\end{eqnarray}

\noindent
where~(\ref{eq:sthgp-subtour-facet-def}) are the subtour inequalities,
and (\ref{eq:sthgp-nonneg-facet-def}) are the non-negativity constraints.
From~(\cite{Warme98}, pages 42--44) we know
that:
(\ref{eq:sthgp-total-degree-equation}) is the only equation
satisfied by all spanning trees, and therefore defines the affine hull
of $\STHGP{n}$;
(\ref{eq:sthgp-subtour-facet-def})
and~(\ref{eq:sthgp-nonneg-facet-def}) define facets of $\STHGP{n}$.

Equation~(\ref{eq:sthgp-total-degree-equation})
can be written more compactly as $x(V) = |V| \,-\, 1$.
Subtour inequality~(\ref{eq:sthgp-subtour-facet-def}) can similarly be
written $x(S) \le |S| \,-\, 1$.
By subtracting~(\ref{eq:sthgp-subtour-facet-def})
from~(\ref{eq:sthgp-total-degree-equation}) and letting $T = V-S$, it
can be shown that
\begin{displaymath}
	x(T) \,+\, x(T:T-S) ~\ge~ |T|,
\end{displaymath}
an alternate form for the subtour inequality that we call an
{\em anti-subtour}.
Subtour $S$ places an upper bound on the weight of edges within $S$,
whereas anti-subtour $T$ places a lower bound on the weight of edges
within and crossing out of $T$.
One can equivalently use whichever form has the fewest non-zeros.

We start with several preliminary results.
Let $X_{\lambda}$ be a Poisson random variable with mean $\lambda$,
e.g., $P( X_{\lambda}=k) ~=~ {{e^{-\lambda} \, \lambda^k} \over {k!}}$
for all integers $k \ge 0$.  The moment generating function for
$X_{\lambda}$ is well known to be (see, e.g.,
\cite{Blake1987}):
\begin{equation}
\label{eq:Poisson-mgf}
	e^{\lambda \, (e^z \,-\, 1)}
	~=~
	\sum_{k \ge 0} E[X_{\lambda}^k] \, {{z^k} \over {k!}}.
\end{equation}

\noindent
The following result about Stirling numbers of the second kind is well
known (see, e.g., \cite{GrahamKnuthPatashnik}):
\begin{equation}
\label{eq:stirling-power-egf}
	e^{\lambda \, (e^z \,-\, 1)}
	~=~
	\sum_{i,k \ge 0} {k \brace i} \, \lambda^i \, {{z^k} \over {k!}}
\end{equation}
which is another interpretation of the moment generating function
for Poisson random variables (\ref{eq:Poisson-mgf}).  From this we
conclude that
\begin{equation}
\label{eq:poisson-moment-as-stirling-sum}
E[X_{\lambda}^k] ~=~ \sum_{i=0}^k {k \brace i} \, \lambda^i.
\end{equation}

\noindent
Equation (\ref{eq:poisson-moment-as-stirling-sum}) allows the Stirling
sums that appear in the results for the spanning tree in hypergraph
polytope to be more conveniently and compactly written using moments
of Poisson random variables.
We assume $X_{\lambda}$ is such a Poisson random variable and use the
more compact moment notation in lieu of the equivalent Stirling sum
throughout.

The exponential generating function for the Bell numbers is well known
to be
\begin{equation}
\label{eq:bell-egf}
e^{e^z \,-\, 1} ~=~ 1 \,+\, \sum_{n \ge 1} \Bell(n) \, {{z^n} \over {n!}}
\end{equation}
(see, e.g., \cite{Stanley},
\cite{Wilf}).

Both indicators (EPR and CD) require that
we count the number of extreme points of $\STHGP{n}$, or in other
words the number of spanning trees in the complete hypergraph on $n$
labeled vertices.
(By labeled, we mean that each vertex is uniquely
distinguishable.
Counting trees over unlabeled vertices means counting only
distinct tree shapes, since there is no way to tell one vertex
from another except by differences in degree, edge topology, asymmetry
or other contextual factors.)

\begin{theorem}[Warme~\cite{Warme98}, page 57]
\label{th:ntrees}
Let $H = (V,E)$ be a complete hypergraph, $n = |V|$, $n \ge 1$.
The number
$t_n$ of distinct spanning trees in $H$ is
\begin{eqnarray}
\label{eq:ntrees-stirling-sum}
t_n	&=& \sum_{i=0}^{n-1} {{n-1} \brace i} \, n^{i-1} \\
\label{eq:ntrees-poisson-moment}
	&=& {1 \over n} \, E \lbrack X_n^{n-1} \rbrack,
\end{eqnarray}
where $X_n$ is a Poisson random variable with mean $n$.
\end{theorem}

\begin{remark}
\label{rem:summation-index-is-num-edges}
The summation index $i$ in~(\ref{eq:ntrees-stirling-sum}) is the
number of edges in the tree.
\end{remark}

\begin{remark}
The $i \,=\, n-1$ term of~(\ref{eq:ntrees-stirling-sum}) is the number
of trees having $n-1$ edges, thereby giving the classic result of
$n^{n-2}$ trees in a conventional graph --- a result known as Caley's
formula~\cite{Cayley}.
\end{remark}

This author first obtained the result in Stirling sum
form~(\ref{eq:ntrees-stirling-sum}).
Lang Withers later noticed the Poisson moment
form~(\ref{eq:ntrees-poisson-moment}).
There is currently no intuitive explanation for why spanning trees in
hypergraphs should have any connection with moments of Poisson random
variables.

The proof of Theorem~\ref{th:ntrees} resulted from collaboration of
this author with W.~D.~Smith, and first appeared in
(\cite{Warme98}, pages 54--59).
(Smith proved Theorems~\ref{th:tree-recurrence} and~\ref{th:tree-egf},
which Warme then used to prove Theorem~\ref{th:ntrees}.)
We include all three of these theorems here, together with their
proofs to make this paper self-contained, and to illustrate all of the
combinatorial techniques used.

For $n \ge 1$, let $h_n$ be the number of {\em rooted} spanning trees
of the complete hypergraph on $n$ labeled vertices.
A rooted tree is a tree in which one particular vertex is identified
as being the {\em root}.
The desired result for unrooted trees is then $t_n = h_n / n$.
Considering rooted trees breaks up the symmetry of the problem and
avoids various automorphisms that would otherwise result.

W.~D.~Smith obtained the following recurrence and generating function
for $h_n$:

\begin{theorem}[W. D. Smith]
\label{th:tree-recurrence}
Let $h_n$ be the number of rooted trees spanning the complete
hypergraph on $n$ labeled vertices.
Then $h_1 = 1$, and for $n > 1$
\begin{equation}
\label{eq:hypertree-recurrence}
h_n = n \sum_{k > 0} {\Bell(k) \over {k!}}
	\sum_{\MyAtop{a_j > 0}%
		      {\sum_{j=1}^k a_j = n-1}}
	{{n-1} \choose {a_1,a_2,\ldots,a_k}}
	\prod_{j=1}^k h_{a_j}.
\end{equation}
\end{theorem}

\begin{MyProof}
The base case is obvious, so assume $n > 1$.
Select a unique root vertex (there are $n$ possible choices).
Now delete the root vertex and every hyperedge incident to the root.
All that remains are the individual subtrees of the root node,
containing a total of $n-1$ vertices.  Each of these subtrees is
itself a rooted tree, the root vertex being the one that was
incident to a deleted edge.  Suppose there are $k$ of these
rooted subtrees.
The vector $a_1,a_2,\ldots,a_k$ indicates how many vertices are in
each of the $k$ subtrees.  We divide by $k!$ since the particular
ordering of the subtrees does not matter.  For each such vector
there are
\begin{displaymath}
{{n-1} \choose {a_1,a_2,\ldots,a_k}}
\end{displaymath}
ways of partitioning the $n-1$ vertices into $k$ non-empty subsets of
sizes $a_1$, $a_2,$ $\ldots,$ $a_k$.
For each subset $j$ of $a_j$ vertices,
there are $h_{a_j}$ distinct rooted subtrees.
Each of the $\Bell(k)$ partitions of the $k$ subtrees
represents a distinct way of hooking the subtrees to the root
using edges.  Let $S_1,S_2,\ldots,S_j$ be such a partition.  Then
the $k$ subtrees are connected to the root using $j$ edges.
The edge for $S_i$ consists of {\em the} root together with the
root vertices of each subtree in $S_i$.
\qed
\end{MyProof}

\begin{remark}
Replacing $\Bell(k)$ with $1$ in~(\ref{eq:hypertree-recurrence}) gives
a recurrence for rooted trees in conventional graphs.
\end{remark}

\begin{theorem}[W. D. Smith]
\label{th:tree-egf}
Let
\begin{equation}
\label{eq:H-egf-def}
H(z) = \sum_{n \ge 1} h_n {z^n \over n!}
\end{equation}
be the exponential generating function for $h_n$.
Then
\begin{equation}
\label{eq:hyper-tree-egf}
H(z) = z \,\, e^{e^{H(z)} - 1}.
\end{equation}
\end{theorem}

\begin{MyProof}
Taking the $k$-th power of~(\ref{eq:H-egf-def}) using the multinomial
expansion yields:
\begin{displaymath}
H(z)^k ~=~
\sum_{n \ge 1}
\left (
	\sum_{\MyAtop{a_j > 0}{\sum_{j=1}^k a_j = n}}
	{n \choose {a_1,a_2,\ldots,a_k}}
	\prod_{j=1}^k h_{a_j}
\right )
{{z^n} \over {n!}}
\end{displaymath}
from which we conclude that
\begin{equation}
\label{eq:multinomial-as-coeff-of-H(z)^k}
	\sum_{\MyAtop{a_j > 0}{\sum_{j=1}^k a_j = n-1}}
	{{n-1} \choose {a_1,a_2,\ldots,a_k}}
	\prod_{j=1}^k h_{a_j}
=	\left [ {{z^{n-1}} \over {(n-1)!}} \right ] \,\, H(z)^k.
\end{equation}
Therefore, if $n > 1$ we can substitute
~(\ref{eq:multinomial-as-coeff-of-H(z)^k}) into
~(\ref{eq:hypertree-recurrence}) to obtain
\begin{eqnarray}
\nonumber
h_n &=& \left [ {z^{n-1} \over (n-1)!} \right ] \,\,
		n \sum_{k > 0} \Bell(k) {H(z)^k \over {k!}} \\
\label{eq:h-recurrence1}
&=&	n! \, \left [ z^{n-1} \right ] \,\,
		\sum_{k > 0} \Bell(k) {H(z)^k \over {k!}}.
\end{eqnarray}
Recall from~(\ref{eq:bell-egf}) that
\begin{displaymath}
1 + \sum_{k \ge 1} \Bell(k) {z^k \over k!} = e^{e^z - 1}.
\end{displaymath}
Substituting into~(\ref{eq:h-recurrence1}) yields:
\begin{displaymath}
h_n
=	n! \, \lbrack z^{n-1} \rbrack \,\,
	(1 + \sum_{k \ge 1} \Bell(k) {H(z)^k \over k!})
=	n! \, \lbrack z^{n-1} \rbrack \,\, e^{e^{H(z)} - 1}
=	n! \, \lbrack z^n \rbrack \,\, z \,\, e^{e^{H(z)} - 1}
\end{displaymath}
which happens to hold at $n = 1$ as well as for $n > 1$.  Therefore
\begin{displaymath}
H(z) = \sum_{n > 0} h_n \,\, {z^n \over n!} = z \,\, e^{e^{H(z)} - 1}.
\tag*{\qed}
\end{displaymath}
\end{MyProof}

\noindent
To obtain a closed-form for $h_n$, we employ the Lagrange inversion
formula:

\begin{MyLemma}[Lagrange inversion formula]
\label{lem:lif}
Let $f(u)$ and $\theta(u)$ be formal power series in $u$, such that
$\theta(0) = 1$.  Then there is a unique formal power series
$u = u(t)$
(about $t = 0$) satisfying
\begin{displaymath}
	u = t \,\, \theta(u).
\end{displaymath}
In addition, the value $f(u(t))$ of $f$ at that root
$u = u(t)$, when expanded in a power series in $t$ about $t=0$,
satisfies
\begin{displaymath}
\lbrack t^n \rbrack \,\, \lbrace f(u(t)) \rbrace
~=~ {1 \over n} \,\, \lbrack u^{n-1} \rbrack \,\,
	\lbrace f'(u) \, \theta(u)^n \rbrace.
\end{displaymath}
\end{MyLemma}

\noindent
Proofs can be found
in~\cite{Wilf}, \cite{SuryaWarnke2023}, \cite{Gessel2016}.
We now move on to the main proof.

\begin{MyProof}[of Theorem~\ref{th:ntrees}]
Apply the Lagrange inversion formula (Lemma~\ref{lem:lif}) to
(\ref{eq:hyper-tree-egf})
with $f(u)=u$ and $\theta(u) = e^{e^u - 1}$:
\begin{eqnarray}
\nonumber
{h_n \over n!}
&=&	\lbrack z^n \rbrack \,\, H(z) \\
\nonumber
&=&	{1 \over n} \,\, \lbrack u^{n-1} \rbrack \,\, \theta(u)^n \\
\label{eq:hyper-tree-intermediate-form}
&=&	{1 \over n} \,\, \lbrack u^{n-1} \rbrack \,\, e^{n (e^u - 1)}.
\end{eqnarray}
From (\ref{eq:hyper-tree-intermediate-form}) we can
use~(\ref{eq:stirling-power-egf})
to go toward Stirling numbers of the second kind:
\begin{eqnarray*}
{h_n \over n!}
&=&	{1 \over n} \,\, \lbrack u^{n-1} \rbrack \,\,
	\sum_{k,i \ge 0} {k \brace i} \, n^i \, {{u^k} \over {k!}} \\
&=&	{1 \over n} \sum_{i \ge 0} {n-1 \brace i} \, n^i
	\, {1 \over {(n-1)!}} \\
&=&	{1 \over {n!}} \sum_{i=0}^{n-1} {{n-1} \brace i} \, n^i \\
\Longrightarrow \\
h_n	&=&  \sum_{i=0}^{n-1} {{n-1} \brace i} \, n^i.
\end{eqnarray*}
Alternatively, from~(\ref{eq:hyper-tree-intermediate-form})
we can use~(\ref{eq:Poisson-mgf})
to go toward toward Poisson moments:
\begin{eqnarray*}
{h_n \over n!}
&=&	{1 \over n} \,\, \lbrack u^{n-1} \rbrack \,\,
	\sum_{k \ge 0} E[X_n^k] \, {{u^k} \over {k!}} \\
&=&	{1 \over n} \, E[X_n^{n-1}] \, {1 \over {(n-1)!}} \\
\Longrightarrow \\
h_n	&=&	 E[X_n^{n-1}].
\end{eqnarray*}
Since the number of unrooted trees $t_n$ for $n$ vertices is
$t_n = h_n / n$, we therefore have:
\begin{align*}
t_n	&= \sum_{i=0}^{n-1} {{n-1} \brace i} \, n^{i-1} \\
	&= {1 \over n} \, E[X_n^{n-1}].
\tag*{\qed}
\end{align*}
\end{MyProof}

\begin{theorem}
\label{th:sthgp-trees-containing-edge-e}
Let $e \in E$ and $k = |e|$.
The number of trees in $H =(V,E)$ that
contain edge $e$ is $g(n,k)$, where
\begin{displaymath}
	g(n,k) ~=~ {k \over n} \, E[X_n^{n-k}].
\end{displaymath}
\end{theorem}

\begin{MyProof}
Remove edge $e$ and the $k$ vertices incident to $e$, leaving
$n-k$ vertices that must be partitioned into $k$ (possibly empty)
subsets, each of which becomes a subtree rooted at its corresponding
vertex of $e$.
Let $b_1,\, b_2,\, \ldots,\, b_k \ge 0$ be the sizes of these subsets,
where $\sum_{i=1}^k b_i = n-k$.
There are ${{n-k} \choose {b_1,\, b_2,\, \ldots,\, b_k}}$ ways of
partitioning the $n-k$ vertices into these subsets.
For each subset $i$ of size $b_i$ there are
$h_{b_i + 1} / (b_i + 1)$ (unrooted) spanning trees over the
corresponding $i$-subset (including the corresponding vertex from $e$).
Let $a_i = b_i + 1$ for $i=1,2,\ldots,k$.
We therefore have
{\allowdisplaybreaks%
\begin{eqnarray*}
g(n,k)
&=&	\sum_{\MyAtop{b_j \ge 0}{\sum_{j=1}^k b_j = n-k}}
	{{n-k} \choose {b_1,\,b_2,\,\ldots,\,b_k}}
	\prod_{j=1}^k {{h_{b_j + 1}} \over {b_j + 1}} \\
&=&	\sum_{\MyAtop{b_j \ge 0}{\sum_{j=1}^k b_j = n-k}}
	{{(n-k)!} \over {b_1!\,b_2!\,\ldots\,b_k!}}
	\prod_{j=1}^k {{h_{b_j + 1}} \over {b_j + 1}} \\
&=&	\sum_{\MyAtop{b_j \ge 0}{\sum_{j=1}^k b_j = n-k}}
	{{(n-k)!} \over {(b_1+1)!\,(b_2+1)!\,\ldots\,(b_k+1)!}}
	\prod_{j=1}^k h_{b_j + 1} \\
&=&	\sum_{\MyAtop{a_j > 0}{\sum_{j=1}^k a_j = n}}
	{{(n-k)!} \over {a_1!\,a_2!\,\ldots\,a_k!}}
	\prod_{j=1}^k h_{a_j} \\
&=&	(n-k)! \,
	\sum_{\MyAtop{a_j > 0}{\sum_{j=1}^k a_j = n}}
	{{n!} \over {a_1!\,a_2!\,\ldots\,a_k!}}
	\, \left [ \prod_{j=1}^k h_{a_j} \right ]
	\, {1 \over {n!}} \\
&=&	(n-k)! \,
	\sum_{\MyAtop{a_j > 0}{\sum_{j=1}^k a_j = n}}
	{n \choose {a_1,\,a_2,\,\ldots,\,a_k}}
	\, \left [ \prod_{j=1}^k h_{a_j} \right ]
	\, {1 \over {n!}} \\
&=&	(n-k)! \, [z^n] \, H(z)^k.
\end{eqnarray*}%
}%
We now apply the Lagrange inversion formula with
$f(u) = u^k$ and $\theta(u) = e^{e^u - 1}$:
{\allowdisplaybreaks%
\begin{align*}
g(n,k)
&=	(n-k)! \, {1 \over n} \, [u^{n-1}] \,
	\lbrace k \, u^{k-1} \, (e^{e^u - 1})^n \rbrace \\
&=	{k \over n} \, (n-k)! \, [u^{n-k}] \,
	\lbrace e^{n(e^u - 1)} \rbrace \\
&=	{k \over n} \, (n-k)! \, [u^{n-k}] \,
	\sum_{j \ge 0} E[X_n^j] \, {{u^j} \over {j!}} \\
&=	{k \over n} \, (n-k)! \, E[X_n^{n-k}] \, {1 \over {(n-k)!}} \\
&=	{k \over n} \,  E[X_n^{n-k}].
\tag*{\qed}
\end{align*}%
}%
\end{MyProof}

\begin{theorem}
\label{th:ntrees=g(n,1)}
The number of distinct spanning trees in $H$ is $g(n,1)$.
\end{theorem}

\begin{MyProof}
By theorem~\ref{th:sthgp-trees-containing-edge-e},
and using (\ref{eq:poisson-moment-as-stirling-sum}):
\begin{displaymath}
g(n,1) = {1 \over n} E[X_n^{n-1}]
= {1 \over n} \sum_{i=0}^{n-1} {{n-1} \brace i} \, n^i
= \sum_{i=0}^{n-1} {{n-1} \brace i} \, n^{i-1}
\end{displaymath}
which matches both results for $t_n$ from theorem \ref{th:ntrees}.
\qed
\end{MyProof}

We now have all of the tools needed to derive EPR
for two facet classes of $\STHGP{n}$.

\subsection{EPR for STHGP Non-Negativity Inequalities}

\begin{theorem}
Let $n \ge 2$, $G = (V,E)$ be the complete hypergraph with $|V|=n$,
$e \in E$ and $k = |e|$.
Then EPR for the non-negativity constraint
\begin{displaymath}
	x_e \,\ge\, 0
\end{displaymath}
of $\STHGP{n}$ is
\begin{displaymath}
	r(n,k) ~=~ 1 \,-\, {{k \, E[X_n^{n-k}]} \over {E[X_n^{n-1}]}}.
\end{displaymath}
\end{theorem}

\begin{MyProof}
The set of extreme points incident to this
inequality are precisely those spanning trees that do {\em not}
contain edge $e$.
By Theorems~\ref{th:sthgp-trees-containing-edge-e}
and~\ref{th:ntrees=g(n,1)},
there are $g(n,k)$ trees that contain edge $e$,
and $g(n,1)$ trees total.
The number of trees that do {\em not} contain edge $e$ is therefore
$g(n,1) \,-\, g(n,k)$.
The fraction of all trees incident to $x_e \ge 0$ is therefore
\begin{align*}
r(n,k)	&=	{{g(n,1) \,-\, g(n,k)} \over {g(n,1)}} \\
	&=	1 \,-\, {{k \, E[X_n^{n-k}]} \over {E[X_n^{n-1}]}}.
\tag*{\qed}
\end{align*}
\end{MyProof}

\subsection{EPR for STHGP Subtour Inequalities}

\begin{theorem}
\label{th:sthgp-subtour-epr}
Let $S \subset V$ and $k = |S|$ such that $2 \le k < n$.
Then EPR for the subtour elimination constraint
\begin{equation}
\label{eq:sthgp-subtour-epr}
	\sum_{e \in E} \max(|e \cap S| \,-\, 1,\, 0) \, x_e
		\,\le\, |S| \,-\, 1
\end{equation}
of $\STHGP{n}$ is
\begin{displaymath}
r(n,k) ~=~
{{n \, \sum_{i=0}^{k-1} {{k-1} \brace i} \, k^i \,
    \left \lbrack
	\sum_{j=1}^{n-k} {{n-k} \choose j}
			 \, j
			 \, E[X_{n-k}^{n-k-j}]
			 \sum_{p=0}^j {j \choose p} E[X_k^p] \, i^{j-p}
    \right \rbrack
  }
  \over
  {k \, (n-k)\,E[X_n^{n-1}]}}.
\end{displaymath}
\end{theorem}
\MyBeginFig
 \begin{center}
  \includegraphics[width=4.5in,clip=]{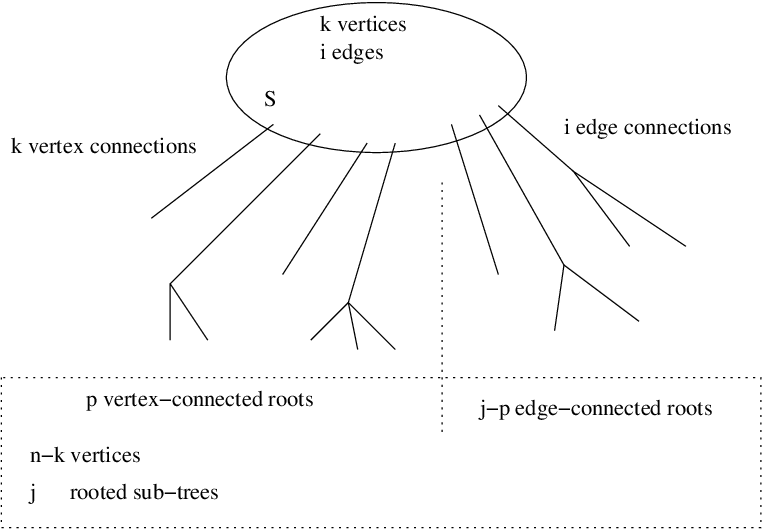}
  \caption{Trees incident to $k$ vertex $\STHGP{n}$ subtour}%
  \label{fig:trees-incident-to-k-vertex-subtour}
 \end{center}
\MyEndFig

\begin{MyProof}
A tree is incident to subtour $S$ when the subhypergraph induced by
$S$ is itself a tree.
We therefore count the trees incident to subtour $S$ by starting with
each tree over the $k$ vertices of $S$, and counting the number of
ways of connecting the other $n-k$ vertices to it in ways that yield a
tree over all $n$ vertices.
% We call the $n-k$ vertices outside of $S$ ``free'' vertices, since
% they need to be connected to the existing tree structure within
% $S$.

Refer to Figure~\ref{fig:trees-incident-to-k-vertex-subtour}.
The number of spanning trees within the $k$ vertices of $S$ is
\begin{displaymath}
	{1 \over k} \sum_{i=0}^{k-1} {{k - 1} \brace i} \, k^i
\end{displaymath}
and the $i$-th term of this sum counts the number of said trees having
exactly $i$ edges.
We partition the remaining $n-k$ vertices into $j$ subtrees, with
$p$ of these connecting to the $k$ vertices within $S$, and the
remaining $j-p$ of these subtrees connecting to the $i$ edges within
$S$.
There are ${j \choose p}$ distinct ways of partitioning the $j$
subtrees into $p$ vertex-connected and $j-p$ edge-connected types.
The expression counting the total number of trees incident to subtour
$S$ therefore has the following structure:
\begin{equation}
\label{eq:trees-over-k-vertex-subtour-template}
f(n,k) ~=~
	{1 \over k}
	\sum_{i=0}^{k-1} {{k-1} \brace i} \, k^i \,
	    \left \lbrack
		\sum_{j=1}^{n-k}
			S(n-k, j) \,
			\sum_{p=0}^j
				{j \choose p}
				V(k,p) \, E(i,j-p)
	    \right \rbrack
\end{equation}
where:
\begin{itemize}
 \item $S(n-k,j)$ is the number of ways of partitioning the remaining
   $n-k$ vertices into an unordered set of $j$ rooted subtrees;
 \item $V(k,p)$ is the number of ways of using edges to connect the
   $k$ vertices of $S$ with $p$ of these rooted subtrees;
 \item $E(i,j-p)$ is the number of ways of augmenting each of the $i$
   edges within the current tree over $S$ to connect to the remaining
   $j-p$ rooted subtrees.  Each of the $i$ edges includes zero or more
   of the subtree root vertices, and each of the $j-p$ subtree
   root vertices is included in exactly one of the $i$ edges.
\end{itemize}
It remains to calculate $S(n-k,j)$, $V(k,p)$, and $E(i,j-p)$.

$S(n-k,j)$ is equivalent to the number of trees over $n-k$ vertices
that contain a particular edge $e$ of cardinality $j$ (and then
removing edge $e$ leaving $j$ distinct rooted subtrees).
There are ${{n-k} \choose j}$ possible choices of such an edge $e$.
This yields
\begin{equation}
\nonumber
\label{eq:S-of-n-k-and-j}
	S(n-k,j) ~=~ {{n-k} \choose j} \, g(n-k, j)
		 ~=~ {{n-k} \choose j} \,
		     {j \over {n-k}} \,
		     E[X_{n-k}^{n-k-j}].
\end{equation}
We can express $V(k,p)$ as
\begin{equation}
\label{eq:V-of-k-and-p-sum-over-vectors-form}
	V(k,p) ~=~ \sum_{\MyAtop{(a_1,\,a_2,\, \ldots,\, a_k) \,\ge\, 0}%
			 {\sum_{q=1}^k a_q = p}}
			{p \choose {a_1,\, a_2,\, \ldots,\, a_k}} \,
			\prod_{q=1}^k \Bell(a_q).
\end{equation}
Vector element $a_q$ indicates how many of the $p$ rooted subtrees
are connected to the $q$-th vertex of $S$.
This partitions those $p$ roots into $k$ bins (each bin $q$
containing the indicated number $a_q$ of roots), and can be done in
${p \choose {a_1,\, a_2,\, \ldots,\, a_k}}$
distinct ways.
For each vertex $q$ in $S$, there are $\Bell(a_q)$ distinct ways of
connecting vertex $q$ of $S$ to the $a_q$ root vertices of their
respective subtrees.
The exponential generating function for $V(k,p)$ as presented
in~(\ref{eq:V-of-k-and-p-sum-over-vectors-form})
is seen to be
\begin{equation}
\label{eq:V-of-k-and-p-egf-template}
	V(k,p) ~=~ \left \lbrack {{z^p} \over {p!}} \right \rbrack \,
			 (f(z))^k
\end{equation}
where $f(z)$ is the exponential generating function for the
Bell numbers~(\ref{eq:bell-egf}).
From~(\ref{eq:V-of-k-and-p-egf-template}) and~(\ref{eq:bell-egf})
we conclude that
{\allowdisplaybreaks%
\begin{eqnarray}
\nonumber
V(k,p)	&=& \left \lbrack {{z^p} \over {p!}} \right \rbrack \,
			 (f(z))^k \\
\nonumber
	&=& \left \lbrack {{z^p} \over {p!}} \right \rbrack \,
		(e^{e^z \,-\, 1})^k \\
\nonumber
	&=& \left \lbrack {{z^p} \over {p!}} \right \rbrack \,
		(e^{k \,(e^z \,-\, 1)}) \\
\label{eq:V-of-k-and-p}
	&=& E[X_k^p].
\end{eqnarray}%
}%

\noindent
We can express $E(i,j-p)$ as
\begin{displaymath}
	E(i,j-p) ~=~ \sum_{\MyAtop{(b_1,\,b_2,\, \ldots,\, b_i) \,\ge\, 0}%
			 {\sum_{q=1}^i b_q = j-p}}
			{{j-p} \choose {b_1,\, b_2,\, \ldots,\, b_i}}.
\end{displaymath}
(There is no product after the multinomial coefficient here because
for each edge $q$, $1 \le q \le i$, we need only choose the $b_q$
subtree roots to include within edge $q$.)
This expression is the multinomial expansion of
$(1 \,+\, 1 \,+\, \cdots \,+\, 1)^{j-p}$
where there are $i$ copies of 1 in the sum being exponentiated.
This yields
\begin{equation}
\label{eq:E-of-i-and-j-p}
	E(i,j-p) ~=~ i^{j-p}.
\end{equation}
We combine~(\ref{eq:trees-over-k-vertex-subtour-template}),
(\ref{eq:S-of-n-k-and-j}), (\ref{eq:V-of-k-and-p}),
(\ref{eq:E-of-i-and-j-p}) and divide by the total number of trees
$E[X_n^{n-1}]/n$ to obtain EPR
\begin{displaymath}
r(n,k) =
{{n \sum_{i=0}^{k-1} {{k-1} \brace i} \, k^i
	\sum_{j=1}^{n-k} {{n-k} \choose j}
			 \, j
			 \, E[X_{n-k}^{n-k-j}]
			 \sum_{p=0}^j {j \choose p} E[X_k^p] \, i^{j-p}
  }
  \over
  {k \, (n-k) \, E[X_n^{n-1}]}}.
\tag*{\qed}
\end{displaymath}
\end{MyProof}

\subsection{Alternate Form for EPR of STHGP Subtour Inequalities}

Although expression~(\ref{eq:sthgp-subtour-epr}) is a closed-form
(assuming you consider binomial coefficients and Stirling numbers to
be closed-forms), it is expensive to evaluate computationally because
of deeply nested summation loops.
(It is tractable for $n=50$, but totally impractical by $n \ge 100$.)
We now present an alternate form that is more computationally
efficient, but not a closed-form (it uses a function $E(i,j)$ defined
via a non-standard recurrence).

Let $S \subset V$ and $k = |S|$ such that $2 \le k < n$.
Note that a tree is incident to subtour $S$ if it contains the
maximum amount of ``edge material'' within $S$.  This happens when the
subhypergraph induced by $S$ is a tree.  We therefore count the number
of trees incident to subtour $S$ by starting with each tree over the
$k$ vertices of $S$, and counting the number of ways of connecting the
other $n-k$ vertices to it in ways that yield a tree over all $n$
vertices.  We call the $n-k$ vertices outside of $S$ ``free''
vertices, since they need to be connected to the existing tree
structure within $S$.

There are two ways that a subtree of free vertices can connect to an
existing tree structure within $S$:
{\em vertex} attachment and {\em edge} attachment.
Let $F \subseteq V \setminus S$ be a set of free vertices and $j=|F|$.
Let $v \in S$.
Subset $F$ is {\em vertex} attached to vertex $v$ simply by
constructing any tree over vertices $F \cup \lbrace v \rbrace$, which
can be done in $g(j+1,1) ~=~ E[X_{j+1}^j]/(j+1)$ different ways.

Alternatively, let $e$ be an existing edge in a tree over $S$.
If $j=0$, then $F=\emptyset$ and there is exactly one way that $F$ can
be attached to edge $e$ (edge $e$ remains unchanged).
Otherwise $j > 0$, and
subset $F$ can be {\em edge} attached to edge $e$ in two distinct ways:
(1) construct a tree for $F$ and include one of its $j$ vertices within
edge $e$ as an additional vertex outside of $S$;
or (2) construct a tree for $F$ and merge one of its edges with edge
$e$.

Define $V(j,k)$ to be the number of ways of using vertex attachment to
connect $j$ free vertices to the $k$ vertices of $S$.
\begin{lemma}
	V(j,k) ~=~ g(j+k,k)
\end{lemma}
\begin{proof}
Vertex attachment does not depend upon the structure of our initial
tree over $S$.
If we simply assume the tree over $S$ is formed by the single edge
$e = S$, then Theorem~\ref{th:sthgp-trees-containing-edge-e}
tells us there are precisely $g(j+k,k)$ trees over the $j$ free
vertices plus the $k$ vertices of $S$ that contain the $k$-vertex edge
$e=S$.
\qed
\end{proof}

Define $E(i,j)$ to be the number of ways of using edge attachment to
connect $j$ free vertices to $i$ existing edges within a tree over
$S$.
\begin{lemma}
$E(i,j)$ satisfies the following recurrence:
\begin{eqnarray}
\label{eq:E(i,0)}
	E(i,0) &=& 1	~~~~\hbox{for all $i \ge 0$} \\
\label{eq:E(0,j)}
	E(0,j) &=& 0	~~~~\hbox{for all $j > 0$} \\
\label{eq:E(1,j)}
	E(1,j) &=& {1 \over j} E[X_j^j]	~~~~\hbox{for all $j > 0$} \\
\label{eq:E(i,j)}
	E(i,j) &=&
		\sum_{m=0}^j {j \choose m} E(i-1,j-m) \, E(1,m)
	~~~~\hbox{for all $j>0$, $i>1$}
\end{eqnarray}
\end{lemma}
\begin{proof}
Base cases~(\ref{eq:E(i,0)}) and~(\ref{eq:E(0,j)}) are clearly true.
For $E(1,j)$ we need to count the number of ways of attaching $j$ free
vertices to a single existing edge $e \subseteq S$.
Method (1) can be done in
$h_j = E[X_j^{j-1}]$ ways, since there are $h_j$ rooted trees over $j$
vertices and we simply include the root vertex within edge $e$.
Method (2) is a bit more involved.
We recall Remark~\ref{rem:summation-index-is-num-edges}, and note
that in a subtree over $F$ consisting of $i$ edges, there are $i$
distinct choices of edge to merge with edge $e$.  This means the
number of ways in which method (2) can be accomplished is
\begin{eqnarray*}
&&	{1 \over j} \sum_{i=0}^{j-1}
		\, i \, {{j-1} \brace i} j^i \\
&=&	{1 \over j} \sum_{i=1}^{j-1}
		\, i \, {{j-1} \brace i} j^i \\
&=&	{1 \over j} \sum_{i=1}^{j-1}
		\left (
			{j \brace i} \,-\, {{j-1} \brace {i-1}}
		\right ) j^i \\
&=&	\left \lbrack
	{1 \over j} \sum_{i=1}^{j-1}
		{j \brace i} j^i
	\right \rbrack
	~-~
	\left \lbrack
	{1 \over j} \sum_{i=1}^{j-1}
		{{j-1} \brace {i-1}} j^i
        \right \rbrack \\
&=&	\left \lbrack
	{1 \over j} \sum_{i=0}^{j-1} {j \brace i} j^i
	\right \rbrack
	~-~ \left \lbrack
		\sum_{i=0}^{j-2} {{j-1} \brace i} j^i
	\right \rbrack \\
&=&	\left \lbrack
	{1 \over j} \sum_{i=0}^j {j \brace i} j^i
	\right \rbrack
	~-~ j^{j-1}
	~-~ \left \lbrack
		\sum_{i=0}^{j-1} {{j-1} \brace i} j^i
	\right \rbrack
	~+~ j^{j-1} \\
&=& {1 \over j} \, E[X_j^j] ~-~ E[X_j^{j-1}]
\end{eqnarray*}
Adding up all of the method (1) and method (2) shows that (for $j >
0$) the total number of edge attachments of $j$ free vertices to a
single edge $e \subseteq S$ is
\begin{equation}
\label{eq:attach-j>1-vertices-to-one-edge}
	E(1,j) ~=~ {1 \over j} \, E[X_j^j]
\end{equation}

We now need to show~(\ref{eq:E(i,j)}).
Let $q_0 = 1$ and for all $n \ge 1$, $q_n ~=~ E(1,n) ~=~ E[X_n^n]/n$.
We define
\begin{equation}
\label{eq:Q(z)-def}
	Q(z) ~=~ \sum_{n \ge 0} q_n {{z^n} \over {n!}}
\end{equation}
to be the exponential generating function for sequence $q_n$.

Now $q_j$ gives the number of edge attachments for $j$ free vertices
to a single fixed edge $e \subset S$.
In general, we must count the number of ways of attaching $j$ free
vertices to $i$ distinct edges $e_i \subseteq S$.
Partition the $j$ free vertices into $i$ subsets of sizes
$a_1,\,a_2,\,\ldots,\,a_i \ge 0$.  This can be done in
${i \choose {a_1,\,a_2,\, \ldots,\, a_i}}$ ways.  Each such partition
results in $q_{a_1} \, q_{a_2} \, \cdots \, q_{a_i}$ distinct
edge attachments.
The total number of edge attachments of $j$ free vertices to $i$
existing edges is therefore:
\begin{eqnarray*}
&&	\sum_{{a_k \ge 0} \atop {\sum_{k=1}^i a_k ~=~ j}}
	{j \choose {a_1,\, a_2,\, \ldots,\, a_i}}
	\prod_{k=1}^i q_{a_i} \\
&=& \lbrack {{z^i} \over {i!}} \rbrack \, Q(z)^j
\end{eqnarray*}
Equation~(\ref{eq:E(i,j)}) can now be seen to be nothing more than a
recursive statement that
\begin{eqnarray*}
&&	Q(z)^j ~=~ Q(z)^{j-1} \, Q(z)^1 \\
\Longrightarrow
&&	[{{z^i} \over {i!}}] \, Q(z)^j
	~=~ [{{z^i} \over {i!}}] \, Q(z)^{j-1} \, Q(z)^1
\end{eqnarray*}
\qed
\end{proof}

The total number of trees incident to subtour $S$ can now be computed
as follows.
We start with a tree over the $k$ vertices of $S$
composed of $i$ edges.
For each $0 \le i \le k-1$, there are
${1 \over k} \, {k-1 \brace i} \, k^i$
trees over $S$ having $i$ edges.
We partition the remaining $n-k$ free vertices into $m$ vertices that
connect via edge attachment to the $i$ edges of the main tree, and
$n-k-m$ vertices that connect via vertex attachment to the $k$
vertices of $S$.
For each $0 \le m \le n-k$ there are ${{n-k} \choose m}$ distinct
partitions of this type.
For each such partition there are
$V(n-k-m,k)$ vertex attachments.
and
$E(i,m)$ edge attachments
The total number of trees incident to subtour $S$ is therefore
\begin{eqnarray}
\nonumber
&&
{1 \over k} \sum_{i=0}^{k-1}
	{{k-1} \brace i} \, k^i
	\, \left \lbrack
		\sum_{m=0}^{n-k}
			{{n-k} \choose m}
			\, V(n-k-m,k)
			\, E(i,m)
	\right \rbrack \\
\nonumber
&=&
{1 \over k} \sum_{i=0}^{k-1}
	{{k-1} \brace i} \, k^i
	\, \left \lbrack
		\sum_{m=0}^{n-k}
			{{n-k} \choose m}
			\, g(n-m,k)
			\, E(i,m)
	\right \rbrack \\
\nonumber
&=&
{1 \over k} \sum_{i=0}^{k-1}
	{{k-1} \brace i} \, k^i
	\, \left \lbrack
		\sum_{m=0}^{n-k}
			{{n-k} \choose m}
			\, {k \over {n-m}} \, E[X_{n-m}^{n-m-k}]
			\, E(i,m)
	\right \rbrack \\
\label{eq:num-secs-ugly}
&=&
	\sum_{i=0}^{k-1}
	{{k-1} \brace i} \, k^i
	\, \left \lbrack
		\sum_{m=0}^{n-k}
			{{n-k} \choose m}
			\, {1 \over {n-m}} \, E[X_{n-m}^{n-m-k}]
			\, E(i,m)
	\right \rbrack
\end{eqnarray}
Dividing by the total number of trees ${1 \over n} E[X_n^{n-1}]$
yields the extreme point ratio
\begin{equation}
	{{n \sum_{i=0}^{k-1}
	{{k-1} \brace i} \, k^i
	\, \left \lbrack
		\sum_{m=0}^{n-k}
			{{n-k} \choose m}
			\, {1 \over {n-m}} \, E[X_{n-m}^{n-m-k}]
			\, E(i,m)
	\right \rbrack}
	\over
	{E[X_n^{n-1}]}}
\end{equation}
Using memoization to compute $E(i,j)$, binomial coefficients, Stirling
numbers and Poisson moments, this form can practically be computed for
at least $n \le 1000$.

We have been unable to obtain a closed-form for $E(i,j)$.
The following two approaches both fail.
\begin{lemma}
$[z^n] \, Q(z) ~=~ [z^n] \, \log(H(z))$ for all $n > 0$.
\end{lemma}
\begin{MyProof}

Apply the Lagrange inversion formula to $\log(H(z))$ with
$f(u) \,=\, \log(u)$ and $\theta(u) \,=\, e^{e^u \,-\, 1}$.
We find that
\begin{eqnarray*}
[z^n] \, \log(H(z))
&=& {1 \over n} \, \left \lbrack u^{n-1} \right \rbrack \,
	\lbrace u^{-1} \, \left ( e^{e^u \,-\, 1} \right )^n \rbrace \\
&=& {1 \over n} \, \left \lbrack u^n \right \rbrack \,
	\lbrace e^{n \, (e^u \,-\, 1)} \rbrace \\
&=& {1 \over n} \, E[X_n^n] \\
&=& q_n ~~~~~~~~\hbox{$\forall \, n > 0$.}
\end{eqnarray*}
\qed
\end{MyProof}
Although this is handy, it does not really lead anywhere for two reasons:
the equality does not hold for $n=0$; even if it did, then applying
the Lagrange inversion formula to $\log(H(z))^i$ for $i \ge 2$ results
in a $\log(u)$ factor for which no valid series expansion exists about
$u=0$.  This allows us to recover~(\ref{eq:E(1,j)}),
but does not provide an alternate form for~(\ref{eq:E(i,j)}).

The second failed approach is to find a simple closed form for
$E(2,j)$, hoping that a pattern emerges for $i>2$ that can be proven
by induction.
Assume $j > 0$, then
\begin{eqnarray*}
E(2,j)
&=&	\sum_{m=0}^j {j \choose m} \, E(1,j-m) \, E(1,m) \\
&=&	{j \choose 0} \, E(1,j) \, E(1,0)
	~+~ {j \choose j} \, E(1,0) \, E(1,j) \\
&&	~~~~ ~+~ \sum_{m=1}^{j-1} {j \choose m}
		\, E(1,j-m) \, E(1,m) \\
&=&	2 \, E(1,j)
	~+~ \sum_{m=1}^{j-1} {j \choose m}
		\, E(1,j-m) \, E(1,m) \\
&=&	{2 \over j} \, E[X_j^j]
	~+~ \sum_{m=1}^{j-1} {j \choose m}
		\, {1 \over {j-m}} \, E[X_{j-m}^{j-m}]
		\, {1 \over m} \, E[X_m^m]
\end{eqnarray*}
which is not a promising pattern to generalize.

We now derive CD formulae for two classes of facets in $\STHGP{n}$.

\subsection{Centroid of $\STHGP{n}$}

\begin{theorem}
\label{th:sthgp-centroid}
Let $n \ge 3$,
$H = (V,E)$ be the complete hypergraph with $n=|V|$,
and $C$ be the centroid of $\STHGP{n}$.  Let $e \in E$
and $k = |e|$.
Then
\begin{displaymath}
	C_e = {{g(n,|e|)} \over {g(n,1)}}
	    = {{k E[X_n^{n-k}]} \over {E[X_n^{n-1}]}}.
\end{displaymath}
\end{theorem}

\begin{MyProof}
The centroid is the sum of the incidence vectors of spanning trees
divided by the number of such incidence vectors.  Component $e$ of an
incidence vector is 1 (or 0) if $e$ is in (or not in) the
corresponding spanning tree.  Therefore the sum of component $e$ over
all incidence vectors is just the number of spanning trees containing
edge $e$, which is $g(n,|e|)$.  The total number of incidence vectors
is $g(n,1)$.
\qed
\end{MyProof}

\subsection{CD for General Polytopes with One Affine Hull Equation}
\label{sec:cd-of-general-polytope-one-equation}

Computing CDs is straightforward in a full-dimensional
polyhedron: one simply uses the formula for distance from a point $C$
to a hyperplane $a \cdot x = b$
\begin{displaymath}
	{{|b \,-\, a \cdot C|} \over {\sqrt{a \cdot a}}}
\end{displaymath}
and the closest point $x$ is
\begin{displaymath}
	x ~=~ C
	  ~+~ \left ( {{b \,-\, a \cdot C} \over {a \cdot a}} \right )
	      \, a
\end{displaymath}
Unfortunately, $\STHGP{n}$ is not full-dimensional, but satisfies a
single equation.
We first present a general method for computing CDs in
a polyhedron having one affine hull equation,
which we use for STHGP and elsewhere in this paper.

We are given a polyhedron $P \subset \R^n$,
$a,c \in \R^n$ and $b,d \in \R$
defining two hyperplanes:
\begin{eqnarray*}
	H_1 &=& \lbrace x \in \R^n : a \cdot x \,=\, b \rbrace, \\
	H_2 &=& \lbrace x \in \R^n : c \cdot x \,=\, d \rbrace.
\end{eqnarray*}
$H_1$ is the hyperplane of the inequality whose CD is to be
calculated.
Hyperplane $H_2$ is the affine hull of $P$.
We now obtain $\hat{a} \in \R^n$, $\hat{b} \in \R$ defining a new
hyperplane
\begin{eqnarray*}
	\hat{H}_1 &=& \lbrace x \in \R^n : \hat{a} \cdot x \,=\, \hat{b} \rbrace
\end{eqnarray*}
such that
\begin{equation}
\label{eq:same-intersection-with-affine-hull}
\hat{H}_1 \cap H_2 = H_1 \cap H_2
\end{equation}
and
\begin{equation}
\label{eq:normal-to-affine-hull}
\hat{a} \cdot c = 0
\end{equation}
so that $\hat{H}_1$ is normal to $H_2$,
as illustrated in
Figure~\ref{fig:non-affine-dist} (where $H$ is our general inequality
hyperplane $H_1$ and $A$ is our affine hull $H_2$) and
Figure~\ref{fig:affine-dist}
(where $H$ is our normalized hyperplane $\hat{H}_1$).

Let
\begin{eqnarray*}
q	&=& a \cdot a, \\
r	&=& a \cdot c, \\
s	&=& c \cdot c.
\end{eqnarray*}
Condition~(\ref{eq:same-intersection-with-affine-hull})
requires that $\hat{a}$ be a linear combination of
$a$ and $c$:
\begin{displaymath}
	\hat{a} = \alpha \, a - \beta \, c
\end{displaymath}
while condition~(\ref{eq:normal-to-affine-hull}) requires
\begin{eqnarray*}
\hat{a} \cdot c &=& (\alpha \, a - \beta \, c) \cdot c
	= \alpha \, a \cdot c - \beta \, c \cdot c
	= \alpha r - \beta s \\
&=&0
\end{eqnarray*}
so we can let $\alpha = s$ and $\beta = r$, so that
\begin{displaymath}
	\hat{a} ~=~ s \, a \,-\, r \, c.
\end{displaymath}
$\hat{a} = 0$ implies that either $H_1 = H_2$ or $H_1 \cap H_2 =
\emptyset$.
We assume $\hat{a} \ne 0$ in the sequel.
Note further that the sense of the given inequality
(e.g., $a \cdot x \le b$)
is preserved (e.g., $\hat{a} \cdot x \le \hat{b}$)
so long as $\alpha > 0$, which is guaranteed in this case
because $\alpha = c \cdot c$.
Moreover, $c \ne 0$ because $c$ defines a hyperplane $H_2$.

We want $\hat{H}_1$ to pass through the intersection of $H_1$ and $H_2$, so
we define a point $p \in H_2$:
\begin{displaymath}
	p = C + \tau \, \hat{a}
\end{displaymath}
where $C \in H_2$ is the centroid and we choose $\tau \in \R$ so that
$p \in H_1$:
\begin{displaymath}
a \cdot p = a \cdot (C + \tau \, \hat{a})
= a \cdot C + \tau \, a \cdot \hat{a} 
= b.
\end{displaymath}
This implies that
\begin{eqnarray*}
\tau &=& {{b - a \cdot C} \over
		{a \cdot \hat{a}}}
=	{{b - a \cdot C} \over
		{a \cdot (s \, a - r \, c)}}
=	{{b - a \cdot C} \over
		{s\,a \cdot a - r\,a \cdot c}} \\
&=&	{{b - a \cdot C} \over
		{q\,s - r^2}}.
\end{eqnarray*}
The distance from $C$ to $H_1 \cap H_2$ is now the same as the
distance from $C$ to $p$, which is just:
\begin{eqnarray}
\nonumber
d_w(H_1,P) &=& ||\tau \, \hat{a}||_2
 =	|\tau| \, \sqrt{\hat{a} \cdot \hat{a}} \\
\nonumber
&=&	|\tau| \, \sqrt{(s\,a-r\,c) \cdot (s\,a-r\,c)} \\
\nonumber
&=&	|\tau| \, \sqrt{q\,s^2 - 2\,s\,r^2 + s\,r^2}
 =	|\tau| \, \sqrt{s\,(q\,s - r^2)} \\
\label{eq:general-centroid-distance-dim-1}
&=&	{{\sqrt{s} \,|b - a \cdot C|} \over
	 {\sqrt{q\,s - r^2}}}.
\end{eqnarray}
If $p \in P$, then $d(H,P) = d_w(H,P)$.

Similar techniques could be developed to handle polyhedra satisfying
any fixed number of equations.

\subsection{Commonly Used Calculations for STHGP Subtours}
\label{sec:sthgp-common-subexprs}

We now apply the methods of
Section~\ref{sec:cd-of-general-polytope-one-equation}
to polytope $P = \STHGP{n}$, with particular emphasis on its subtour
inequalities.
The various functions and closed-forms derived in this section are
used in several subsequent sections.
We define the following functions:
\begin{eqnarray*}
	b(n) &=& n \, 2^{n-1}, \\
	c(n) &=& 2^n \,-\, 1, \\
	d(n) &=& 1 \,+\, n\,2^{n-1} \,-\, 2^n ~=~ b(n) \,-\, c(n), \\
	\alpha(n) &=& (n^2 \,-\, 3n \,+\, 4) \, 2^{n-2} \,-\, 1, \\
	\beta(n,k) &=& \gamma(n,k) \,+\, b(n \,-\, k) \, d(k), \\
	\gamma(n,k) &=& 2^{n-k} \, \alpha(k).
\end{eqnarray*}

\noindent
We use (but do not prove) the following identities:

\begin{displaymath}
	\sum_{i=0}^n {n \choose i} \, i ~=~ n \, 2^{n-1}
	~~~~~~~~~~~~~~~~~~ (n \ge 0),
\end{displaymath}

\iffalse
\begin{MyProof}
\begin{align*}
	\sum_{i=0}^n {n \choose i} \, i
&=	\sum_{i=1}^n {n \choose i} \, i
	= \sum_{i=1}^n {{n!} \over {i! \, (n-i)!}} \, i
	= n \, \sum_{i=1}^n {{(n-1)!} \over {(i-1)! \, (n-i)!}} \\
&=	n \, \sum_{i=0}^{n-1} {{(n-1)!} \over {i! \, (n-1-i)!}}
	= n \, \sum_{i=0}^{n-1} {{n-1} \choose i} \\
&=	n \, 2^{n-1}
\tag*{\qed}
\end{align*}
\end{MyProof}
\fi

\begin{displaymath}
	\sum_{i=0}^n {n \choose i} \, i^2
	~=~ n \, (n+1) \, 2^{n-2}
	~~~~~~~~ (n \ge 0).
\end{displaymath}

\iffalse
\begin{MyProof}
\begin{align*}
	\sum_{i=0}^n {n \choose i} \, i^2
&=	\sum_{i=1}^n {n \choose i} \, i^2
	= \sum_{i=1}^n {{n!} \over {i! \, (n-i)!}} \, i^2
	= n \, \sum_{i=1}^n {{(n-1)!} \over {(i-1)! \, (n-i)!}} \, i \\
&=	n \, \sum_{i=0}^{n-1} {{(n-1)!} \over {i! \, (n-1-i)!}}
			\, (i+1)
	= n \, \sum_{i=0}^{n-1} {{n-1} \choose i} \, (i + 1) \\
&=	n \, \sum_{i=0}^{n-1} {{n-1} \choose i} \, i +
		n \, \sum_{i=0}^{n-1} {{n-1} \choose i} \\
&=	n \, (n-1) \,2^{n-2} + n \, 2^{n-1}
\tag*{\qed}
\end{align*}
\end{MyProof}
\fi

Let $P = \STHGP{n}$, and
$H_1$ be the hyperplane of a subtour
inequality~(\ref{eq:sthgp-subtour-facet-def})
for some $S \subset V$ with $k = |S|$.
In equation $a \cdot x = b$, $a$ therefore satisfies
\begin{displaymath}
	a_e ~=~ \min(|e \cap S| \,-\, 1,\, 0)
	~~~~~~~ \forall e \in E,
\end{displaymath}
and $b = k-1$.
Hyperplane $H_2$ is the affine-hull
equation~(\ref{eq:sthgp-total-degree-equation})
for which our equation $c \cdot x = d$ has
$c_e = |e|-1$ for all $e \in E$,
and $d = n - 1$.

We now evaluate $\alpha = s = c \cdot c$ as a function of $n$,
summing $c \cdot c$ by edges of cardinality $i$:
{\allowdisplaybreaks%
\begin{eqnarray}
\nonumber
c \cdot c
&=&	\sum_{i=2}^n {n \choose i} (i \,-\, 1)^2 \\
\nonumber
&=&	-1 \,+\, \sum_{i=0}^n {n \choose i} (i \,-\, 1)^2 \\
\nonumber
&=&	-1 \,+\, \sum_{i=0}^n {n \choose i} i^2
	\,-\, 2 \sum_{i=0}^n {n \choose i} i
	\,+\, \sum_{i=0}^n {n \choose i} \\
\nonumber
&=&	-1 \,+ n \, (n \,+\, 1) \, 2^{n-2}
	\,-\, 2 (n \, 2^{n-1})
	\,+\, (2^n) \\
\nonumber
&=&	(n^2 \,-\, 3n \,+\, 4) \, 2^{n-2} \,-\, 1 \\
\label{eq:sthgp-c-dot-c}
\Longrightarrow~~~~\hfil{}
c \cdot c
&=&	\alpha(n).
\end{eqnarray}%
}%

\noindent
We next evalulate quantity $q = a \cdot a$ as a function of $n$ and $k$.
The only non-zero coefficients are for those edges $e$ having at least
2 vertices in the subtour $S$.
Consider that there are ${k \choose i}$ subsets of $S$ having
cardinality $i$.
For each such subset there are $2^{n-k}$ choices for that portion of
$e$ that resides outside of $S$.
The coefficient is $i-1$ for each such edge.
This yields:
{\allowdisplaybreaks%
\begin{eqnarray}
\nonumber
a \cdot a
&=&	\sum_{i=2}^k {k \choose i} (i \,-\, 1)^2 \, 2^{n-k} \\
\nonumber
&=&	2^{n-k} \sum_{i=2}^k {k \choose i} (i \,-\, 1)^2 \\
\nonumber
&=&	2^{n-k} \, \alpha(k) \\
\Longrightarrow~~~~\hfil{}
a \cdot a
&=&	\gamma(n,k).
\end{eqnarray}%
}%

\noindent
We next evaluate $\beta = r = a \cdot c$ as a function of $n$ and $k$.
We consider edge $e$ having $j$ vertices inside of $S$ and
$i$ vertices outside of $S$:
{\allowdisplaybreaks%
\begin{eqnarray}
\nonumber
a \cdot c
&=&	\sum_{j=2}^k {k \choose j} (j \,-\, 1)
	  \sum_{i=0}^{n-k} {{n - k} \choose i} (i \,+\, j \,-\, 1) \\
\nonumber
&=&	\sum_{j=2}^k {k \choose j} (j \,-\, 1)
	\left \lbrack
	  \sum_{i=0}^{n-k} {{n - k} \choose i} \, i
	  \,+\,
	  (j \,-\, 1) \sum_{i=0}^{n-k} {{n - k} \choose i}
	\right \rbrack \\
\nonumber
&=&	\sum_{j=2}^k {k \choose j} (j \,-\, 1)
	\left \lbrack
		(n \,-\, k) \, 2^{n-k-1} \,+\, (j \,-\, 1) \, 2^{n-k}
	\right \rbrack \\
\nonumber
&=&	\sum_{j=2}^k {k \choose j} (j \,-\, 1)
	\left \lbrack
		n \,-\, k \,+\, 2j \,-\, 2
	\right \rbrack
	\, 2^{n-k-1} \\
\nonumber
&=&	2^{n-k-1}
	\sum_{j=2}^k {k \choose j}
		(j \,-\, 1)
		\, (n \,-\, k \,+\, 2j \,-\, 2) \\
\nonumber
&=&	2^{n-k}
	\sum_{j=2}^k {k \choose j}
		(j \,-\, 1)^2
	\,+\, (n \,-\, k) \, 2^{n-k-1}
	\sum_{j=2}^k {k \choose j} (j \,-\, 1) \\
\nonumber
&=&	\gamma(n,k)
	\,+\, (n \,-\, k) \, 2^{n-k-1}
	\left \lbrack
	1 \,+\, 
	\sum_{j=0}^k {k \choose j} (j \,-\, 1)
	\right \rbrack \\
\nonumber
&=&	\gamma(n,k)
	\,+\, (n \,-\, k) \, 2^{n-k-1}
	\left \lbrack
	1
	\,+\, \sum_{j=0}^k {k \choose j} j
	\,-\, \sum_{j=0}^k {k \choose j}
	\right \rbrack \\
\nonumber
&=&	\gamma(n,k)
	\,+\, (n \,-\, k) \, 2^{n-k-1} (1 \,+\, k\,2^{k-1} \,-\, 2^k) \\
\nonumber
&=&	\gamma(n,k) \,+\, (n \,-\, k) \, 2^{n-k-1} \, d(k) \\
\Longrightarrow~~~~\hfil{}
a \cdot c
&=&	\beta(n,k).
\end{eqnarray}%
}%

\subsection{CD for STHGP Non-negativity Inequalities}
\label{sec:sthgp-nonneg-cd}

\begin{theorem}
Let $n \ge 2$, $G = (V,E)$ be the complete hypergraph with $|V|=n$,
$e \in E$ and $k = |e|$.
Then the (weak) CD for the non-negativity
constraint
\begin{equation}
\label{eq:sthgp-non-neg-ineq}
	x_e \,\ge\, 0
\end{equation}
of $\STHGP{n}$ is
\begin{displaymath}
d_w(n,k) ~=~
   {{k \, E[X_n^{n-k}]} \over {E[X_n^{n-1}]}}
   \,
   \sqrt{{\alpha(n)}
    \over {\alpha(n) \,-\, (k-1)^2}}
\end{displaymath}
where
$\alpha(n) ~=~ (n^2-3n+4) \, 2^{n-2} \,-\, 1$.
\end{theorem}

\begin{MyProof}
We use~(\ref{eq:general-centroid-distance-dim-1}) and results from
Section~\ref{sec:sthgp-common-subexprs}.
$H_1$ is equation $a \cdot x = b$ corresponding to the hyperplane of
the non-negativity
inequality~(\ref{eq:sthgp-non-neg-ineq}).
Hyperplane $H_2$ is equation $c \cdot x = d$ corresponding
to~(\ref{eq:sthgp-total-degree-equation}).
Theorem~\ref{th:sthgp-centroid} gives us the centroid $C$.
We need to compute the constants $q$, $r$, $s$ and also the dot
product $a \cdot C$.

We have $b = 0$, $q = a \cdot a = 1$ and $r = a \cdot c = k-1$.
From (\ref{eq:sthgp-c-dot-c}) we have
$s \,=\, c \cdot c \,=\, \alpha(n)$.
The dot product $a \cdot C$ is just $g(n,k) / g(n,1)$.

Putting this together yields the following CD for
non-negativity facets:
\begin{align*}
d_w(n,k)
&=	{{\sqrt{s} \,|b - a \cdot C|} \over
	 {\sqrt{q\,s - r^2}}}
=	{{g(n,k)} \over {g(n,1)}} \,
	{{\sqrt{\alpha(n)}} \over {\sqrt{\alpha(n) - (k-1)^2}}} \\
&=	{{k \, E[X_n^{n-k}]} \over {E[X_n^{n-1}]}} \,
	\sqrt{{\alpha(n)}
	 \over
	 {\alpha(n) - (k-1)^2}}.
\tag*{\qed}
\end{align*}
\end{MyProof}

\noindent
For fixed $n$, this distance decreases with $k$.
This can be verified computationally for small $n$.
When $n$ is large, the square root is very nearly 1, approaching it
from above.
The first ratio decreases rapidly with $k$.
In practice one works with a subhypergraph of $H$ having a relatively
small, fixed number of edges, so that no separation of these
inequalities is needed.
If they were to be separated, however, one would prefer to locate
violations $x_e < 0$ such that $|e|$ is large, since CD predicts these
to be stronger.

\subsection{CD for STHGP Subtour Inequalities}

\begin{theorem}
\label{th:sthgp-subtour-cd}
Let $n \ge 2$, $G = (V,E)$ be the complete hypergraph with $|V|=n$,
$S \subset V$ such that $|S| \ge 2$, and $k = |S|$.
Then the (weak) CD for the subtour inequality
\begin{equation}
\label{eq:sthgp-subtour-ineq}
	\sum_{e \in E} \min(|e \cap S| - 1, 0) \, x_e \le |S| - 1
\end{equation}
of $\STHGP{n}$ is
\begin{displaymath}
d_w(n,k) ~=~
  {{(n-k) \, t(n,k)} \over {E[X_n^{n-1}]}} \,
	\sqrt{{\alpha(n) \over {\alpha(n) \, \gamma(n,k) \,-\, \beta^2(n,k)}}},
\end{displaymath}
where
{\allowdisplaybreaks%
\begin{eqnarray*}
t(n,k) &=& E[X_n(X_n+1)^{n-2}] \,-\, E[X_n^k(X_n+1)^{n-k-1}], \\
\alpha(n) &=& (n^2 \,-\, 3n \,+\, 4) \, 2^{n-2} \,-\, 1, \\
\beta(n,k) &=& \gamma(n,k) \,+\, b(n \,-\, k) \, d(k), \\
\gamma(n,k) &=& 2^{n-k} \, \alpha(k), \\
b(n) &=& n \, 2^{n-1}, \\
d(n) &=& 1 \,+\, n\,2^{n-1} \,-\, 2^n.
\end{eqnarray*}%
}%
\end{theorem}

\noindent
The proof of Theorem~\ref{th:sthgp-subtour-cd}
requires several preliminary results.
We use (but do not prove) the following identities:
\begin{equation}
\label{eq:binom-sum-mi}
	\sum_{i=0}^k {k \choose i} \, i \, m^{k-i}
	~=~ k \, (m + 1)^{k-1},
\end{equation}
\begin{equation}
\label{eq:binom-sum-m2i}
	\sum_{i=0}^k {k \choose i} \, i^2 \, m^{k-i}
	~=~ k(k-1)(m+1)^{k-2} \, + \, k(m+1)^{k-1}.
\end{equation}
We will also need several lemmas:
\begin{MyLemma}
\label{lem:sthgp-subtour-cd-lemma-1}
\begin{displaymath}
	\sum_{j=0}^{n-k} {{n-k} \choose j} (i+j) m^{n-i-j}
	~=~
	(m+1)^{n-k-1} m^{k-i} \left [ i (m+1) + n - k \right ].
\end{displaymath}
\end{MyLemma}
\begin{MyProof}
{\allowdisplaybreaks%
\begin{align*}
\lefteqn{\sum_{j=0}^{n-k} {{n-k} \choose j} (i+j) m^{n-i-j}} \\
&= m^{k-i} \sum_{j=0}^{n-k} {{n-k} \choose j} (i+j) m^{n-k-j} \\
&= m^{k-i}
	\left [
	i \sum_{j=0}^{n-k} {{n-k} \choose j} m^{n-k-j}
	+ \sum_{j=0}^{n-k} {{n-k} \choose j} j \, m^{n-k-j}
	\right ] \\
&= m^{k-i} \left [ i (m+1)^{n-k} + (n-k) (m+1)^{n-k-1} \right ] \\
&= (m+1)^{n-k-1} m^{k-i} \left [ i (m+1) + n - k \right ].
\tag*{\qed}
\end{align*}%
}%
\end{MyProof}
\begin{MyLemma}
\label{lem:moment_X+1}
Let $X_{\lambda}$ be a Poisson random variable with parameter
$\lambda$, and let $n \ge 0$.  Then
\begin{displaymath}
	E[(X_{\lambda} + 1)^n]
	~=~ {1 \over {\lambda}} E[X_{\lambda}^{n+1}].
\end{displaymath}
\end{MyLemma}

\begin{MyProof}
{\allowdisplaybreaks%
\begin{align*}
E[(X_{\lambda}+1)^n]
&= \sum_{i \ge 0} {{{(i+1)}^n\,e^{-\lambda}\,\lambda^i} \over {i!}}
	= \sum_{i \ge 1} {{i^n\,e^{-\lambda}\,\lambda^{i-1}} \over
			  {(i-1)!}} \\
&=	{1 \over {\lambda}} \,
	\sum_{i \ge 1} {{i^{n+1}\,e^{-\lambda}\,\lambda^i} \over
			{i!}}
	=
	{1 \over {\lambda}} \,
	\sum_{i \ge 0} {{i^{n+1}\,e^{-\lambda}\,\lambda^i} \over
			{i!}} \\
&=	{1 \over {\lambda}} \, E[X_{\lambda}^{n+1}].
\tag*{\qed}
\end{align*}%
}%
\end{MyProof}
\begin{MyLemma}
\label{lem:sthgp-subtour-cd-lemma-2}
\begin{eqnarray*}
&&	\sum_{i=2}^k {k \choose i} (i - 1)
		\left [ i (m+1) + n - k \right ] m^{k-i} \\
&=&
	(n-k)m^k + k(n-1)(m+1)^{k-1} - (n-k)(m+1)^k.
\end{eqnarray*}
\end{MyLemma}
\begin{MyProof}
{\allowdisplaybreaks%
\begin{eqnarray*}
\lefteqn{\sum_{i=2}^k {k \choose i} (i - 1)
	\left [ i (m+1) + n - k \right ] m^{k-i}} \\
&=&	(n-k)m^k
	+ \sum_{i=0}^k {k \choose i} (i-1)
	\left [ i (m+1) + n - k \right ] m^{k-i} \\
&=&	(n-k)m^k
	+ \sum_{i=0}^k {k \choose i}
	\left [ (m+1)i^2 + (n-k-m-1) i - (n - k) \right ] m^{k-i} \\
&=&	(n-k)m^k
	+ (m+1) \sum_{i=0}^k {k \choose i} i^2 m^{k-i} \\
&&	\mbox{}
	+ (n-k-m-1) \sum_{i=0}^k {k \choose i} i m^{k-i}
	- (n-k) \sum_{i=0}^k {k \choose i} m^{k-i}
\end{eqnarray*}%
}%
which by~(\ref{eq:binom-sum-mi}) and~(\ref{eq:binom-sum-m2i}) yields
{\allowdisplaybreaks%
\begin{align*}
&=	(n-k)m^k
	+ (m+1) \left [ k(k-1)(m+1)^{k-2} + k(m+1)^{k-1} \right ] \\
&~~~	\mbox{}
	+ (n-k-m-1) k (m+1)^{k-1}
	- (n-k) (m+1)^k \\
&= (n-k) \, m^k \,+\, k \, \left (
		(k \,-\, 1) \,+\, (m \,+\, 1) \,+\,
		(n \,-\, k \,-\, m \,-\, 1)
	\right ) \, (m \,+\, 1)^{k-1} \\
&~~~	\,-\, (n \,-\, k) \, (m \,+\, 1)^k \\
&=	(n-k)m^k + k(n-1)(m+1)^{k-1} - (n-k)(m+1)^k.
\tag*{\qed}
\end{align*}%
}%
\end{MyProof}
\begin{MyLemma}
\label{lem:sthgp-subtour-cd-lemma-3}
Let $n \ge 2$.  Then
\begin{displaymath}
	E[X_n^n] \,-\, E[X_n^{n-1}]
	~=~
	n \, E[X_n \, (X_n + 1)^{n-2}].
\end{displaymath}
\end{MyLemma}
\begin{MyProof}
{\allowdisplaybreaks%
\begin{align*}
E[X_n^n] \,-\, E[X_n^{n-1}]
&=	E[(X_n - 1)X_n^{n-1}] \\
&=	\sum_{i \ge 0}
		{{(i - 1) i^{n-1} e^{-n} n^i} \over {i!}}
~=~	\sum_{i \ge 1}
		{{(i - 1) i^{n-1} e^{-n} n^i} \over {i!}} \\
&=	\sum_{i \ge 0}
		{{i(i + 1)^{n-1} e^{-n} n^{i+1}} \over {(i+1)!}}
~=~	n \, \sum_{i \ge 0}
		{{i(i + 1)^{n-2} e^{-n} n^i} \over {i!}} \\
&=	n \, E[X_n \, (X_n + 1)^{n-2}]
\tag*{\qed}
\end{align*}%
}%
\end{MyProof}

\begin{MyProof}[of Theorem~\ref{th:sthgp-subtour-cd}]
We use~(\ref{eq:general-centroid-distance-dim-1}) and results from
Section~\ref{sec:sthgp-common-subexprs}.
$H_1$ is given by equation $a \cdot x = b$ corresponding to hyperplane
of~(\ref{eq:sthgp-subtour-ineq}).
$H_2$ is given by equation $c \cdot x = d$ correpsonding
to~(\ref{eq:sthgp-total-degree-equation}).
Theorem~\ref{th:sthgp-centroid} gives the centroid $C$.
Using the results from Section~\ref{sec:sthgp-common-subexprs}
we have
\begin{eqnarray*}
	q &=& a \cdot a = \gamma(n,k), \\
	r &=& a \cdot c = \beta(n,k), \\
	s &=& c \cdot c = \alpha(n).
\end{eqnarray*}
The final quantity needed is $a \cdot C$.
For given edge $e$, we have $C_e = g(n,|e|) / g(n,1)$.
Consider an edge $e$ having $i$ vertices in $S$ and $j$ vertices not
in $S$.
This gives:
{\allowdisplaybreaks%
\begin{eqnarray*}
\lefteqn{a \cdot C} \\
&=&	\sum_{i=2}^k {k \choose i} (i - 1)
		\sum_{j=0}^{n-k} {{n-k} \choose j}
			{{g(n,i+j)} \over {g(n,1)}} \\
&=&	{1 \over {g(n,1)}}
	\sum_{i=2}^k {k \choose i} (i - 1)
		\sum_{j=0}^{n-k} {{n-k} \choose j}
			{{i+j} \over n} E[X_n^{n-i-j}] \\
&=&	{1 \over {n \, g(n,1)}}
	\sum_{i=2}^k {k \choose i} (i - 1)
		\sum_{j=0}^{n-k} {{n-k} \choose j}
			(i+j) E[X_n^{n-i-j}] \\
&=&	{1 \over {n \, g(n,1)}}
	\sum_{i=2}^k {k \choose i} (i - 1)
		\sum_{j=0}^{n-k} {{n-k} \choose j}
			(i + j)
			\sum_{m \ge 0} {{m^{n-i-j} e^{-n} n^m} \over
				{m!}} \\
&=&	{{e^{-n}} \over {n \, g(n,1)}}
	\sum_{m \ge 0} {{n^m} \over {m!}}
	  \sum_{i=2}^k {k \choose i} (i - 1)
	    \sum_{j=0}^{n-k} {{n-k} \choose j} (i+j) m^{n-i-j}.
\end{eqnarray*}%
}%
By Lemma~\ref{lem:sthgp-subtour-cd-lemma-1} this yields
{\allowdisplaybreaks%
\begin{eqnarray*}
\lefteqn{a \cdot C} \\
&=&	{{e^{-n}} \over {n \, g(n,1)}}
	\sum_{m \ge 0} {{n^m} \over {m!}}
	  \sum_{i=2}^k {k \choose i} (i - 1)
	  (m+1)^{n-k-1} m^{k-i} \left [ i (m+1) + n - k \right ] \\
&=&	{{e^{-n}} \over {n \, g(n,1)}}
	\sum_{m \ge 0} {{(m+1)^{n-k-1} n^m} \over {m!}}
	  \sum_{i=2}^k {k \choose i} (i - 1)
	    \left [ i (m+1) \,+\, n \,-\, k \right ] \, m^{k-i}
\end{eqnarray*}%
}%
which by Lemma~\ref{lem:sthgp-subtour-cd-lemma-2} yields
{\allowdisplaybreaks%
\begin{eqnarray*}
\lefteqn{a \cdot C} \\
&=&	{{e^{-n}} \over {n \, g(n,1)}}
	\sum_{m \ge 0} {{(m+1)^{n-k-1} n^m} \over {m!}}
	  \begin{array}[t]{l}
	    \left [ (n-k)m^k + k(n-1)(m+1)^{k-1} \right. \\
	    \left. \hbox{~~~} - (n-k)(m+1)^k \right ]
	  \end{array} \\
&=&
	{1 \over {n \, g(n,1)}}
	  \begin{array}[t]{l}
	    \left [
		(n-k) \sum_{m \ge 0} {{m^k(m+1)^{n-k-1} e^{-n} n^m}
				      \over {m!}}
	    \right. \\
		\hbox{~~~} + k(n-1) \sum_{m \ge 0}
			{{(m+1)^{n-2} e^{-n} n^m} \over {m!}} \\
	    \left.
		\hbox{~~~} - (n-k) \sum_{m \ge 0}
			{{(m+1)^{n-1} e^{-n} n^m} \over {m!}}
	    \right ]
	  \end{array} \\
&=&	{1 \over {n \, g(n,1)}}
	  \begin{array}[t]{l}
		\left [ (n-k) E[X_n^k(X_n+1)^{n-k-1}] \right. \\
			\hbox{~~~} + k(n-1) E[(X_n + 1)^{n-2}] \\
		\left.	\hbox{~~~} - (n-k) E[(X_n + 1)^{n-1}]
			\,\, \right ].
	  \end{array}
\end{eqnarray*}%
}%
Using two applications of Lemma~\ref{lem:moment_X+1} this yields
{\allowdisplaybreaks%
\begin{eqnarray*}
&=&	{1 \over {n \, g(n,1)}}
	    \left [
		(n-k) E[X_n^k(X_n+1)^{n-k-1}]
		\,+\, {{k(n-1)} \over n} \, E[X_n^{n-1}]
	    \right. \\
&&	~~~~~~~~~~~~~~~~
	    \left.
		\,-\, {{(n-k)} \over n} E[X_n^n]
	    \right ] \\
&=&	{1 \over {n^2 \, g(n,1)}}
	    \left [
		n(n-k) \, E[X_n^k(X_n+1)^{n-k-1}]
		\,+\, k(n-1) \, E[X_n^{n-1}]
	    \right. \\
&&	~~~~~~~~~~~~~~~~~
	\left.
		\,-\, (n-k) \, E[X_n^n] \,\,
	    \right ].
\end{eqnarray*}%
}%

\noindent
From this we conclude
{\allowdisplaybreaks%
\begin{eqnarray*}
\lefteqn{b - a \cdot C} \\
&=&	k-1 - a \cdot C \\
&=&	{{n^2 \, (k-1) \, g(n,1)} \over {n^2 \, g(n,1)}}
	- a \cdot C \\
&=&	{1 \over {n^2 \, g(n,1)}}
	\left \lbrack
		n \, (k-1) \, E[X_n^{n-1}]
		- n(n-k) E[X_n^k(X_n+1)^{n-k-1}]
	\right. \\
&&	~~~~~~~~~~~~~~~
	\left.
	  - k(n-1) E[X_n^{n-1}]
	  + (n-k) E[X_n^n]
	\right \rbrack \\
&=&	{{(n-k) \, \left (
		E[X_n^n] - E[X_n^{n-1}]
		- n E[X_n^k(X_n+1)^{n-k-1}]
	  \right )}
	 \over {n^2 \, g(n,1)}} \\
&=&	{{(n-k)} \over {n \, E[X_n^{n-1}]}} \,
	  \left [
		E[X_n^n] - E[X_n^{n-1}]
		- n E[X_n^k(X_n+1)^{n-k-1}
	  \right ]
\end{eqnarray*}%
}%
which by Lemma~\ref{lem:sthgp-subtour-cd-lemma-3} yields
{\allowdisplaybreaks%
\begin{eqnarray}
\nonumber
\lefteqn{b - a \cdot C} \\
\nonumber
&=&	{{(n-k)} \over {n \, E[X_n^{n-1}]}} \,
	  \left [
		n E[X_n (X_n+1)^{n-2}]
		- n E[X_n^k(X_n+1)^{n-k-1}]
	  \right ] \\
\nonumber
&=&	{{(n-k) \, E[X_n (X_n+1)^{n-2} - X_n^k(X_n+1)^{n-k-1}]}
	 \over {E[X_n^{n-1}]}} \\
\label{eq:N1_dot_C_d}
&=&	{{(n-k) \, t(n,k)} \over {E[X_n^{n-1}]}}.
\end{eqnarray}%
}%

\noindent
The weak CD $d_w(n,k)$
for subtour $S \subset V$ with $|S| = k$ is therefore
\begin{equation}
\label{eq:sthgp-subtour-cd-proof-final-form}
	d_w(n,k) ~=~
	{{(n - k) \, t(n,k)} \over {E[X_n^{n-1}]}} \,
	\sqrt{{{\alpha(n)} \over
	       {\alpha(n) \, \gamma(n,k) \,-\, \beta^2(n,k)}}}.
\end{equation}
\qed
\end{MyProof}

\noindent
We note that the $(n-k)$ factor is a decreasing
function of $k$, and that $t(n,k)$ is zero when $k=1$
and strictly increases.  Therefore, (\ref{eq:N1_dot_C_d}) is zero when
$k=1$ and when $k=n$, and is positive and concave down in between.

Equation~(\ref{eq:sthgp-subtour-cd-proof-final-form})
can be easily evaluated
using~(\ref{eq:poisson-moment-as-stirling-sum})
and the following Lemma:

\begin{MyLemma}
\label{lem:moment_Xm(X+1)n}
Let $X_{\lambda}$ be a Poisson random variable with parameter
$\lambda$, and let $n,m \ge 0$.  Then
\begin{displaymath}
	E[X_{\lambda}^m (X_{\lambda} + 1)^n] ~=~
	\sum_{j=0}^n {n \choose j} E[X_{\lambda}^{m+j}].
\end{displaymath}
\end{MyLemma}

\begin{MyProof}
{\allowdisplaybreaks%
\begin{align*}
E[X_{\lambda}^m (X_{\lambda} + 1)^n]
&= \sum_{i \ge 0} {{i^m (i+1)^n e^{-\lambda} \lambda^i} \over {i!}}
	= \sum_{i \ge 0} {{i^m e^{-\lambda} \lambda^i} \over {i!}}
		\sum_{j=0}^n {n \choose j} i^j \\
&= \sum_{j=0}^n {n \choose j}
	\sum_{i \ge 0} {{i^{m+j} e^{-\lambda} \lambda^i} \over
			{i!}}
= \sum_{j=0}^n {n \choose j} E[X_{\lambda}^{m+j}].
\tag*{\qed}
\end{align*}%
}%
\end{MyProof}

\subsection{Subtour Angles}
\label{sec:subtour-angles}

Let $H = (V,E)$ be a complete hypergraph.
The corresponding subtour relaxation is:
\begin{eqnarray}
\label{eq:angle-affine-hull}
 \sum_{e \in E} (|e| \,-\, 1) \, x_e
	&=& |V| \,-\, 1, \\
\label{eq:angle-subtour-S}
 \sum_{e \in E} \max(0, |e \cap S| \,-\, 1) \, x_e
	&\le& |S| \,-\, 1
	~~~~~~~ \forall S \subset V,\, |S| \ge 2, \\
\nonumber
 x_e	 &\ge& 0
	~~~~~~~ \forall e \in E.
\end{eqnarray}
Equation~(\ref{eq:angle-affine-hull}) is the affine hull,
and~(\ref{eq:angle-subtour-S}) are the subtour inequalities.

Let $n = |V|$.
We require that $n \ge 3$ throughout this section so that subtour
inequalities exist.
Let $S_1,\, S_2 \subset V$ such that $|S_1| \ge 2$ and $|S_2| \ge 2$.
In this section we study the angle formed between the subtour
inequalities corresponding to $S_1$ and $S_2$.
There are several reasons to study these angles:
\begin{itemize}
  \item When the interior (feasible region) angle is small (i.e., near
    zero), we should expect the combination of inequalities $S_1$ and
    $S_2$ to be quite strong in practice.
  \item When the interior angle is large (i.e., near $\pi$ radians),
    we should expect this combination to provide relatively little
    improvement in strength over that of $S_1$ or $S_2$ individually.
  \item When this angle approaches either zero or $\pi$ radians, the
    corresponding cutting planes become nearly parallel, which would
    impact the condition number of the LP basis matrix (when both
    inequalities are simultaneously binding).
\end{itemize}
We note that the angle between any two such hyperplanes should always
be measured within the affine hull of the polyhedron, as these angles
can otherwise be made arbitrarily close to zero or $\pi$ radians by
rotating the hyperplanes about their intersection with the affine
hull.

Let the subtour inequality associated with
$S_1$ be $a_1 \cdot x \le b_1$,
and the subtour inequality associated with
$S_2$ be $a_2 \cdot x \le b_2$.
Let $\hat a_1,\, \hat a_2$
be the corresponding (scaled) projections of
$a_1,\, a_2$
onto the affine hull~(\ref{eq:angle-affine-hull}),
as described in Section~\ref{sec:cd-of-general-polytope-one-equation}
above.
Let $\phi$ be the angle between vectors $\hat a_1$ and $\hat a_2$.
Then
\begin{eqnarray}
\nonumber
	&& \hat a_1 \cdot \hat a_2 ~=~
		||\hat a_1|| \, ||\hat a_2|| \, \cos \phi \\
\label{eq:cos-phi}
\Longrightarrow
	&& \cos \phi ~=~ {{\hat a_1 \cdot \hat a_2} \over
			    {||\hat a_1|| \, ||\hat a_2||}}.
\end{eqnarray}
The interior angle that we seek is then
\begin{equation}
\label{eq:theta=pi-phi}
	\theta ~=~ \pi \,-\, \phi.
\end{equation}

\noindent
For given $S_1$ and $S_2$, let
\begin{eqnarray*}
	n &=& |V|, \\
	p &=& |S_1 \setminus S_2|, \\
	q &=& |S_2 \setminus S_1|, \hbox{ and} \\
	r &=& |S_1 \cap S_2|.
\end{eqnarray*}
We seek a function
\begin{equation}
\label{eq:f(n,p,q,r)-dot-prod-equivalence}
f(n,p,q,r)
	~=~ \cos \phi
	~=~ {{\hat a_1 \cdot \hat a_2} \over
			    {||\hat a_1|| \, ||\hat a_2||}}
\end{equation}
so that
\begin{displaymath}
	\theta(n,p,q,r) ~=~ \pi \,-\, \mathrm{acos} (f (n,p,q,r)).
\end{displaymath}

\noindent
We define
\begin{eqnarray*}
\mu(n,k) &=&
	\alpha(n) \, \gamma(n,k) \,-\, \beta^2(n,k), \\
w(n,p,q,r) &=&
	2^{n-p-q-r} \,
	\left \lbrack
		2^{p+q}\alpha(r)
		\,+\, d(p)\,d(q)
		\,+\, b(p)\,b(q)\,c(r)
		\,+\,b(p+q)\,d(r)
	\right \rbrack.
\end{eqnarray*}

\begin{theorem}
\label{th:f(n,p,q,r)}
\begin{equation}
\label{eq:f(n,p,q,r)}
f(n,p,q,r) ~=~
	{{	\alpha(n) \, w(n,p,q,r)
		\,-\, \beta(n,p+r) \, \beta(n,q+r)}
	\over
	{\sqrt{\mu(n,p+r) \, \mu(n,q+r)}}}.
\end{equation}
\end{theorem}

\noindent
To prove this, we will need the following lemma:
\begin{MyLemma}
\label{lem:a1-dot-a2=w(n,p,q,r)}
\begin{displaymath}
a_1 \cdot a_2 ~=~ w(n,p,q,r).
\end{displaymath}
\end{MyLemma}

\noindent
We first prove Theorem~\ref{th:f(n,p,q,r)} and then prove
Lemma~\ref{lem:a1-dot-a2=w(n,p,q,r)}.

\begin{MyProof}[of Theorem~\ref{th:f(n,p,q,r)}]
From Section~\ref{sec:sthgp-common-subexprs} above, given a $k$-vertex
subtour hyperplane $a \cdot x = b$, we project $a$ onto the affine
hull as:
\begin{displaymath}
	\hat{a} ~=~ \alpha(n) \, a \,-\, \beta(n,k) \, c
\end{displaymath}
whose magnitude is:
{\allowdisplaybreaks%
\begin{eqnarray}
\nonumber
||\hat a||
&=&	||\alpha(n) \, a \,-\, \beta(n,k) \, c|| \\
\nonumber
&=&	\sqrt{(\alpha(n) \, a \,-\, \beta(n,k) \, c) \cdot
	      (\alpha(n) \, a \,-\, \beta(n,k) \, c)} \\
\nonumber
&=&	\sqrt{\alpha^2(n) \, (a \cdot a)
	      \,-\, 2 \, \alpha(n) \, \beta(n,k) \, (a \cdot c)
	      \,+\, \beta^2(n,k) \, (c \cdot c)} \\
\nonumber
&=&	\sqrt{\alpha^2(n) \, \gamma(n,k)
	      \,-\, 2 \, \alpha(n) \, \beta^2(n,k)
	      \,+\, \alpha(n) \, \beta^2(n,k)} \\
\nonumber
&=&	\sqrt{\alpha(n) \,
		\left \lbrack
			\alpha(n) \, \gamma(n,k)
			\,-\, \beta^2(n,k)
		\right \rbrack} \\
\label{eq:angle-ahat-magnitude}
||\hat a||
&=&	\sqrt{\alpha(n) \, \mu(n,k)}.
\end{eqnarray}%
}%

\noindent
We now return to the original subtours $S_1$ and $S_2$, with
corresponding inequalities
$a_1 \cdot x \le b_1$ and
$a_2 \cdot x \le b_2$.
Recalling the parameterization in terms of $n,\, p,\, q,\,$ and $r$,
the projected versions of $a_1$ and $a_2$ are:
\begin{eqnarray}
\label{eq:angle-project-a1}
	\hat a_1 &=& \alpha(n) \, a_1 \,-\, \beta(n,p \,+\, r) \, c, \\
\label{eq:angle-project-a2}
	\hat a_2 &=& \alpha(n) \, a_2 \,-\, \beta(n,q \,+\, r) \, c.
\end{eqnarray}
Their magnitudes (which follow from~(\ref{eq:angle-ahat-magnitude})) are:
\begin{eqnarray*}
||\hat a_1||	&=&	\sqrt{\alpha(n) \, \mu(n,p+r)}, \\
||\hat a_2||	&=&	\sqrt{\alpha(n) \, \mu(n,q+r)}.
\end{eqnarray*}
We will need the following dot products (refer to
Section~\ref{sec:sthgp-common-subexprs}):
{\allowdisplaybreaks%
\begin{eqnarray*}
a_1 \cdot a_1	&=&	\gamma(n,p+r), \\
a_1 \cdot c	&=&	\beta(n,p+r), \\
a_2 \cdot a_2	&=&	\gamma(n,q+r), \\
a_2 \cdot c	&=&	\beta(n,q+r).
\end{eqnarray*}%
}%
The main dot product we need is:
{\allowdisplaybreaks%
\begin{eqnarray*}
\hat a_1 \cdot \hat a_2
&=&	(\alpha(n) \, a_1 \,-\, \beta(n,p+r) \, c) \, \cdot \,
	(\alpha(n) \, a_2 \,-\, \beta(n,q+r) \, c) \\
&=&	\alpha^2(n) \, (a_1 \cdot a_2)
	\,-\, \alpha(n) \, \beta(n,q+r) \, (a_1 \cdot c) \\
&&	\,-\, \alpha(n) \, \beta(n,p+r) \, (a_2 \cdot c)
	\,+\, \beta(n,p+r) \, \beta(n,q+r) \, (c \cdot c) \\
&=&	\alpha^2(n) \, (a_1 \cdot a_2)
	\,-\, \alpha(n) \, \beta(n,q+r) \, \beta(n,p+r) \\
&&	\,-\, \alpha(n) \, \beta(n,p+r) \, \beta(n,q+r)
	\,+\, \alpha(n) \, \beta(n,p+r) \, \beta(n,q+r) \\
\Longrightarrow~~~~\hfil{}
\hat a_1 \cdot \hat a_2
&=&	\alpha(n) \,
	\left \lbrack
		\alpha(n) \, (a_1 \cdot a_2)
		\,-\, \beta(n,q+r) \, \beta(n,p+r)
	\right \rbrack.
\end{eqnarray*}%
}%
The function we desire is then
{\allowdisplaybreaks%
\begin{eqnarray*}
f(n,p,q,r)
&=&	{{\hat a_1 \cdot \hat a_2} \over
	 {||\hat a_1|| \, ||\hat a_2||}} \\
&=&	{{	\alpha(n) \,
	\left \lbrack
		\alpha(n) \, (a_1 \cdot a_2)
		\,-\, \beta(n,q+r) \, \beta(n,p+r)
	\right \rbrack}
	\over
	{\sqrt{\alpha(n) \, \mu(n,p+r)}
	 \,
	\sqrt{\alpha(n) \, \mu(n,q+r)}}} \\
&=&	{{\alpha(n) \,
	\left \lbrack
		\alpha(n) \, (a_1 \cdot a_2)
		\,-\, \beta(n,q+r) \, \beta(n,p+r)
	\right \rbrack}
	\over
	{\alpha(n)
	 \,
	 \sqrt{\mu(n,p+r) \, \mu(n,q+r)}}} \\
&=&	{{	\alpha(n) \, (a_1 \cdot a_2)
		\,-\, \beta(n,p+r) \, \beta(n,q+r)}
	\over
	{\sqrt{\mu(n,p+r) \, \mu(n,q+r)}}}.
\end{eqnarray*}%
}%
which by Lemma~\ref{lem:a1-dot-a2=w(n,p,q,r)} yields:
\begin{displaymath}
f(n,p,q,r) ~=~
	{{	\alpha(n) \, w(n,p,q,r)
		\,-\, \beta(n,p+r) \, \beta(n,q+r)}
	\over
	{\sqrt{\mu(n,p+r) \, \mu(n,q+r)}}}.
\tag*{\qed}
\end{displaymath}
\end{MyProof}

\noindent
We now prove the more involved
Lemma~\ref{lem:a1-dot-a2=w(n,p,q,r)}.

\begin{MyProof}[of Lemma~\ref{lem:a1-dot-a2=w(n,p,q,r)}]
Partition $V$ into four disjoint sets as follows:
{\allowdisplaybreaks%
\begin{eqnarray*}
V_0 &=& \lbrace v \in V : (v \notin S_1) \wedge (v \notin S_2) \rbrace, \\
V_1 &=& \lbrace v \in V : (v \in S_1) \wedge (v \notin S_2) \rbrace, \\
V_2 &=& \lbrace v \in V : (v \notin S_1) \wedge (v \in S_2) \rbrace, \\
V_3 &=& \lbrace v \in V : (v \in S_1) \wedge (v \in S_2) \rbrace
\end{eqnarray*}%
}%
so that $S_1 \,=\, V_1 \cup V_3$ and $S_2 \,=\, V_2 \cup V_3$.
(The V subscripts can be thought of as a 2-bit binary number, where
the most-significant bit indicates membership within $S_2$ and the
least-significant bit indicates membership within $S_1$.)

In similar fashion,
the edges $E$ can be partitioned into 15 classes
$E_1,\, \ldots,\, E_{15}$
where the subscript is now a 4-bit binary number where the $2^i$ bit
is one if $|e \cap V_i| \ge 1$ and zero otherwise.
(Note that edges have at least two vertices, so there is no edge class
$E_0$.)
The dot product $a_1 \cdot a_2$ is a sum of one term for each edge
$e \in E$.
We will evaluate $a_1 \cdot a_2$ by developing an expression
$K_i = \sum_{e \in E_i} a_1(e) \, a_2(e)$ that gives the partial sum
for each edge class $i$, $(1 \le i \le 15)$.
The full dot-product is then $a_1 \cdot a_2 = \sum_{i=1}^{15} K_i$.

The first observation is that vertices $V_0$
do not affect the coefficient values for any edge $e$ in
either $a_1$ or $a_2$.
Let $m\,=\,|V_0|\,=\,n-p-q-r$.
For $e \subseteq V \setminus V_0$, we define
$\xi(e) ~=~ \lbrace e \cup \bar{e} : \bar{e} \subseteq V_0 \rbrace$
Thus, $\xi(e)$ partitions $E$ into disjoint sets, each having $2^m$
members.
This gives rise to the following relationships between the even and
odd partial sums (for $1 \le i \le 7$):
\begin{eqnarray*}
&&	K_{2i+1} ~=~ (2^m \,-\, 1) \, K_{2i} \\
\Longrightarrow
&&	K_{2i} \,+\, K_{2i+1} ~=~ 2^m \, K_{2i}.
\end{eqnarray*}
It therefore suffices to sum the $K_{2i}$ and multiply by $2^m$.
Since there is no edge class $E_0$, however, we must consider edge
class $E_1$ separately.
The sum $K_1 = 0$ because for all $e \in E_1$, we have
$e \subseteq V_0$
which means that the $e$ components of $a_1$ and $a_2$ are both zero.

We now consider each $(E_{2i},K_{2i})$ pair, for $1 \le i \le 7$.
$K_2$ is zero because for all $e \in E_2$, the $e$ component of
$a_2$ is zero.
$K_4$ is likewise zero: for all $e \in E_4$, the $e$ component of
$a_1$ is zero.
Therefore:
\begin{displaymath}
a_1 \cdot a_2 ~=~
2^{n-p-q-r} \,
	(K_6 \,+\, K_8 \,+\, K_{10} \,+\, K_{12} \,+\, K_{14})
\end{displaymath}

\noindent
Edges $e \in E_6$ have vertices in $V_1$ and $V_2$, giving rise to
dot product terms of the form
$(i - 1) \, (j - 1)$, where
$i = |e \cap V_1|$ and $j = |e \cap V_2|$:
{\allowdisplaybreaks%
\begin{eqnarray*}
K_6
&=&	\sum_{i=1}^p {p \choose i}
	 \sum_{j=1}^q {q \choose j}
		(i \,-\, 1) \, (j \,-\, 1) \\
&=&	\left \lbrack
	\sum_{i=1}^p {p \choose i}
	 (i \,-\, 1)
	\right \rbrack
	\,
	\left \lbrack
	 \sum_{j=1}^q {q \choose j}
		(j \,-\, 1)
	\right \rbrack \\
&=&
	(1 \,+\, p \, 2^{p-1} \,-\, 2^p) \,
	(1 \,+\, q \, 2^{q-1} \,-\, 2^q) \\
\Longrightarrow~~~~\hfil{}
K_6
&=&	d(p) \, d(q).
\end{eqnarray*}%
}%

\noindent
Edges $e \in E_8$ have all their vertices in $V_3$.
Their dot product term is $(i - 1)^2$, where $i=|e|$:
\begin{eqnarray*}
K_8
&=&	\sum_{i=2}^r {r \choose i} (i \,-\, 1)^2 \\
K_8
&=&	\alpha(r).
\end{eqnarray*}

\noindent
Edges $e \in E_{10}$ have vertices in $V_3$ and $V_1$.
Let $e \in E_{10}$,
$i = |e \cap V_1|$, and
$j = |e \cap V_3|$.
The corresponding dot-product term is $(i + j - 1) \, (j - 1)$:
{\allowdisplaybreaks%
\begin{eqnarray*}
K_{10}
&=&	\sum_{i=1}^p {p \choose i}
	 \sum_{j=1}^r {r \choose j}
		(i \,+\, j \,-\, 1) \, (j \,-\, 1) \\
&=&	\sum_{i=1}^p {p \choose i}
	 \left \lbrack
	  \sum_{j=1}^r {r \choose j}
		(j \,-\, 1)^2
	  \,+\, i
	  \sum_{j=1}^r {r \choose j}
		(j \,-\, 1)
	 \right \rbrack \\
&=&	\sum_{i=1}^p {p \choose i}
	 \left \lbrack
	  \alpha(r)
	  \,+\, i \,	\left (
				1
				\,+\,
				\sum_{j=0}^r {r \choose j}
					(j \,-\, 1)
			\right )
	 \right \rbrack \\
&=&	\sum_{i=1}^p {p \choose i}
	 \left \lbrack
	  \alpha(r)
	  \,+\, i \, (1 \,+\, r \, 2^{r-1} \,-\, 2^r)
	 \right \rbrack \\
&=&	\alpha(r) \,
	\sum_{i=1}^p {p \choose i}
	\,+\,
	(1 \,+\, r \, 2^{r-1} \,-\, 2^r) \,
	\sum_{i=0}^p {p \choose i} i \\
&=&	\alpha(r) \,
	(2^p \,-\, 1)
	\,+\,
	(1 \,+\, r \, 2^{r-1} \,-\, 2^r) \, p \, 2^{p-1} \\
&=&	(2^p \,-\, 1) \,
	\alpha(r)
	\,+\,
	p \,
	2^{p-1} \,
	(1 \,+\, r \, 2^{r-1} \,-\, 2^r) \\
\Longrightarrow~~~~\hfil{}
K_{10}
&=&	c(p) \,	\alpha(r) \,+\, b(p) \, d(r).
\end{eqnarray*}%
}%

\noindent
Edges $e \in E_{12}$ have vertices in $V_3$ and $V_2$.
Let $e \in E_{12}$,
$i = |e \cap V_2|$, and
$j = |e \cap V_3|$.
The corresponding dot-product term is
$(j - 1) \, (i + j - 1)$:
{\allowdisplaybreaks%
\begin{eqnarray*}
K_{12}
&=&	\sum_{i=1}^q {q \choose i}
	 \sum_{j=1}^r {r \choose j}
		(j \,-\, 1) \, (i \,+\, j \,-\, 1) \\
&=&	\sum_{i=1}^q {q \choose i}
	 \left \lbrack
	  \sum_{j=1}^r {r \choose j}
		(j \,-\, 1)^2
	  \,+\, i \,	\sum_{j=1}^r {r \choose j}
				(j \,-\, 1)
	 \right \rbrack \\
&=&	\sum_{i=1}^q {q \choose i}
	 \left \lbrack
	  \alpha(r)
	  \,+\, i \,	\left (1 \,+\,
				\sum_{j=0}^r {r \choose j}
					(j \,-\, 1)
			\right )
	 \right \rbrack \\
&=&	\sum_{i=1}^q {q \choose i}
	 \left \lbrack
	  \alpha(r)
	  \,+\, i \, (1 \,+\, r \, 2^{r-1} \,-\, 2^r)
	 \right \rbrack \\
&=&	\alpha(r) \,
	\sum_{i=1}^q {q \choose i}
	\,+\,
	(1 \,+\, r \, 2^{r-1} \,-\, 2^r) \,
	\sum_{i=0}^q {q \choose i} i \\
&=&	\alpha(r) \,
	(2^q \,-\, 1)
	\,+\,
	(1 \,+\, r \, 2^{r-1} \,-\, 2^r) \, q \, 2^{q-1} \\
&=&	(2^q \,-\, 1) \, \alpha(r)
	\,+\,
	q \, 2^{q-1} \, (1 \,+\, r \, 2^{r-1} \,-\, 2^r) \\
\Longrightarrow~~~~\hfil{}
K_{12}
&=&	c(q) \, \alpha(r) \,+\, b(q) \, d(r).
\end{eqnarray*}%
}%

\noindent
Edges $e \in E_{14}$ have vertices in $V_3$, $V_2$ and $V_1$.
Let $e \in E_{14}$,
$i = |e \cap V_1|$,
$j = |e \cap V_2|$, and
$k = |e \cap V_3|$.
The corresponding dot-product term is
$(i \,+\, k \,-\, 1) \, (j \,+\, k \,-\, 1)$:
{\allowdisplaybreaks%
\begin{eqnarray*}
K_{14}
&=&	\sum_{i=1}^p {p \choose i}
	 \sum_{j=1}^q {q \choose j}
	  \sum_{k=1}^r {r \choose k}
		(i \,+\, k \,-\, 1) \, (j \,+\, k \,-\, 1) \\
&=&	\sum_{i=1}^p {p \choose i}
	 \sum_{j=1}^q {q \choose j}
	  \left \lbrack
	   \sum_{k=1}^r {r \choose k}
		(k - 1)^2
	   +
	   (i + j)
	   \sum_{k=1}^r {r \choose k}
		(k - 1)
	  \right . \\
&&	  \left .
	~~~~~~~~~~~~~~~~~~~~~~~~
	   +
	   i j
	   \sum_{k=1}^r {r \choose k}
	  \right \rbrack \\
&=&	\sum_{i=1}^p {p \choose i}
	 \sum_{j=1}^q {q \choose j}
	  \left \lbrack
	   \alpha(r)
	   \,+\,
	   (i \,+\, j) \,
	   (1 \,+\, r \, 2^{r-1} \,-\, 2^r)
	   \,+\,
	   i \, j\,
	   (2^r \,-\, 1)
	  \right \rbrack \\
&=&	\sum_{i=1}^p {p \choose i}
	 \left \lbrack
	  (\alpha(r)
	   \,+\,
	   i \, (1 \,+\, r \, 2^{r-1} \,-\, 2^r))
	  \sum_{j=1}^q {q \choose j}
	 \right \rbrack \\
&&	 \,+\,
	\sum_{i=1}^p {p \choose i}
	 \left \lbrack
	   (1 \,+\, r \, 2^{r-1} \,-\, 2^r
	    \,+\,
	    i \, (2^r \,-\, 1)) \,
	  \sum_{j=1}^q {q \choose j} j
	 \right \rbrack \\
&=&	\sum_{i=1}^p {p \choose i}
	 \left \lbrack
	  (\alpha(r)
	   \,+\,
	   i \, (1 \,+\, r \, 2^{r-1} \,-\, 2^r)) \,
	  (2^q \,-\, 1)
	 \right \rbrack \\
&&	 \,+\,
	\sum_{i=1}^p {p \choose i}
	 \left \lbrack
	   (1 \,+\, r \, 2^{r-1} \,-\, 2^r
	    \,+\,
	    i \, (2^r \,-\, 1)) \,
	  (q \, 2^{q-1})
	 \right \rbrack \\
&=&	\left \lbrack
		(\alpha(r) \, (2^q \,-\, 1)
		\,+\,
		q \, 2^{q-1} \,
		(1 \,+\, r \, 2^{r-1} \,-\, 2^r))
		\sum_{i=1}^p {p \choose i}
	\right \rbrack \\
&&	\,+\,
	\left \lbrack
		((2^q \,-\, 1) \, (1 \,+\, r\,2^{r-1} \,-\, 2^r)
		\,+\,
		q \, 2^{q-1} \,
		(2^r \,-\, 1))
		\sum_{i=1}^p {p \choose i} i
	\right \rbrack \\
&=&	\left \lbrack
		(\alpha(r) \, (2^q \,-\, 1)
		\,+\,
		q \, 2^{q-1} \,
		(1 \,+\, r \, 2^{r-1} \,-\, 2^r)) \,
		(2^p \,-\, 1)
	\right \rbrack \\
&&	\,+\,
	\left \lbrack
		((2^q \,-\, 1) \, (1 \,+\, r\,2^{r-1} \,-\, 2^r)
		\,+\,
		q \, 2^{q-1} \,
		(2^r \,-\, 1)) \,
		p \, 2^{p-1}
	\right \rbrack \\
&=&	(p \, 2^{p-1} \, (2^q \,-\, 1)
	 \,+\, q \, 2^{q-1} \, (2^p \,-\, 1)) \,
	(1 \,+\, r\,2^{r-1} \,-\, 2^r) \\
&&
	~~~~
	\,+\,
	p \, 2^{p-1} \, q \, 2^{q-1} \, (2^r \,-\, 1)
	\,+\,
	(2^p \,-\, 1) \, (2^q \,-\, 1) \, \alpha(r) \\
\Longrightarrow~~~~\hfil{}
K_{14}
&=&	(b(p) \, c(q) \,+\, b(q) \, c(p)) \, d(r)
	\,+\, b(p) \, b(q) \, c(r)
	\,+\, c(p) \, c(q) \, \alpha(r).
\end{eqnarray*}%
}%

\noindent
We add the partial sums:
{\allowdisplaybreaks%
\begin{eqnarray*}
&&	K_6 \,+\, K_8 \,+\, K_{10} \,+\, K_{12} \,+\, K_{14} \\
&=&
	d(p) \, d(q)
	\,+\, \alpha(r)
	\,+\, c(p) \, \alpha(r) \,+\, b(p) \, d(r)
	\,+\, c(q) \, \alpha(r) \,+\, b(q) \, d(r) \\
&&	~~
	\,+\,
	(b(p) \, c(q) \,+\, b(q) \, c(p)) \, d(r)
	\,+\, b(p) \, b(q) \, c(r) \,+\, c(p) \, c(q) \, \alpha(r) \\
&=&
	(1 \,+\, c(p) \,+\, c(q) \,+\, c(p) \, c(q)) \, \alpha(r)
	\,+\, d(p) \, d(q) \\
&&	~~
	\,+\, b(p) \, b(q) \, c(r)
	\,+\,	(b(p) \,+\, b(q)
		 \,+\, b(p) \, c(q) \,+\, b(q) \, c(p)
		) \, d(r) \\
&=&
	(2^{p+q}) \, \alpha(r)
	\,+\, d(p) \, d(q)
	\,+\, b(p) \, b(q) \, c(r)
	\,+\,	b(p+q) \, d(r).
\end{eqnarray*}%
}%
This yields:
\begin{align*}
a_1 \cdot a_2 &=
	2^{n-p-q-r} \,
	\left \lbrack
		2^{p+q}\alpha(r)
		\,+\, d(p)\,d(q)
		\,+\, b(p)\,b(q)\,c(r)
		\,+\,b(p+q)\,d(r)
	\right \rbrack \\
&= w(n,p,q,r).
\tag*{\qed}
\end{align*}
\end{MyProof}

\noindent
We now examine the special case of the angle between complementary
subtours $S$ and $V-S$.
Let $n = |V|$ and $k = |S|$.
This case is given by
$p=k$, $q=n-k$ and $r=0$,
for which we observe that
$\alpha(r) = c(r) = d(r) = 0$, giving
\begin{eqnarray*}
h(n,k)
&=&	f(n,k,n-k,0) \\
&=&	{{	\alpha(n) \, d(k) \, d(n-k)
		\,-\, \beta(n,k) \, \beta(n,n-k)}
	\over
	{\sqrt{\mu(n,k) \, \mu(n,n-k)}}}
\end{eqnarray*}

\noindent
We have
\begin{displaymath}
	\lim_{n \to \infty} h(2n, n)
 ~=~ \lim_{n \to \infty}
	{{-n^2 \, 2^{4n \,+\, 1} \,+\, O(2^{4n \,-\, 1})}
	 \over
	 { n^2 \, 2^{4n \,+\, 1} \,+\, O(2^{4n \,-\, 1})}}
 ~=~ -1
\end{displaymath}
which means the angle $\phi$ between the complementary subtour normal
vectors $\hat{a_1}$ and $\hat{a_2}$ goes to $\pi$ and the interior
angle $\theta$ goes to zero.
This has important implications upon the condition number of the LP
basis matrix when both of these subtours are simultaneously binding.

We note that calculating such angles does not require enumerating the
extreme points of the polyhedron, only calculating dot products
regarding the two inequalities (and the affine hull when present).

\subsection{Large and Small Subtours are Effective Together}
\label{sec:larg+small=strong}

We now explore reasons to expect that combining subtours of both
small and large cardinality should be very strong in practice.

Let $S,T \subset V$, $k_1=|S|$ and $k_2=|T|$ such that
$S \cap T = \emptyset$,
$S$ is a violated subtour,
$T$ is a violated anti-subtour,
with $k_1$ and $k_2$ similarly small.
Subtour $x(S) \le |S|-1$ places an upper bound on the total weight of
edges within $S$, while anti-subtour $x(T) + x(T:V-T) \ge |T|$ places
a lower bound on the total weight of edges within $T$ and crossing cut
$(T:V-T)$.
Together with equation $x(V) = |V|-1$ (which serves as a
``conservation of edge weight'' constraint), combining these two
inequalities together forces the surplus of edge weight within $S$ to
directly offset the deficit of edge weight at $T$ and vice versa.
Only the difference between these violations (positive or negative)
affects the rest of the problem.

In contrast, the ``push'' action of small subtours alone spreads the
surplus outward, generally forming larger, weaker violated subtour
violations.
Furthermore, using only small subtours, it can take many such ``push''
operations to finally satisfy a physically distant edge weight
deficit.

Further consider strength obtained by the angle at which subtour $S$
and anti-subtour $T$ meet.
Figure~\ref{fig:subtour-angles-5} shows this angle is just over 90
degrees, and increases as $k_1$ and $k_2$ increase.
Although this does not provide substantial angular strength, we still
have two individually very strong inequalities, strongly coupled.

The angular component substantially improves, however, when $S$ and
$T$ share vertices.
The case where $S=T$ (a complementary subtour) is particularly
strong:
EPR and CD indicate that both subtours are individually very
strong, and Figure~\ref{fig:complementary-subtour-angles} shows the
angle between them is also very small, making this combination
{\em triply} strong.

Consider also case where
$|S| \,=\, 2k$,
$|T| \,=\, 2k$, and
$|S \cap T| = k$,
for small $k$.
This gives subtours of size $2k$ and $n-2k$ which EPR and CD indicate
are both very strong.
The angle between them, is $\theta(n,\, k,\, n-3k,\, k)$ as
illustrated in Figure~\ref{fig:subtour-angles-6},
which are just over 60 degrees, which is moderately strong.

Similar combinations of two small subtours do not provide similar
angular strength.
As shown in Figure~\ref{fig:subtour-angles-4}, two independent
subtours of size $k$ yield angles just under 90 degrees.
As shown in Figure~\ref{fig:subtour-angles-7}, a pair of subtours of
size $2k$ that overlap by $k$ vertices yield angles just below 120
degrees.
(These angles go as low as 60 degrees, but only with very weak EPR and
CD indicators.)

The interaction of subtours with large and small cardinality is
therefore predicted to be significantly stronger than pairs of
subtours having small cardinality.

\section{Spanning Tree in Graph Polytope: $STGP(n)$}
\label{sec:stgp-polytope}

Throughout this section, let $K_n = G = (V,E)$ be a complete graph,
$n = |V|$,
$m = |E| = {n \choose 2}$,
and $n \ge 2$.

The spanning tree in graph polytope is characterized as follows:
\begin{eqnarray}
\nonumber
	x(V) &=&	|V| \,-\, 1, \\
\label{eq:spanning-tree-subtour}
	x(S) &\le&	|S| \,-\, 1
	~~~~~~~\hbox{$\forall S \subset V$, $|S| \ge 2$}, \\
\nonumber
	x_e &\ge& 0 ~~~~\hbox{$\forall e \in E$}.
\end{eqnarray}
It is well-known that the extreme points of this set are the integral
0-1 incidence vectors of all spanning trees in $K_n$.

\begin{theorem}[Cayley's formula]
The number of distinct spanning trees of $K_n$ is
\begin{displaymath}
	n^{n-2}.
\end{displaymath}
\end{theorem}
This is a classic result often attributed to Cayley~\cite{Cayley},
even though Cayley's paper references this result in an earlier paper
by Borchardt~\cite{Borchardt}.
There is reason to believe Kirchoff may have been aware of this
formula much earlier (circa 1840), though no actual evidence has yet
surfaced.

\subsection{EPR for STGP Subtours}

\begin{theorem}
Let $n \ge 3$,
$S \subset V$ and $k = |S|$ such that $k \ge 2$.
The number of spanning trees of $K_n$ incident to subtour $S$ is:
\begin{displaymath}
	\mathrm{ni}(n,k) ~=~ k^{k \,-\, 1} \, n^{n \,-\, k \,-\, 1}.
\end{displaymath}
\end{theorem}

\begin{MyProof}
Every spanning tree incident to subtour $S$ forms a tree within $S$.
By Caley's formula, there are $k^{k-2}$ distinct such trees.
We partition the remaining $n-k$ vertices into $k$ bins, the $i$-th
bin containing $a_i$ vertices for $1 \le i \le k$.
For each such vector $(a_1,\, a_2,\, \ldots,\, a_k)$ there are
\begin{displaymath}
	\left ( {\MyAtop{n \,-\, k}%
		 {a_1,\, a_2,\, \ldots,\, a_k}}
	\right )
\end{displaymath}
ways of choosing said partition.
For each bin $i$ we construct a spanning tree consisting of the
$i$-th vertex of $S$, plus the vertices in bin $i$.
By Caley's formula, there are $(a_i \,+\,1)^{a_i \,-\,1}$ distinct
such trees.
Let $b_i \,=\, a_i \,+\, 1$ for all $i$.
The total number of spanning trees incident to subtour $S$ is
therefore
{\allowdisplaybreaks%
\begin{eqnarray}
\nonumber
\mathrm{ni}(n,k) &=&
	k^{k-2}
	\sum_{\MyAtop{a_1,\, a_2,\, \ldots,\, a_k \,\ge\, 0}%
	      {a_1 + a_2 + \cdots + a_k = n-k}}
	\left ( {\MyAtop{n \,-\, k}%
		 {a_1, a_2, \ldots, a_k}}
	\right )
	\prod_{i=1}^k (a_i + 1)^{a_i-1} \\
\nonumber
&=&	k^{k-2} \,
	\sum_{\MyAtop{a_1,\, a_2,\, \ldots,\, a_k \,\ge\, 0}%
	      {a_1 \,+\, a_2 \,+\, \cdots \,+\, a_k \,=\, n-k}}
	{{(n \,-\, k)!} \over
	 {a_1!\, a_2!\, \ldots \, a_k!}} \,
	\prod_{i=1}^k (a_i \,+\, 1)^{a_i-1} \\
\nonumber
&=&	k^{k-2} \,
	{{(n \,-\, k)!} \over {n!}} \\
\nonumber
&&~~~~~
	\sum_{\MyAtop{a_1, \ldots, a_k \ge 0}%
	      {a_1 + \cdots + a_k = n-k}}
	{{n!} \over
	 {(a_1 + 1)! \, \ldots \, (a_k + 1)!}} \,
	\prod_{i=1}^k (a_i + 1)^{a_i} \\
\nonumber
&=&	k^{k-2} \,
	{{(n \,-\, k)!} \over {n!}} \,
	\sum_{\MyAtop{b_1,\, b_2,\, \ldots,\, b_k \,\ge\, 1}%
	      {b_1 \,+\, b_2 \,+\, \cdots \,+\, b_k \,=\, n}}
	{{n!} \over  {b_1!\, b_2!\, \ldots \, b_k!}} \,
	\prod_{i=1}^k b_i^{b_i-1} \\
\label{eq:st-epr-proof-form-1}
&=&	k^{k-2} \,
	{{(n - k)!} \over {n!}} \,
	\sum_{\MyAtop{b_1, b_2, \ldots, b_k \ge 1}%
	      {b_1 + b_2 + \cdots + b_k = n}}
	\left ( {\MyAtop{n}{b_1, b_2, \ldots, b_k}}
	\right )
	\prod_{i=1}^k b_i^{b_i-1}.
\end{eqnarray}%
}%
It is known~\cite{Wilf} that the number of {\em rooted} spanning trees
in a graph is $n^{n-1}$ and its exponential generating function $H(z)$
satisfies
\begin{displaymath}
	H(z) ~=~ z \, e^{H(z)}.
\end{displaymath}
We recognize the summation of (\ref{eq:st-epr-proof-form-1}) to
be
\begin{displaymath}
	\lbrack {{z^n} \over {n!}} \rbrack H(z)^k
\end{displaymath}
so that
\begin{displaymath}
\mathrm{ni}(n,k) ~=~
	k^{k-2} \,
	{{(n \,-\, k)!} \over {n!}} \,
	\lbrack {{z^n} \over {n!}} \rbrack \, H(z)^k.
\end{displaymath}
We now invoke the Lagrange inversion formula (Lemma~\ref{lem:lif})
with
$f(u) \,=\,u^k$ and $\theta(u) \,=\, e^u$:
{\allowdisplaybreaks%
\begin{align*}
\mathrm{ni}(n,k) &=
	k^{k-2} \,
	{{(n \,-\, k)!} \over {n!}} \,
	\lbrack {{z^n} \over {n!}} \rbrack \, H(z)^k \\
&=
	k^{k-2} \,
	(n \,-\, k)! \,
	\lbrack z^n \rbrack \, H(z)^k \\
&=
	k^{k-2} \,
	(n \,-\, k)! \,
	\lbrack z^n \rbrack \, f(u(z)) \\
&=
	k^{k-2} \,
	(n \,-\, k)! \,
	{1 \over n} \,
	\lbrack u^{n-1} \rbrack \, f'(u) \, \theta(u)^n) \\
&=
	k^{k-2} \,
	(n \,-\, k)! \,
	{1 \over n} \,
	\lbrack u^{n-1} \rbrack \, k\,u^{k-1} \, e^{n \, u} \\
&=
	k^{k-1} \,
	(n \,-\, k)! \,
	{1 \over n} \,
	\lbrack u^{n-k} \rbrack \, e^{n \, u} \\
&=
	k^{k-1} \,
	(n \,-\, k)! \,
	{1 \over n} \,
	\lbrack u^{n-k} \rbrack \,
	\sum_{i \ge 0} {{n^i \, u^i} \over {i!}} \\
&=
	k^{k-1} \,
	(n \,-\, k)! \,
	{1 \over n} \,
	{{n^{n-k}} \over {(n \,-\, k)!}} \\
&=
	k^{k-1} \,
	n^{n-k-1}.
\tag*{\qed}
\end{align*}%
}%
\end{MyProof}

\begin{theorem}
Let $n \ge 3$, $S \subset V$ and $k = |S|$ such that $k \ge 2$.
EPR of subtour $S$
(equation (\ref{eq:spanning-tree-subtour})) is:
\begin{displaymath}
	\mathrm{EPR}(n,k) ~=~ (k \,/\, n)^{k-1}. 
\end{displaymath}
\end{theorem}

\begin{MyProof}
EPR is the ratio of the number of incident
spanning trees to the total number of spanning trees
\begin{displaymath}
\mathrm{EPR(n,k)}
	= {{\mathrm{ni}(n,k)} \over {n^{n-2}}}
	= {{k^{k-1} \, n^{n-k-1}} \over {n^{n-2}}}
	= (k \,/\, n)^{k-1}.
\tag*{\qed}
\end{displaymath}
\end{MyProof}

\subsection{Centroid of $\STGP{n}$}

\begin{theorem}
Let $n \ge 3$.  The centroid $C$ of the Spanning Tree in Graph
Polytope $\mathrm{STGP}(n)$ satisfies
\begin{displaymath}
	C_e ~=~ {2 \over n} ~~~~~\forall e \in E.
\end{displaymath}
\end{theorem}

\begin{MyProof}
Let $e \in E$.
The number of spanning trees containing edge $e$ is just
$\mathrm{ni}(n,2) ~=~ 2 \, n^{n-3}$.
Dividing by the total number of trees gives
\begin{displaymath}
C_e	= {{2 \, n^{n-3}} \over {n^{n-2}}}
	= {2 \over n}.
\tag*{\qed}
\end{displaymath}
\end{MyProof}

\subsection{CD of STGP Subtour Inequalities}

\begin{theorem}
Let $n \ge 3$, $S \subset V$ and $k = |S|$ such that $k \ge 2$.
CD squared of subtour $S$
(equation (\ref{eq:spanning-tree-subtour})) is:
\begin{displaymath}
\mathrm{CD}(n,k)^2 ~=~
	{{2 \, (k \,-\, 1) \, (n \,-\, 1) \, (n \,-\, k)} \over
	 {k \, n \, (n \,+\, k \,-\, 1)}}.
\end{displaymath}
\end{theorem}

\begin{MyProof}
We use the distance formula squared for polyhedra having a single
affine hull equation~(\ref{eq:general-centroid-distance-dim-1}).
\begin{displaymath}
	\mathrm{dist}^2 ~=~ {{s \, (b_1 \,-\, a \cdot C)^2} \over
			     {q \, s \,-\, r^2}},
\end{displaymath}
where in this case
$q = {k \choose 2}$,
$r = {k \choose 2}$,
$s = {n \choose 2}$,
$b_1 = k - 1$,
and $a \cdot C = {2 \over n} \, {k \choose 2}$.
This readily simplifies to
\begin{equation}
\label{eq:stgp-subtour-cdist}
\mathrm{cdist}(n,k)^2 ~=~
	{{2 \, (k \,-\, 1) \, (n \,-\, 1) \, (n \,-\, k)} \over
	 {k \, n \, (n \,+\, k \,-\, 1)}}.
\end{equation}
\qed
\end{MyProof}

\subsection{Unexpected Behavior of CD for STGP Subtours}
\label{sec:unexpected-stgp-cd-behavior}

\begin{figure}[!ht]
 \begin{center}
  \begin{minipage}[t]{2.25in}
   \includegraphics[width=2.25in,clip=]{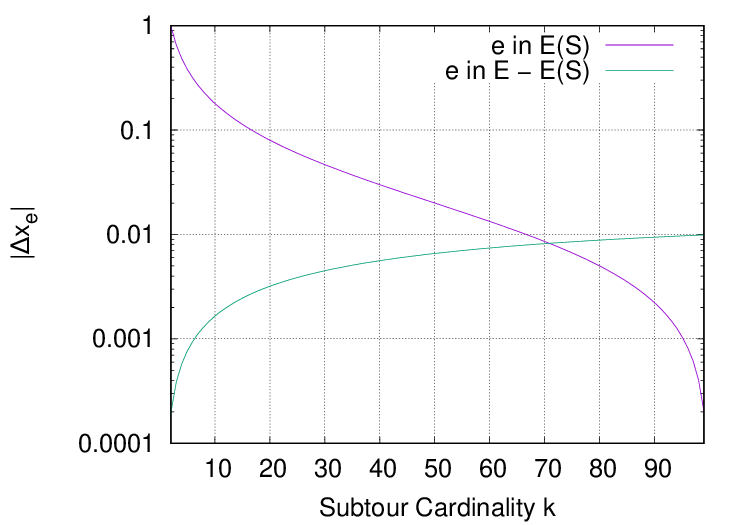}
   \caption{$\STGP{n}$ subtour CD $|\Delta x_e|$ values as function of
     subtour size $k$, for $n=100$}
   \label{fig:stgp-dx-funcs}
  \end{minipage}
  \quad
  \begin{minipage}[t]{2.25in}
   \includegraphics[width=2.25in,clip=]{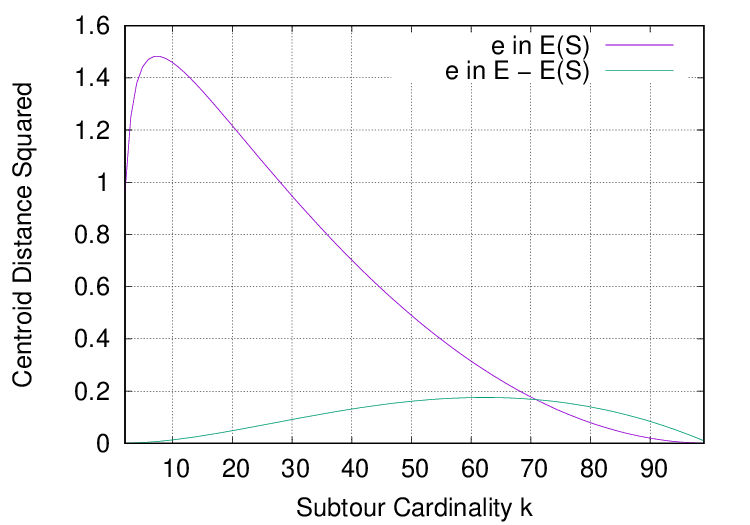}
   \caption{$\STGP{n}$ subtour CD squared summed separately over edges
	$e \in E(S)$ and $e \in E - E(S)$, as functions of subtour
	size $k$, for $n=100$}
   \label{fig:stgp-dist2-funcs}
  \end{minipage}
 \end{center}
\end{figure}

We now investigate:
(1) why CD is so large for STGP subtours of small cardinality, and
(2) why CD is so different between STHGP and STGP.
The behavior of the STGP
formula~(\ref{eq:stgp-subtour-cdist})
itself is clear from:
(a) the $2(n-1)/n$ factors do not involve $k$ and approach the
constant $2$ as $n$ increases;
(b) the factors $(k-1)/k$ are also relatively constant, staying in the
interval $[ 1/2, 1)$.
The dominant behavior therefore arises from the factors
\begin{displaymath}
	{{n \,-\, k} \over {n \,+\, k \,-\, 1}}
\end{displaymath}
for which the numerator decreases linearly with $k$, while the
denominator increases linearly with $k$.
We would like to understand the geometry of how this arises within
polytope STGP, however.
We use the various results leading up
to~(\ref{eq:general-centroid-distance-dim-1}).
The vector between closest point $p$ and centroid $C$ is
\begin{displaymath}
\Delta x
	~=~	p \,-\, C
	~=~	\tau \, \hat{a}
\end{displaymath}
where in the present case
\begin{displaymath}
	\tau	~=~ {4 \over {k \, n \, (n \,+\, k \,-\, 1)}},
	~~~~
	\hat{a}	~=~ {n \choose 2} a \,-\, {k \choose 2} c
\end{displaymath}
where $a$ is the left-hand-side of the subtour and $c$ is the
left-hand-side of the affine hull equation.
This gives
\begin{displaymath}
	\Delta x_e ~=~
	\begin{cases}
		{{2 \, (n \,-\, k)} \over {k \, n}} &
		 \text{for all $e \in E(S)$;} \\
		\\
		{{- \, 2 \, (k \,-\, 1)} \over {n \, (n \,+\, k \,-\, 1)}} &
		 \text{for all $e \in E - E(S)$.}
	\end{cases}
\end{displaymath}
Figure~\ref{fig:stgp-dx-funcs} plots these two $|\Delta x_e|$
coordinate
values for $n=100$ as functions of subtour size $k$ with a log scale.
This plot clearly shows that as $k$ increases, the $E(S)$ edge
components get closer to $C$, whereas the components for other edges
grow more distant.
We also have $\mathrm{cdist}(n,k)^2 = \sum_{e \in E} {\Delta x}_e^2$.
Splitting this into two sums, one over $e \in E(S)$ and the other over
$e \in E \,-\, E(S)$ yields:
\begin{align*}
	{{2 \, (k \,-\, 1) \, (n \,-\, k)^2} \over
	 {k \, n^2}}
	~~~~~~~~~~
	& \text{for sum over edges $e \in E(S)$}, \\
	{{2 \, (k \,-\, 1)^2 \, (n \,-\, k)} \over
	 {n^2 \, (n \,+\, k \,-\, 1)}}
	~~~~~~~~~~
	& \text{for sum over edges $e \in E \,-\, E(S)$}.
\end{align*}
Figure~\ref{fig:stgp-dist2-funcs} plots these two partial sums for
$n=100$ as functions of subtour size $k$.
This clearly shows the edges $e \in E(S)$ are responsible for most
of the large CD seen for small $k$.

To shed some light on why CD behaves so differently
between subtours of $\STHGP{n}$ and $\STGP{n}$, we now
repeat this computation for subtours of $\STHGP{n}$, for which we
have:
\begin{displaymath}
\Delta x ~=~ p \,-\, C ~=~ \tau \, \hat{a},
\end{displaymath}
where (using (6.30))
\begin{align*}
\tau	&= {\MyFrac{b \,-\, a \cdot C}%
		   {q \, s \,-\, r^2}}
	~=~ {\MyFrac{(n \,-\, k) \, t(n,k)}%
	     {(\alpha(n) \, \gamma(n,k) \,-\, \beta^2(n,k))
	      \, E[X_n^{n-1}]}}
	~=~ {\MyFrac{(n \,-\, k) \, t(n,k)}%
	     {\mu(n,k) \, E[X_n^{n-1}]}}, \\
\hat{a}	&= \alpha(n) \, a \,-\, \beta(n,k) \, c.
\end{align*}

Let $n=|V|$, $S \subset V$ such that $k=|S| \ge 2$, $e \in E$,
$p = |e \cap S|$ and $q = |e \cap (V-S)| = |e|-p$.
Then
\begin{align*}
\Delta x	&= \tau \, (\alpha(n) \, a \,-\, \beta(n,k) \, c) \\
\Delta x_e	&= \Delta x(n,\, k,\, p,\, q) \\
		&= \tau \, \left \lbrack
			 \alpha(n) \, \max(p \,-\, 1,\, 0)
			 \,-\,
			 \beta(n,k) \, \max(p \,+\, q \,-\, 1,\, 0)
		   \right \rbrack.
\end{align*}
For $n = 100$, the minimum and maximum values of $|\Delta x(n,k,p,q)|$
occur at:
\begin{eqnarray*}
	|\Delta x(100, 99, 98, 1)| &=& 4.9135920550309504 \,\times\, 10^{-34} \\
	|\Delta x(100, 2, 100, 0)| &=& 4.0791412569526042 \,\times\, 10^{-28}
\end{eqnarray*}
so that $|\Delta x_e|$ values vary by a factor of at most 830175.
The partial sum of CD squared for all edges $e \in E(S)$ is:
{\allowdisplaybreaks
\begin{flalign}
\nonumber
&\sum_{e \in E(S)}
	(\Delta x_e)^2 \\
\nonumber
&= \sum_{p=2}^k {k \choose p} \sum_{q=0}^{n-k} {{n-k} \choose q}
	\Delta x(n,\, k,\, p,\, q)^2 \\
\nonumber
&= \tau^2 \sum_{p=2}^k {k \choose p} \sum_{q=0}^{n-k} {{n-k} \choose q}
	(\alpha(n) \, (p-1) \,-\, \beta(n,k) \,(p+q-1))^2 \\
\nonumber
&= \tau^2 \sum_{p=2}^k {k \choose p}
	\left \lbrack
	\beta^2(n,k) \sum_{q=0}^{n-k} {{n-k} \choose q} \, q^2
	\right . \\
\nonumber
&	\left .
	~~~~~~~~~~~~~~~~~~
	 \,+\,
	 2 \, (p-1) \,\beta(n,k) \, (\beta(n,k) \,-\, \alpha(n))
		\sum_{q=0}^{n-k} {{n-k} \choose q} \, q
	\right . \\
\nonumber
&	\left .
	~~~~~~~~~~~~~~~~~~
	\,+\,
	 (p-1)^2 \, (\alpha(n) \,-\, \beta(n,k))^2
		\sum_{q=0}^{n-k} {{n-k} \choose q}
	\right \rbrack \\
\nonumber
&= \tau^2 \sum_{p=2}^k {k \choose p}
	\left \lbrack
	\beta^2(n,k) (n-k) \, (n-k+1) \, 2^{n-k-2}
	\right . \\
\nonumber
&	~~~~~~~~~~~~~~~~~~
	 \,+\,
	 (p-1) \, \beta(n,k) \, (\beta(n,k) \,-\, \alpha(n))
		\, (n-k) \, 2^{n-k} \\
\nonumber
&	\left .
	~~~~~~~~~~~~~~~~~~
	\,+\,
	 (p-1)^2 \, (\alpha(n) \,-\, \beta(n,k))^2 \, 2^{n-k}
	\right \rbrack \\
\nonumber
&= \tau^2 \, (\alpha(n) \,-\, \beta(n,k))^2 \, 2^{n-k}
	\sum_{p=2}^k {k \choose p} \, p^2 \\
\nonumber
&~~~~ -\, \tau^2
	\, (\alpha(n) \,-\, \beta(n,k))
	\, ((n \,-\, k \,-\, 2) \, \beta(n,k) \,+\, 2 \, \alpha(n))
	\, 2^{n-k}
	\sum_{p=2}^k {k \choose p} \, p \\
\nonumber
&~~~~ +\, \tau^2 \,
	\left \lbrack
		(n \, (n \,-\, 2k \,-\, 3) \,+\, k^2 \,+\, 3k \,+ 4) \, \beta^2(n,k)
        \right . \\
\nonumber
&	\left .
~~~~~~~~~~~~
		+\, 4\,(n \,-\, k \,-\, 2) \, \alpha(n) \, \beta(n,k)
		\,+\, 4 \, \alpha^2(n)
	\right \rbrack
	\, 2^{n-k-2}
	\sum_{p=2}^k {k \choose p} \\
\nonumber
&= \tau^2 (\alpha(n) \,-\, \beta(n,k))^2 \, 2^{n-k}
	k \, ((k + 1) \, 2^{k-2} \,-\, 1) \\
\nonumber
&~~~~ -\, \tau^2 \,
	\left \lbrack
		(\alpha(n) \,-\, \beta(n,k))
		\, (n \,-\, k \,-\, 2) \, \beta(n,k)
		\,+\, 2 \, \alpha(n)
	\right \rbrack
	\, 2^{n-k}
	\, k \, (2^{k-1} \,-\, 1) \\
\nonumber
&~~~~ +\, \tau^2 \,
	\left \lbrack
		(n \, (n \,-\, 2k \,-\, 3) \,+\, k^2 \,+\, 3k \,+\, 4)
		   \, \beta^2(n,k)
	\right . \\
\nonumber
&	\left .
~~~~~~~~~~~~
		\,+\, 4 \, (n \,-\, k \,-\ 2) \, \alpha(n) \, \beta(n,k)
		\,+\, 4 \, \alpha^2(n)
	\right \rbrack
	\, 2^{n-k-2}
	\, (2^k \,-\, k \,-\, 1) \\
\label{eq:sthgp-dx-sum-E(S)}
&= \tau^2 \, \left \lbrack
	((k^2 \,-\, 3k \,+\, 4) \, 2^{n-2} \,-\, 2^{n-k})
	\, \alpha^2(n)
	\right . \\
\nonumber
&
~~~~~~~~
	\,+\, (((-k \,+\, 2) \, n \,+\, k \,-\, 4) \, 2^{n-1}
	       \,-\, (n \,-\, k \,-\, 2) \, 2^{n-k})
	      \, \alpha(n) \, \beta(n,k) \\
\nonumber
&
~~~~~~~~
	  \,+\, ((n^2 \,-\, 3n \,+\, 4) \, 2^{n-2} \\
\nonumber
&	\left .
~~~~~~~~~~~~
		 \,-\,
		 (k^3 \,-\, 2 n k^2 \,+\, (n^2 \,-\, n \,+\, 3)\, k
		  \,+\, n^2 \,-\, 3n \,+\, 4) \, 2^{n-k-2})
	 \, \beta^2(n.k)
	\right \rbrack
\end{flalign}
}
The partial sum of CD squared for all edges $e \in (E - E(S))$ is:
{\allowdisplaybreaks
\begin{flalign}
\nonumber
&\sum_{e \in (E - E(S))}
	(\Delta x_e)^2 \\
\nonumber
&=	\sum_{p=0}^1 {k \choose p} \sum_{q=0}^{n-k} {{n-k} \choose q}
		\Delta x(n,\, k,\, p,\, q)^2 \\
\nonumber
&=	\tau^2 \sum_{q=0}^{n-k} {{n-k} \choose q}
		(\beta(n,k) \max(q-1, 0))^2
	+ k \tau^2 \sum_{q=0}^{n-k} {{n-k} \choose q}
		(\beta(n,k) \max(q, 0))^2 \\
\nonumber
&=	\tau^2 \, \beta^2(n,k) \sum_{q=2}^{n-k} {{n-k} \choose q}
		(q-1)^2
	\,+\, k \, \tau^2 \, \beta^2(n,k) \sum_{q=0}^{n-k} {{n-k} \choose q}
		q^2 \\
\label{eq:sthgp-dx-sum-E-E(S)}
&=	\tau^2 \, \beta^2(n,k) \,
	\left \lbrack
		\alpha(n-k)
		\,+\, k \, (n-k) \, (n-k+1) \, 2^{n-k-2}
	\right \rbrack.
\end{flalign}
}

\begin{figure}[!ht]
 \begin{center}
  \begin{minipage}[t]{2.125in}
   \includegraphics[width=2.125in,clip=]{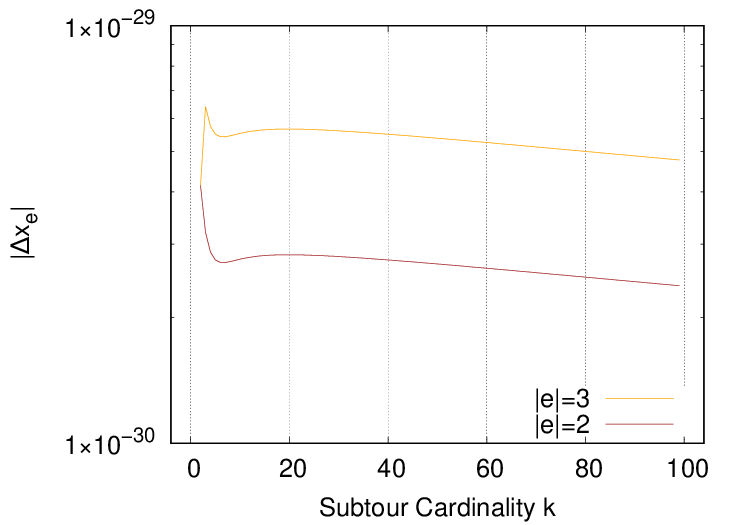}
   \caption{$\STHGP{n}$ subtour CD $|\Delta x_e|$ values for edge
	$e \in E$, $2 \le |e| \le 3$ within subtour of size $k$, $n = 100$}
   \label{fig:sthgp-dx-S-2-and-3}
  \end{minipage}
  \quad
  \begin{minipage}[t]{2.125in}
   \includegraphics[width=2.125in,clip=]{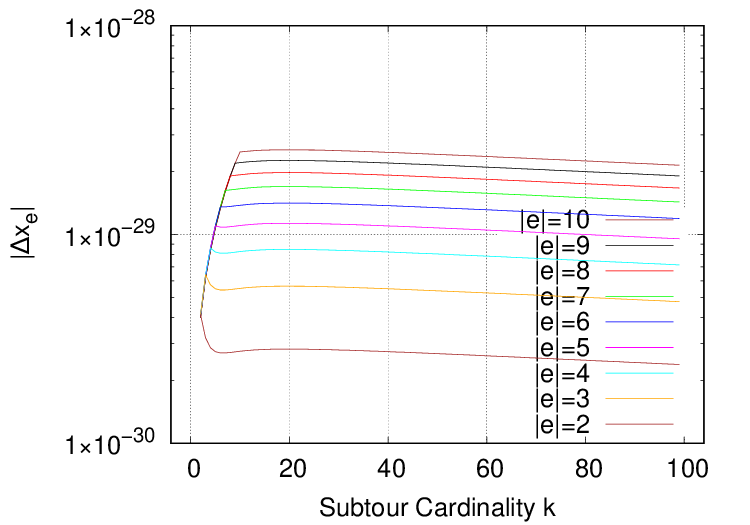}
   \caption{$\STHGP{n}$ subtour CD $|\Delta x_e|$ values for edge
	$e \in E$, $2 \le |e| \le 10$ within subtour of size $k$, $n = 100$}
   \label{fig:sthgp-dx-S-2-to-10}
  \end{minipage}
 \end{center}
\end{figure}

\begin{figure}[!ht]
 \begin{center}
  \begin{minipage}[t]{2.125in}
   \includegraphics[width=2.125in,clip=]{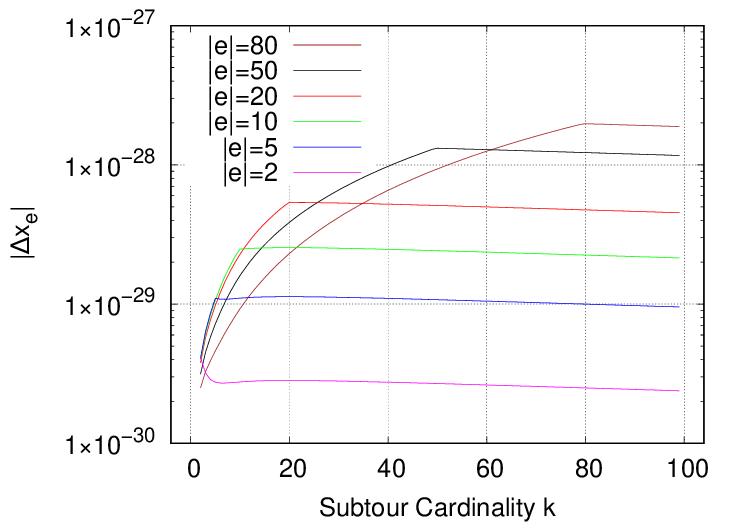}
   \caption{$\STHGP{n}$ subtour CD $|\Delta x_e|$ values for edge
	$e \in E$, $|e| \in \lbrace 2,\, 5,\, 10,\, 20,\, 50,\, 80 \rbrace$ within
	subtour of size $k$, $n = 100$}
   \label{fig:sthgp-dx-S-2-to-80}
  \end{minipage}
  \quad
  \begin{minipage}[t]{2.125in}
   \includegraphics[width=2.125in,clip=]{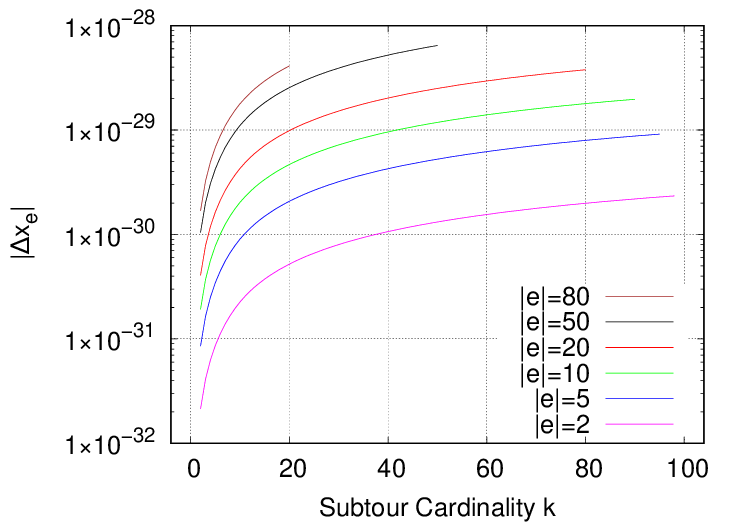}
   \caption{$\STHGP{n}$ subtour CD $|\Delta x_e|$ values for edge
	$e \in E$,
	$|e| \in \lbrace 2,\, 5,\, 10,\, 20,\, 50,\, 80 \rbrace$
	having $p=0$ vertices in common with the subtour, $n = 100$}
   \label{fig:sthgp-dx-0-2-to-80}
  \end{minipage}
 \end{center}
\end{figure}

Figure~\ref{fig:sthgp-dx-S-2-and-3} plots the $|\Delta x_e|$ values
as functions of subtour size $k$ for edges $e \in E(S)$, for $n=100$.
There is one plot for an edge of cardinality 2 and another for an edge
of cardinality 3.
(These are plots of $\Delta x_e(n,k,p,q)$, where
$n=100$, $2 \le k \le n-1$, $p = \min(k,|e|)$ and $q = |e| - p$
so that the subtour of $k$ vertices contains the maximum number of
vertices from $e$.)
Figure~\ref{fig:sthgp-dx-S-2-to-10}
plots similar curves for $2 \le |e| \le 10$.
Figure~\ref{fig:sthgp-dx-S-2-to-80}
plots similar curves for $|e| \in \lbrace 5, 10, 20, 50, 80 \rbrace$.
Note that although there is unusual behavior for small $k$, the
overall trend is that $|\Delta x|$ increases with $k$ until $e$ is
completely contained within the subtour, and decreases thereafter.

\begin{figure}[!ht]
 \begin{center}
  \begin{minipage}[t]{2.125in}
   \includegraphics[width=2.125in,clip=]{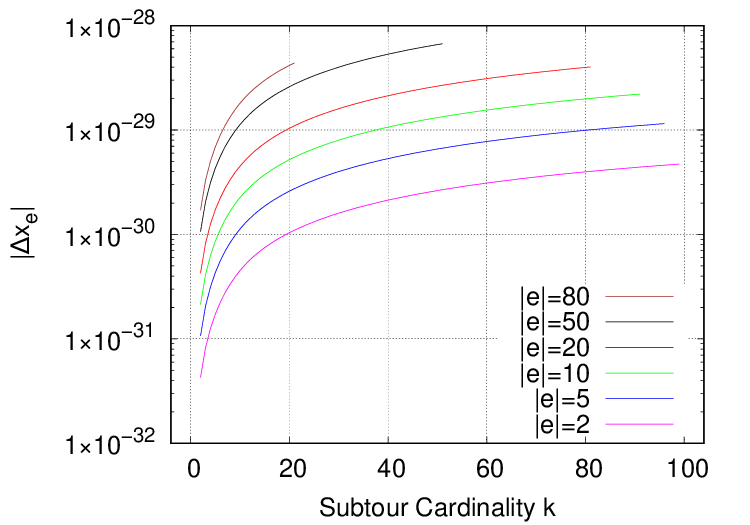}
   \caption{$\STHGP{n}$ subtour CD $|\Delta x_e|$ values for edge
	$e \in E$,
	$|e| \in \lbrace 2,\, 5,\, 10,\, 20,\, 50,\, 80 \rbrace$
	having $p=1$ vertex in common with the subtour, $n = 100$}
   \label{fig:sthgp-dx-1-2-to-80}
  \end{minipage}
  \quad
  \begin{minipage}[t]{2.125in}
   \includegraphics[width=2.125in,clip=]{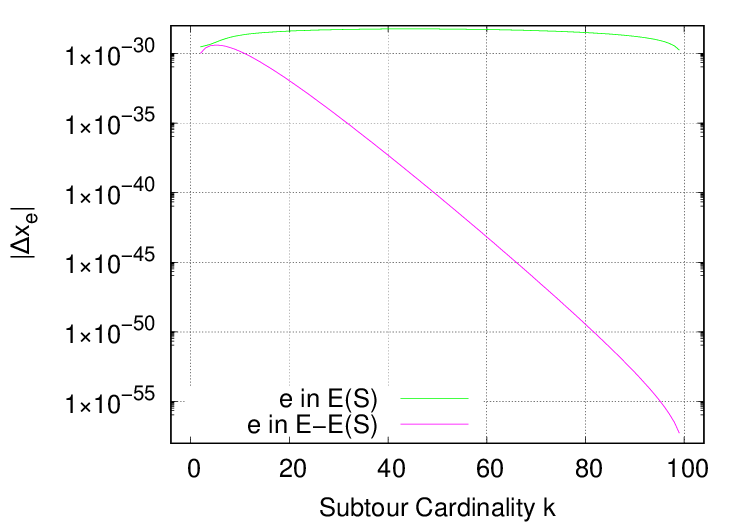}
   \caption{$\STHGP{n}$ subtour CD $\sum_e (\Delta x_e)^2$ sums:
	for edges $e \in E(S)$ and $E \in E-E(S)$, $n = 100$}
   \label{fig:sthgp-dx-sums}
  \end{minipage}
 \end{center}
\end{figure}

Figure~\ref{fig:sthgp-dx-0-2-to-80} plots the $|\Delta x_e|$
values for edges $e$ having $p=0$ vertices in common with the
subtour.
Figure~\ref{fig:sthgp-dx-1-2-to-80} plots the $|\Delta x_e|$
values for edges $e$ having $p=1$ vertex in common with the
subtour.
The trend for such edges $e \in E-E(S)$ is that $|\Delta x_e|$
increases with subtour cardinality $k$.

Figure~\ref{fig:sthgp-dx-sums} plots the sum of $(\Delta x_e)^2$ over
all $e \in E(S)$ and
all $e \in E-E(S)$.
The edges $e \in E(S)$ clearly dominate the summation of $(\Delta
x_e)^2$ over all $e \in E$.
To a first approximation (in both $\STGP{n}$ and $\STHGP{n}$),
the $|\Delta x_e|$ magnitudes generally decrease with $k$ for
$e \in E(S)$ and increase with $k$ for $e \in E-E(S)$.
As shown above for $n=100$, in $\STHGP{n}$ the magnitudes of the
individual $|\Delta x_e|$ vary by at most a factor of 830175, so their
squares vary by at most a factor of 689190530625.
The shape of the $\mathrm{CD}^2$ curve therefore appears to be
dominated by the sheer number of terms (exponential) of each
cardinality, governed by the bell curve shape of the binomial
coefficients.
(Note that ${100 \choose 50} = 100891344545564193334812497256$.)
In contrast, $\STGP{n}$ has only a polynomial (quadratic) number of
edges (4950 edges at $n = 100$), so that for small $k$, the large
$|\Delta x_e|$ values dominate the summation.

These numeric examples suggest that the surprising disagreement of
CD indicator regarding the subtours of STHGP vs STGP can be
attributed to dimensionality of the respective spaces.
With STHGP, the relatively small $\Delta x_e^2$ values appear to be
dominated by the sheer number of such terms being summed.
(Binomial coefficients give the number of subsets of each size, giving
CD curve for STHGP subtours more of a Bell curve shape.)
We need a more formal asymptotic argument regarding the STHGP subtours
to fully answer this question.

It is important to note that EPR indicator is completely
insensitive to the dimensionality of the space.
Because CD and EPR significantly disagree regarding the strength
of STGP subtours, both cannot be correct.
Since EPR gives a similar picture for both STHGP and STGP subtours
while CD gives a conflicting picture, and because the conflict appears
to be driven by the dimensionality of the space (exponential versus
polynomial), we conclude that CD is likely to be a less reliable
indicator of strength than EPR.

\FloatBarrier

\section{Validation}
\label{sec:validation}

We independently validated our indicator formulae for small $n$
using software available under an open source license
~\cite{WarmeIndicatorsBundle},~\cite{WarmeCALIB}.
(This software also produces most plots and figures herein.)
The first validation step is to recursively enumerate all extreme
points $X$ for the prospective polytope $P$
(e.g., tours for $\TSP{n}$ and spanning trees for $\STHGP{n}$ and
$\STGP{n}$).
We validated $3 \le n \le 12$ for TSP,
$3 \le n \le 9$ for STHGP,
and $3 \le n \le 9$ for STGP and their respective inequalities.

For EPRs, choose a representative inequality from
the class to validate and simply count how many members of $X$
satisfy the contraint with equality.

We validated the coordinates of the centroid directly via
equation~(\ref{eq:centroid-def}).

CD squared formulae were validated by solving quadratic programming
formulation~(\ref{eq:qp-objective})--(\ref{eq:convex-bounds})
using CPLEX 12.5.1.
(Note that for $n \ge 6$ there are ${{n+1} \choose {n-6}}$ valid TSP
3-toothed comb inequality configurations, and we validated every such
configuration.)
Numeric results from CPLEX for CD squared usually matched the
exact solutions to at least 7 digits of precision.
(The few numerical outliers exceeding this threshold were for STHGP
non-negativity inequalities with $n=k$ for $6 \le n \le 9$.  It should
be possible to significantly improve numerical accuracy by either
reducing the CPLEX optimality tolerance, or by pre-loading a CPLEX
basis corresponding to the exact closest point coordinates.)
As expected, the normal and weak CDs were the same for non-negativity
and subtour contraints of the TSP, subtours of STGP, and different for
all other inequalities studied.

To verify the formula for STHGP subtour angles it suffices to
verify the definition of $f(n,p,q,r)$ as given in
Equation~(\ref{eq:f(n,p,q,r)}).
(We note that both the numerator and denominator squared
of~(\ref{eq:f(n,p,q,r)}) are integers.)
We do this using
Equation~(\ref{eq:f(n,p,q,r)-dot-prod-equivalence}).
Given $n$, $p$, $q$ and $r$, our validation software constructs
explicit
$S_1 = \lbrace 1,\, \ldots,\, p+r \rbrace$ and
$S_2 = \lbrace p+1,\, \ldots,\, p+q+r \rbrace$
having these properties, and the corresponding subtour hyperplane
equations $a_1.x = b_1$ and $a_2.x = b_2$.
Using~(\ref{eq:angle-project-a1})
and~(\ref{eq:angle-project-a2})
it computes projections $\hat{a}_1$ and $\hat{a}_2$ onto the affine
hull.
After verifying that $\hat{a}_1.c = 0$ and  $\hat{a}_2.c = 0$,
it uses the simple dot product calculations
in~(\ref{eq:f(n,p,q,r)-dot-prod-equivalence}) to compute
$f(n,p,q,r)$ from ``first principles.''
The numerator and denominator squared of this can be compared directly
to those obtained with~(\ref{eq:f(n,p,q,r)}).
We checked all valid combinations of $(n,p,q,r)$
for $3 \le n \le 15$:
the closed-form angular formula and ``first principle'' computations
yielded identical results in all cases.

Finally, for all $3 \le n \le 100$ and $k \le 2 < n$,
we verified that
adding~(\ref{eq:sthgp-dx-sum-E(S)})
and~(\ref{eq:sthgp-dx-sum-E-E(S)})
yields the same rational number as the square
of~(\ref{eq:sthgp-subtour-cd-proof-final-form}).

\section{Indicators Not Infallible}
\label{sec:counter-examples}

In this section we present concrete polytopes demonstrating that
neither of the indicators studied in this paper are infallible.
We demonstrate this by showing each indicator to be inconsistent with
the solid angle indicator, wherein strength is expressed as the ratio
$a_f / S$, where $a_f$ is the solid angle subtended (from the
centroid) by facet $f$, and $S$ is the solid angle of the entire sphere.

\subsection{Example Polytope Where EPR Performs Poorly}
\label{sec:EPR-counter-example}

We present a polytope in $\R^3$ for which the facet having the largest
EPR is in fact the weakest of all facets
according to the solid angle indicator.
This polytope is inspired by, and similar in structure to a
 ``brilliant cut'' diamond.
The polytope is illustrated in
Figures~\ref{fig:brilliant-cut-side-view} (side view)
and~\ref{fig:brilliant-cut-top-view} (top view), which utilize the
formal geometric terminology for brilliant cut diamonds.
\FloatBarrier
\MyBeginFig
 \begin{center}
  \begin{minipage}[t]{2.25in}
   \includegraphics[width=2.25in,clip=]{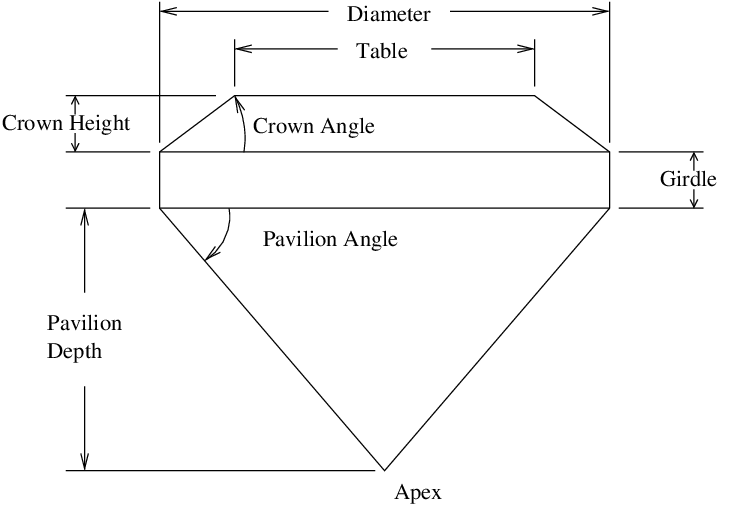}
   \caption{Example polytope where EPR performs poorly (side view)}%
   \label{fig:brilliant-cut-side-view}
  \end{minipage}
  \quad
  \begin{minipage}[t]{2.25in}
   \includegraphics[width=2.25in,clip=]{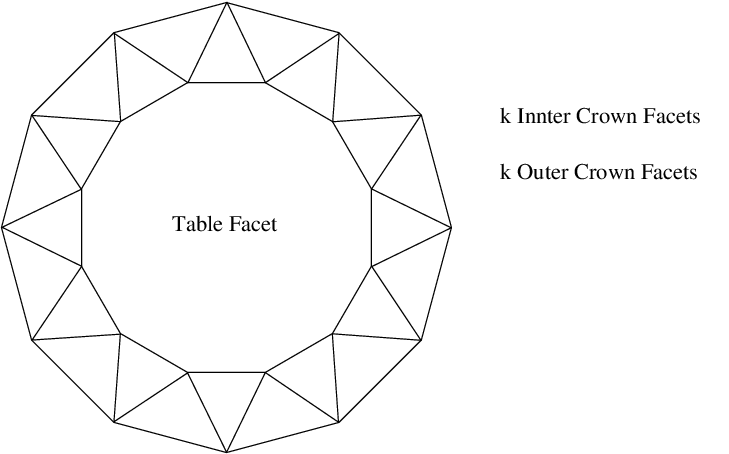}
   \caption{Example polytope where EPR performs poorly (top view)}%
   \label{fig:brilliant-cut-top-view}
  \end{minipage}
 \end{center}
\MyEndFig
\FloatBarrier
\noindent
The table facet of this polytope is a regular $k$-vertex polygon,
where $k$ is large.
The crown consists of $k$ inner facets and $k$ outer facets, each
triangular.
The girdle consists of $k$ {\em upper} facets (each a triangle
sharing an edge with the crown), and $k$ {\em lower} facets (each a
triangle sharing an edge with the pavilion).
The pavilion consists of $k$ triangular facets that meet at the
apex.
We require that the crown angle be some very small angle
$\epsilon > 0$.
The table facet is incident to $k$ extreme points, where $k$ is large.
All other facets have 3 extreme points.
The table facet therefore has the largest EPR, thereby
predicting the table facet to be the strongest of all facets.
In reality, the table facet is the weakest of all facets, because the
cone of objective values for which one or more of the table facet's
extreme points is optimal is a very skinny cone --- the solid angle
subtended by this cone is a vanishingly small fraction of the solid
angle of the entire 3-d sphere.
The crown, girdle and pavilion facets are much stronger, each being
optimal over a cone having solid angle that is a significant fraction
of the entire 3-d sphere.
This polytope is easily extended to higher dimensions.

The centroid $C$ lies along the axis through the apex, and can be made
to have an arbitrarily large distance from the table facet by moving
the apex away from the table facet.
Under these conditions, the pavilion facets are closest to $C$,
the table facet is far from $C$, and CD
correctly labels the table facet to be weak.

\subsection{Example Polytope Where CD Performs Poorly}
\label{sec:CD-counter-example}

We present a polytope in $\R^d$ for which the facets having the
smallest CD are in fact the weakest of all facets according to the
solid angle indicator.
See Figure~\ref{fig:cdist-counter-example} for an illustration of this
polytope in $d=3$ dimensions.
Let $P(d)$ be the symmetric $d$-dimensional simplex.
We construct $P'(d)$ from $P(d)$ as follows:
truncate a tiny section of the polytope at each of its $d+1$
vertices, replacing each of them with a copy of a small, yet very
complex polytope $V$.  (Each copy of $V$ must be suitably translated
and rotated into proper position and alignment.)
Thus $P'(d)$ is just $P(d)$ with its corners rounded off (but in a
flat-faced, polyhedral fashion).
\begin{center}
 \begin{minipage}[t]{2.75in}
  \includegraphics[width=2.75in,clip=]{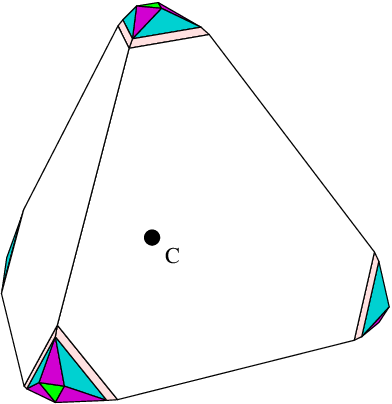}
  \captionof{figure}{Example polytope where CD performs poorly}%
  \label{fig:cdist-counter-example}
\end{minipage}
\end{center}
The portion truncated from each vertex is isomorphic to the original
$P(d)$ (scaled by some small $\epsilon > 0$).
Each of the $d+1$ original facets $F$ are copies of $P(d-1)$, and
after truncation, each of $F$'s $d$ vertices is replaced by $d-1$
vertices within $P'(d)$.
Each original facet $F$ of $P(d)$ becomes incident to $d(d-1)$
vertices in $P'(d)$.
(This is a {\em very small} number of extreme points, so far as things
usually go in the realm of high-dimensional polyhedra.)
We can arrange the following attributes of each grafted polytope
$V$:
\begin{enumerate}
 \item They are all congruent and symmetric so that the centroid of
   $P'(d)$ is the same as the centroid of $P(d)$.
 \item Let $F$ be a facet of $P(d)$, $F'$ be the corresponding facet
   of $P'(d)$, and $H'$ be a hyperplane defining some facet $Q$ of $V$
   such that $\mathrm{dim}(F' \cap Q) \,=\, d-2$, i.e., $F'$ and $Q$
   intersect to form a ridge of $P'(d)$.
   Let $n(H')$ and $n(F')$ be the unit normal vectors defining $H'$ and
   $F'$, respectively.
   We can arrange for the angle between $n(H')$ and $n(F')$ to be some
   very small $\epsilon > 0$.
\end{enumerate}
Let $\F$ be the set of facets of $P(d)$
and $\F'$ be the corresponding set of facets of $P'(d)$.
$P'(d)$ now has the following attributes:
\begin{enumerate}
 \item Facets $\F'$ have the smallest CD among all facets of $P'(d)$.
   (All of the other facets of $P'(d)$ reside out near where the
   vertices of $P(d)$ used to be.)
 \item There is a cone $D$ of direction vectors around $n(F')$, whose
   extreme rays are the $n(H')$.
   Let $c \in \mathrm{interior}(D)$.
   Then using $c$ as an objective function of an LP over $P'(d)$
   results in an optimal solution $x$ that is an extreme point (or
   other higher dimensional face) of $F'$.
   (If $c = n(F')$, then all $x \in F'$ are optimal solutions.)
   By construction, the solid angle subtended by $D$ is a vanishingly
   small fraction of the complete sphere's solid angle.
 \item This implies that if we choose a unit objective vector $c$
   uniformly at random, then the probability that the optimal LP
   solution $x$ is incident to $F'$ is vanishingly small.
   (The solid angle indicator labels $F'$ as being very weak.)
\end{enumerate}
Although CD labels $F'$ as the strongest
of all the facets of $P'(d)$, we see that in fact the facets $F' \in
\F'$ are very weak according to the solid angle indicator.
The two indicators disagree in their appraisal of the strength of
$F'$ with respect to polytope $P'(d)$:
CD incorrectly labels $F'$ to be strong, while
EPR correctly labels $F'$ to be weak ($F'$ has
very few incident extreme points, $d(d-1)$ of them.)
Polytope $P'(d)$ is an example where there appears to be a lack of the
desired (negative) correlation between the two indicators.
(CDs and EPRs are both numerically small for facets $F'$.)

It is an interesting open problem to construct a single polyhedron
upon which both EPR and CD indicators perform poorly.
\FloatBarrier

\section{Future Research}
\label{sec:open-problems}

We now mention some areas for future research and open problems.

\subsection{Other Polyhedra and Inequalities}
\label{sec:other-polytopes-and-inequalities}

This field is completely open.
EPR and/or CD indicators for essentially all other
polyhedron/facet pairs are currently unknown.
Each new indicator result could reveal previously unknown properties
of the underlying problem, some of which can be exploited in
seperation, strengthening or other algorithms.
The same can be said regarding calculated angles between various pairs
of inequalities.

Other indicators similar to EPR and CD should also be studied more
closely.
We should definitely heed Galileo's words:
\begin{quotation}
``Measure what is measurable, and make measurable what is not so.''
\end{quotation}

\subsection{Ranking Distinct Inequality Classes}
\label{sec:ranking}

Given some polyhedron $P$, and two facet classes $A$ and $B$,
EPR and CD indicators should be just as applicable for comparing
the relative strength of two inequalities from $A$ as for comparing
the relative strength of one inequality from $A$ versus another
inequality from $B$.
Consider $P = \TSP{n}$, $A =$ subtours, and $B =$ 3-toothed combs.
For fixed $n$ EPR (or CD) indicator could be used to obtain a
complete strength ranking of $A \cup B$.
There are many subtours that are stronger than most 3-toothed combs,
and there are some 3-toothed combs that are stronger than certain
other subtours.
Knowing where these boundaries lie could be useful.

\subsection{Characterizing Weakest Inequalities}
\label{sec:weakest}

For the subtours of $\STHGP{n}$, obtain closed-forms $f(n)$ such that
$k = f(n)$ gives the cardinality $k$ of the weakest subtour for
$n$-vertex hypergraphs.
There would be one $f(n)$ for EPR and another for CD.
Such characterization of the weakest subtour would, for example, allow
one to choose whether to strengthen by reduction or augmentation.

For each of these, compute
\begin{displaymath}
	\lim_{n \rightarrow \infty} {{f(n)} \over n}.
\end{displaymath}

\subsection{Bounds on Disagreement}
\label{sec:disagreement-bounds}

Well-known bounds exist on the ``disagreement'' between various vector
norms:
\begin{align*}
{\hskip 1.0in}
&	||x||_2 &&\le ||x||_1 &&\le \sqrt{n}\,||x||_2,
	{\hskip 2.5in} \\
&	||x||_{\infty} &&\le ||x||_2 &&\le \sqrt{n}\,||x||_{\infty}, \\
&	||x||_{\infty} &&\le ||x||_1 &&\le n\,||x||_{\infty}.
\end{align*}
One could certainly ask what the maximum disagreement is between two
indicators $I_1$ and $I_2$ (e.g., EPR and CD), especially over certain
classes of polyhedra and/or inequalities.
Since indicators can have the opposite sense (like EPR and CD) and can
differ markedly in scale, this is probably best handled in a
``comparative'' sense as we do in
Section~\ref{sec:compare-indicators} --- characterizing the extent to
which indicators $I_1$ and $I_2$ can rank specific inequalties $H_1$
and $H_2$ differently.
A good place to begin would be zero-one polytopes.

\section{Conclusion}
\label{sec:conclusion}

We have intensively studied the notion of quantitative indicators for
the strength of inequalities (relative to an underlying polyhedron),
and defined several concrete indicators.
We have applied two of these indicators (EPR and CD) to several
well-known facet classes of the Traveling Salesman, Spanning Tree in
Hypergraph and Spanning Tree in Graph polytopes.
With the exception of the surprising CD results for STGP subtours,
the indicators agree strongly with each other, and
corroborate the notions of strength observed empirically during
computations.

Although not infallible, we assert that these indicators can usefully
be employed as an a priori guide to the creation of highly tuned,
{\em effective} separation algorithms and constraint strengthening
procedures.
Knowing exactly where the {\em strong} constraints reside eliminates
the guesswork and experimentation that is often required to achieve
superior performance from these algorithms.
The second paper presents strong computational evidence supporting
these claims~\cite{WarmeArxivIndicators2}.

\begin{credits}
\subsubsection{\ackname}
I am grateful to Martin Zachariasen, Pawel Winter, Jeffrey Salowe,
Karla Hoffman, Maurice Queyranne, Warren D. Smith and Bill Cook:
this work would not exist without their helpful discussions and
encouragement.
I also thank Group W for their support: this work has no direct
connection to their business.
Proverbs 25:2, Soli Deo Gloria.

\subsubsection{\discintname}
The author has financial and intellectual property interests in
GeoSteiner.
\end{credits}

\bibliographystyle{splncs04}
\bibliography{phd,steiner,local}

% This prevents the final page of the bibliography from being stretched
% vertically with pronounced vertical gaps between each entry.
\vskip 0.1in plus 1filll

\eject

\strut
\vskip 0.25in
\noindent
{\huge \bf Revision History}

\vskip 0.25in

\noindent
03/21/2024: Original version.

\vskip 0.25in

\noindent
07/08/2025: Version 2:
\begin{itemize}
  \item Improvements from review of condensed journal version.
  \item Cite additional work.
  \item Improvements to formatting and sectioning.
  \item Retract incorrect EPR for TSP 3-toothed combs.
  \item Add alternate EPR for STHGP subtours.
  \item Add Computational Highlights section.
\end{itemize}

\end{document}